\documentclass[aos]{imsart}
\usepackage{amsthm,amsmath,natbib, verbatim}
\usepackage{amssymb,graphicx,pstricks}
\bibpunct{[}{]}{,}{n}{ }{,}
\startlocaldefs
\def\nlab#1{#1^{(n)}}
\def\cM{\mathcal{M}}

\newcommand\sgn{\mathop{\rm sgn}}
\newcommand\R{{\mathbb{\mathchar"52}}} 
\def\argmin{\mathop{\rm argmin}}

\newtheorem{theorem}{Theorem}

\newtheorem{lemma}{Lemma}

\let\hat\widehat
\let\tilde\widetilde
\let\sss\scriptscriptstyle

\def\Norm#1{\left\Vert#1\right\Vert}
\def\norm#1{\Vert#1\Vert}

\def\qed{\hskip1pt $\;\;\scriptstyle\Box$}

\def\todo#1{
  \ifhmode%
  \unskip%
  {\dimen1=\baselineskip \divide\dimen1 by 2 }\fi%
  {\marginpar{\renewcommand{\baselinestretch}{0.8}%
    \footnotesize\raggedright #1}}}
\endlocaldefs
\newtheorem{cor}{Corollary}[section]

\newtheorem{definition}{Definition}[section]
\newcommand{\bed}{\begin{definition}}
\newcommand{\eed}{\end{definition}}
\newcommand{\beq}{\begin{equation}}
\newcommand{\eeq}{\end{equation}}

\newcommand{\eps}{\epsilon}

\newcommand{\bitem}{\begin{itemize}}
\newcommand{\eitem}{\end{itemize}}

\newcommand{\goto}{\rightarrow}

\newcommand{\beqn}{\begin{equation}}
\newcommand{\eeqn}{\end{equation}}
\newcommand{\balign}{\begin{align}}
\newcommand{\ealign}{\end{align}}

\newcommand{\bphi}{\bar{\Phi}}

\def\limsup{\mathop{\overline{\rm lim}}}
\def\liminf{\mathop{\underline{\rm lim}}}
\usepackage{amssymb}
\begin{document}

\begin{frontmatter}

\title{Revisiting Marginal Regression}
\runtitle{Marginal Regression}



\author{\fnms{Christopher R.} \snm{Genovese}\ead[label=e1]{chris@stat.cmu.edu}}
\address{5000 Forbes Ave\\ Pittsburgh, PA 15213\\ \printead{e1}}
\author{\fnms{Jiashun} \snm{Jin}\ead[label=e2]{jiashun@stat.cmu.edu}}
\address{5000 Forbes Ave\\ Pittsburgh, PA 15213\\ \printead{e2}}
\and
\author{\fnms{Larry} \snm{Wasserman}\corref{}\ead[label=e3]{larry@stat.cmu.edu}}
\address{5000 Forbes Ave\\ Pittsburgh, PA 15213\\ \printead{e3}}
\affiliation{Carnegie Mellon University}
\runauthor{Genovese et al}
\today

\begin{abstract}
The lasso has become an important practical tool for high dimensional regression
as well as the object of intense theoretical investigation.
But despite the availability of efficient algorithms, the lasso remains computationally
demanding in regression problems where the number of variables
vastly exceeds the number of data points.
A much older method, marginal regression, largely displaced by the lasso,
offers a promising alternative in this case.
Computation for marginal regression is practical even when the dimension is very high.
In this paper, we study the relative performance of the lasso and marginal regression for regression problems
in three different regimes: (a)~exact reconstruction in the noise-free and noisy
cases when design and coefficients are fixed, (b)~exact reconstruction in the noise-free
case when the design is fixed but the coefficients are random,
and (c)~reconstruction in the noisy case where performance is measured by the number
of coefficients whose sign is incorrect.

In the first regime, we compare the conditions for exact reconstruction of the two procedures, 
find examples where each procedure succeeds while the other fails, and characterize the advantages and disadvantages of each.
In the second regime, we derive   conditions under which marginal regression  will provide exact reconstruction with high probability.
And in the third regime, we derive rates of convergence for the procedures
and offer a new partitioning of the ``phase diagram,'' 
that shows when exact or Hamming reconstruction is effective.

In addition to theoretical investigation, 
we present simulations showing that in practice,
marginal regression and the lasso can have comparable performance,
while the computational advantages of marginal regression make
it feasible for much larger problems.

\end{abstract}



\end{frontmatter}

\newpage

\section{Introduction}

A central theme in recent work on regression
is that sparsity plays a critical role in
effective high-dimensional inference.
Consider a regression model,
\begin{equation} \label{eq::basic-regression}
Y = X \beta + z,
\end{equation}
with 
response 
$Y = (Y_1,\ldots, Y_n)^T$,
$n\times p$ design matrix $X$,
coefficients $\beta = (\beta_1,\ldots,\beta_p)^T$,
and noise variables
$z = (z_1,\ldots,z_n)^T$.
Loosely speaking, this model is high-dimensional when $p \gg n$
and is sparse when many components of $\beta$ equal zero.

An important problem in this context is variable selection:
determining which components of $\beta$ are non-zero.
For general $\beta$, the problem is underdetermined,
but recent results have demonstrated that under particular conditions on $X$, to be discussed below,
sufficient sparsity of $\beta$ allows 
(i)~exact reconstruction of $\beta$ in the noise-free case \cite{Tropp:2004}
and
(ii)~consistent selection of the non-zero coefficients in the noisy-case
\cite{Chen-Donoho:1998,  CWX, Candes,  Donoho:06a, Donoho-Elad:2003,  Fan-Lv:2008,  Fuchs:2005,  Knight-Fu:2000, Meinshausen-Buhlmann:2006,  Tropp:2004, Wainwright:2006, Zhao-Yu:2006, Zou}.   Many of these results are based on showing that under sparsity constraints,
a convex optimization problem that controls the $\ell^1$ norm of the coefficients
has the same solution as an (intractable) combinatorial optimization problem that
controls the number of non-zero coefficients.

In practice, the lasso \cite{Tibshirani:1996, Chen-Donoho:1998} has become one of the main tools for sparse high-dimensional variable selection,
due both to its computational simplicity and its direct connection to these theoretical results.
The lasso estimator in the regression problem is defined by
\begin{equation}  \label{Definelasso} 
\hat\beta_{\sss\text{lasso}} = \argmin_\beta \Norm{Y- X\beta}_2^2 + \lambda \norm{\beta}_1,
\end{equation}
where $\norm{\beta}_1 = \sum_j |\beta_j|$ and $\lambda \geq 0$ is a regularization parameter
that must be specified.
The lasso gives rise to a convex optimization problem and thus is computationally tractable
even for moderately large problems.
Indeed, the LARS algorithm \cite{Efron-etal:2004} can compute the entire solution path as a function of
$\lambda$ in $O(p^3 + n p^2)$ operations.
Gradient descent algorithms for the lasso are  faster  in practice, but have the same computational complexity. 
For very large $p$,  the lasso remains computationally demanding.

A much older and computationally simpler method for variable selection
is marginal regression (also called correlation learning, simple thresholding \cite{Donoho:06a}, and sure screening \cite{Fan-Lv:2008}),  
in which the outcome variable is regressed on each covariate separately.  
To compute the marginal regression estimates for variable selection,
we begin by computing the marginal regression coefficients
which, assuming $X$ has been standardized, are
\begin{equation}
\hat\alpha \equiv X^T Y.
\end{equation}
Then, we threshold $\hat\alpha$ using the tuning parameter $t  > 0$:
\begin{equation} \label{Definemarginal} 
\hat\beta_j = \hat\alpha_j 1\{|\hat\alpha_j|  \ge t\}.
\end{equation}
This requires $O(n p)$ operations, two orders faster than the lasso for $p \gg n$,
and so is tractable for much larger problems.

The lasso has mostly displaced marginal regression in practice.
But the computational advantage for large problems prompts a second look.
Tibshirani and Witten \citet{Tibshirani-Witten:2009} have found that marginal regression sometimes
outperforms the lasso in predictive error.
Here we revisit marginal regression as a tool for variable selection
and ask whether there is any strong reason to prefer the lasso.
If marginal regression exhibits comparable performance, theoretically and empirically,
then it offers a plausible alternative to the lasso.
Put another way: because of its simplicity, marginal regression only needs to tie to win. 

In this paper, we study the relative performance of the lasso and marginal regression
in three different regimes.
In Section \ref{sec:main}, we compare the conditions that guarantee exact variable selection
in the noise-free case and briefly compare the conditions for consistent variable selection (i.e. sparsistency)
in the noisy case.
The two sets of conditions are generally overlapping, and we give examples where each procedure
fails while the other succeeds.
One advantage of the lasso is that, given a fixed matrix $X$, the conditions for its success
hold over a larger class of $\beta$'s than that of marginal regression.  On the other hand,
marginal regression has a larger tolerance for collinearity than does the lasso
and is somewhat easier to tune, as we illustrate in Section \ref{subsec:noise}. 
 
In Section \ref{sec:exten}, we consider the regime where the design matrix $X$ is fixed but
the coefficient vector $\beta$ is randomly generated.
We find conditions such that marginal regression performs well with overwhelming probability.
The main condition, which we call {\it faithfulness},  is  closely related to both the Faithfulness Condition  of \cite{Meinshausen-Buhlmann:2006}
and the Incoherence Condition  of \cite{Donoho-Elad:2003}.  
The Incoherence Condition depends only on $X$ and is thus checkable in practice, but it aims to control  the worst case
so is quite conservative.
The Faithfulness  Condition  of \cite{Meinshausen-Buhlmann:2006}
is relatively less stringent but depends on the unknown support of the parameter vector.
Our version of the Faithfulness Condition  strikes a compromise between the two.

Although exact variable selection has been the focus of many studies in the literature,
it is rare in practice to select exactly the right variables,
so it is natural to measure performance in terms of the deviation from exact selection.
In Section \ref{sec:hamm}, we study the convergence rates of the two procedures
in Hamming distance between $\sgn({\beta})$ and $\sgn({\hat\beta})$.
Our main result in this section is a new partition of the parameter space into three
regions I--III.
In the interior of region I, exact variable selection is possible (asymptotically),
and both procedures achieve this given properly chosen tuning parameters.  
In region II, it is possible to have a variable selection  procedure that recovers most  relevant variables, but not all of them. 
And in Region III,  successful variable selection is impossible,
and the optimal Hamming distance is asymptotically equivalent to the total number of relevant variables. 
 
Finally, in Section \ref{sec:simul}  
we present simulation studies showing that marginal regression
and the lasso perform comparably over a range of parameters.    
Section \ref{sec:proof} gives the proofs of all theorems and 
lemmas in the order they appear.


{\em Notation.}
For a real number $x$, let $\sgn(x)$ be -1, 0, or 1 when $x < 0$, $x = 0$, and $x > 0$;
and for a vector $u\in\R^k$, define $\sgn(u) = (\sgn(u_1),\ldots,\sgn(u_k))^T$.
We will use $\norm{\cdot}$, with various subscripts, to denote vector and matrix norms,
and $|\cdot|$ to represent absolute value, applied component-wise when applied to vectors.
With some abuse of notation, we will write $\min u$ ($\min |u|$) to denote the minimum (absolute)
component of a vector $u$.
Inequalities between vectors are to be understood component-wise as well.

\section{Noise-Free Conditions for Exact Variable Selection} \label{sec:main} 
Consider a sequence of regression problems
with deterministic design matrices, indexed by sample size $n$,
\begin{equation} \label{eq::regression-seq}
\nlab{Y} = \nlab{X} \nlab{\beta} + \nlab{z}. 
\end{equation}
Here,   $\nlab{Y}$ and $\nlab{z}$ are $n\times 1$ response and noise vectors,
respectively,
$\nlab{X}$ is an $n \times \nlab{p}$ matrix and $\nlab{\beta}$ is a
$\nlab{p} \times 1$ vector, where we typically assume $\nlab{p} \gg n$.
We assume that $\nlab{\beta}$ is sparse in the sense that it has $\nlab{s}$ nonzero components
where $\nlab{s} \ll \nlab{p}$.   
By rearranging $\nlab{\beta}$ without loss of generality,
we can partition each $\nlab{X}$ and $\nlab{\beta}$ into ``signal'' and ``noise'' pieces, 
corresponding to the non-zero or zero coefficients, as follows:
\begin{equation}
\nlab{X} = \left( \nlab{X}_S, \nlab{X}_N \right) \qquad \nlab{\beta} = \left(
\begin{array}{c}
\beta_S \\
\beta_N
\end{array} \right). 
\end{equation}
In fact, we assume that $\nlab{\beta}_S \in \cM^{\nlab{s}}_{\nlab{\rho}}$
for a sequence $\nlab{\rho} > 0$ (and not converging to zero too quickly) with
\begin{equation} \label{eq::Mspace}
\cM^k_a = \left\{ x = (x_1, \ldots, x_k)^T  \in\R^k:\;  \mbox{$ |x_j| \ge a$   for all $1 \leq j \leq k$}\right\},
\end{equation}
for positive integer $k$ and $a > 0$.
This commonly used condition on $\nlab{\beta}$ ensures that the non-zero components are
not too close to zero to be indistinguishable.
Finally, define the Gram matrix $\nlab{C} =   ({\nlab{X}})^T\nlab{X}$ and partition this as
\begin{equation}
\nlab{C} = 
\left(\begin{array}{cc}
\strut\nlab{C}_{SS} & \nlab{C}_{SN} \\\noalign{\smallskip}
\strut\nlab{C}_{NS} & \nlab{C}_{NN}
\end{array} 
\right),
\end{equation}
where of course $\nlab{C}_{NS} = (\nlab{C}_{SN})^{T}$.   Except in Sections \ref{sec:hamm}--\ref{sec:simul},  we suppose  $X^{(n)}$ is normalized so that all diagonal coordinates of $C^{(n)}$ are $1$.  

These $\nlab{}$ superscripts become tedious, so for the remainder of the
paper, we suppress them unless necessary to show variation in $n$.
The quantities $X$, $C$, $p$, $s$, $\rho$, as well as the tuning
parameters $\lambda$ (for the lasso; see (\ref{Definelasso})) and $t$ (for marginal regression; see  (\ref{Definemarginal})) 
are all thus implicitly dependent on $n$. 
We use $\cM_\rho$ to denote the space $\cM^{\nlab{s}}_{\nlab{\rho}}$.

We will begin by specifying conditions on $C$, $\rho$, $\lambda$, and $t$
such that in the noise-free case, exact reconstruction of $\beta$ is possible
for the lasso or marginal regression, for all $\beta_S \in \cM_\rho$.
These in turn lead to conditions on $\nlab{C}$, $\nlab{p}$, $\nlab{s}$, $\nlab{\rho}$, $\nlab{\lambda}$, and $\nlab{t}$
such that in the case of homoscedastic Gaussian noise, the non-zero coefficients
can be selected consistently, meaning that for all sequences $\nlab{\beta}_S\in \cM^{\nlab{s}}_{\nlab{\rho}} \equiv \cM_\rho$,
\begin{equation} \label{eq::sparistency}
P\biggl(  \left|\sgn(\nlab{\hat\beta})\right| = \left|\sgn(\nlab\beta)\right| \biggr)  \to 1,
\end{equation}
as $n\to\infty$. (This property was dubbed \emph{sparsistency} by Pradeep Ravikumar \cite{Ravi}.)
Our goal is to compare these conditions.
In this section, we focus on the noise-free case,  and keep the discussion on the noise case brief.   


\subsection{Exact reconstruction conditions for the lasso in the noise-free case} 

We begin by considering three conditions in the noise-free case
that are now standard in the literature on the lasso:

\begin{itemize}
\item[] {\it
{\bf Condition E}.~The minimum eigenvalue of $C_{SS}$ is positive. 
}
\medskip
\item[] {\it
{\bf Condition I}.~(Irrepresentableness)
$$
\max \left| C_{NS} C_{SS}^{-1} \,\sgn(\beta_S) \right| \le 1.
$$
}
\item[]
{\it
{\bf Condition J}.~
$$
\min \left| \beta_S \,-\, \lambda C_{SS}^{-1}\, \sgn(\beta_S) \right| > 0.
$$
}
\end{itemize}
Because $C_{SS}$ is symmetric and non-negative definite, Condition E is equivalent
to $C_{SS}$ being invertible. Later we will strengthen this condition.
A critical feature of Condition I is that it only depends on the sign pattern,
as we will see.

For the noise-free case,  Wainwright \cite[Lemma 1]{Wainwright:2006} shows
that assuming Condition E, conditions I and J are necessary and sufficient
for the existence of a lasso solution $\hat\beta$ with tuning parameter $\lambda$
such that 
$$
\sgn(\hat\beta) = \sgn(\beta).
$$
(See also \cite{Zhao-Yu:2006}). 
Note that this result is stronger than correctly selecting the non-zero coefficients,
as it gets the signs correct as well.

Maximizing the left-hand side of Condition I considers all $2^s$ sign patterns
and gives $\norm{C_{NS} C_{SS}^{-1}}_\infty$,
the maximum-absolute-row-sum matrix norm.
It follows that Condition I holds for all $\beta_S\in\cM_\rho$ if and only if
$\norm{C_{NS} C_{SS}^{-1}}_\infty \le 1$.
Similarly, one way to ensure that Condition J holds over $\cM_\rho$
is to require that every component of $\lambda C_{SS}^{-1}\, \sgn(\beta_S)$
be less than $\rho$.
The maximum component of this vector over $\cM_\rho$ equals $\lambda \norm{C_{SS}^{-1}}_\infty$,
which must be less than $\rho$.
A simpler relation, in terms of the smallest  eigenvalue  of $C_{SS}$ is 
\begin{equation}
\frac{\sqrt{s}}{\text{eigen}_{\sss\text{min}}(C_{SS})} = 
\sqrt{s} \norm{C_{SS}^{-1}}_2 \ge \norm{C_{SS}^{-1}}_\infty \ge 
\norm{C_{SS}^{-1}}_2 = \frac{1}{\text{eigen}_{\sss\text{min}}(C_{SS})},
\end{equation}
where the inequality follows from the symmetry of $C_{SS}$ and standard norm inequalities.

Stronger versions of the above conditions will be useful.
\begin{itemize}
\item[] {\it
{\bf Condition E'}.~The minimum eigenvalue of $C_{SS}$ is no less than $\lambda_0 > 0$,
where $\lambda_0$ does not depend on $n$.
}
\medskip
\item[] {\it
{\bf Condition I'}.~
$$
 \norm{C_{NS} C_{SS}^{-1}}_\infty \le 1 - \eta,
$$
for $0 < \eta < 1$ small and independent of $n$.
}
\medskip
\item[]
{\it
{\bf Condition J'}.~
$$
 \lambda < \frac{\rho}{\norm{C_{SS}^{-1}}_\infty}.
$$
}
\end{itemize}
Note that under Condition E', Condition J' can be replaced by the stronger
condition $\lambda < \rho \lambda_0/\sqrt{s}$.

\begin{theorem} \label{thm::lasso-Mrho}
In the noise-free case, Conditions E' (or E), I' (or I), and J'
imply that for all $\beta_S\in\cM_\rho$,
there exists a lasso solution $\hat\beta$ with $\sgn(\hat\beta) = \sgn(\beta)$.
\end{theorem}

The conditions for exact reconstruction can be weakened. For instance,
Conditions E, I, and J' are also sufficient for exact reconstruction.
But we chose these forms because they transition nicely to the noisy case.
In a later section, we will discuss  Wainwright's result  \cite{Wainwright:2006} 
showing that a slight extension of Conditions E', I', and J' gives
sparsistency in the case of homoscedastic Gaussian noise.

\subsection{Exact reconstruction conditions for marginal regression in the noise-free case}

As above, define $\hat\alpha = X^T Y$ and define
 $\hat\beta$ by $\hat\beta_j = \hat\alpha_j 1\{|\hat\alpha_j| \ge t \}$, $1 \leq j \leq p$.  
For exact reconstruction with marginal regression,
we require
that $\hat\beta_j\ne 0$ whenever $\beta_j \ne 0$,
or equivalently $|\hat\alpha_j| \ge t$ whenever $\beta_j \ne 0$.
In the literature on causal inference,
this assumption is called
\emph{faithfulness}
\cite{Spirtes-etal:1993}
and is also   used in  \cite{Buhlmann, Fan-Lv:2008}.
The faithfulness assumption has received much
criticism \cite{Robins-etal:2003}.
The usual justification for faithfulness assumptions is that
if $\beta$ is selected at random from some distribution,
then faithfulness holds with high probability.
The criticism in \cite{Robins-etal:2003} is that
results which hold under faithfulness cannot hold in any uniform sense.

We write
$$
\hat\alpha = 
\left(
\begin{array}{ll}
\hat\alpha_S \\
\hat\alpha_N
\end{array} 
\right). 
$$
By elementary algebra 
$$
\hat\alpha =  
\left(
\begin{array}{ll}
X_S^T X_S \beta_S \\
X_N^T X_S \beta_S
\end{array} 
\right). 
$$
It follows directly that
\begin{itemize}
\item[]
{\it
{\bf Condition F}.~(Faithfulness)
\begin{equation} \label{eq::corrlearning} 
\max |C_{NS} \beta_S|  \,<\,  \min |C_{SS} \beta_S| 
\end{equation} 
}
\end{itemize}
is required to correctly identify the non-zero coefficients.
We call this the {\it Faithfulness Condition} 
even though it is technically different from the standard
definition of faithfulness above.
We thus have: 

\begin{lemma} 
Condition F is necessary and sufficient for exact reconstruction with marginal regression.
\end{lemma}

Unfortunately, as the next theorem shows, Condition F cannot
hold for all $\beta_S\in\cM_\rho$.
Applying the theorem to $C_{SS}$ shows that for any $\rho > 0$,
there exists a $\beta_S\in\cM_\rho$ that violates equation (\ref{eq::corrlearning}).

\begin{theorem} \label{thm::marginal-Mrho}
Let $C$ be an $s\times s$ positive definite, symmetric matrix that is not diagonal.
Then for any $\rho > 0$, there exists a $\beta\in \cM^s_\rho$
such that $\min |C \beta| = 0$.
\end{theorem}

Despite the seeming pessimism of Theorem \ref{thm::marginal-Mrho},  
we need to be cautious about over-interpreting  this result.
Since  $C \beta \equiv (Y, X)$,   
what  Theorem \ref{thm::marginal-Mrho}  says is that,   if we fix $X$ and  let $Y = X \beta$ ranges through all possible $\beta  \in \cM^s_\rho$,   then there exists  a $Y$ such that  $\min |(Y, X)|  = 0$.   However, both $X$ and $Y$ are observed, and if $\min |(Y,X)| > 0$, one can rule out the result of  Theorem \ref{thm::marginal-Mrho}.
Although Theorem \ref{thm::marginal-Mrho} is sufficient but not necessary for failure of marginal regression,
this mitigates the pessimism of the result.

\subsection{Comparison of the exact reconstruction conditions in the noise-free case}

In this section, we compare conditions for exact reconstruction in the noise-free
case required for the lasso and marginal regression.
We will see that at the level of individual $\beta$'s,
the conditions are generally overlapping and very closely related.
Although the conditions for marginal regression do not hold uniformly over any $\cM_\rho$,
they have the advantage that they do not require invertibility of $C_{SS}$
and hence are less sensitive to small eigenvalues.

We illustrate with a few examples, each in a subsection. 
\subsubsection{For an individual $\beta$,   the condition for the lasso and that  for marginal regression are generally overlapping}   
\begin{figure}
\begin{centering}
\includegraphics[height = 2.5 in,width = 3.5 in]{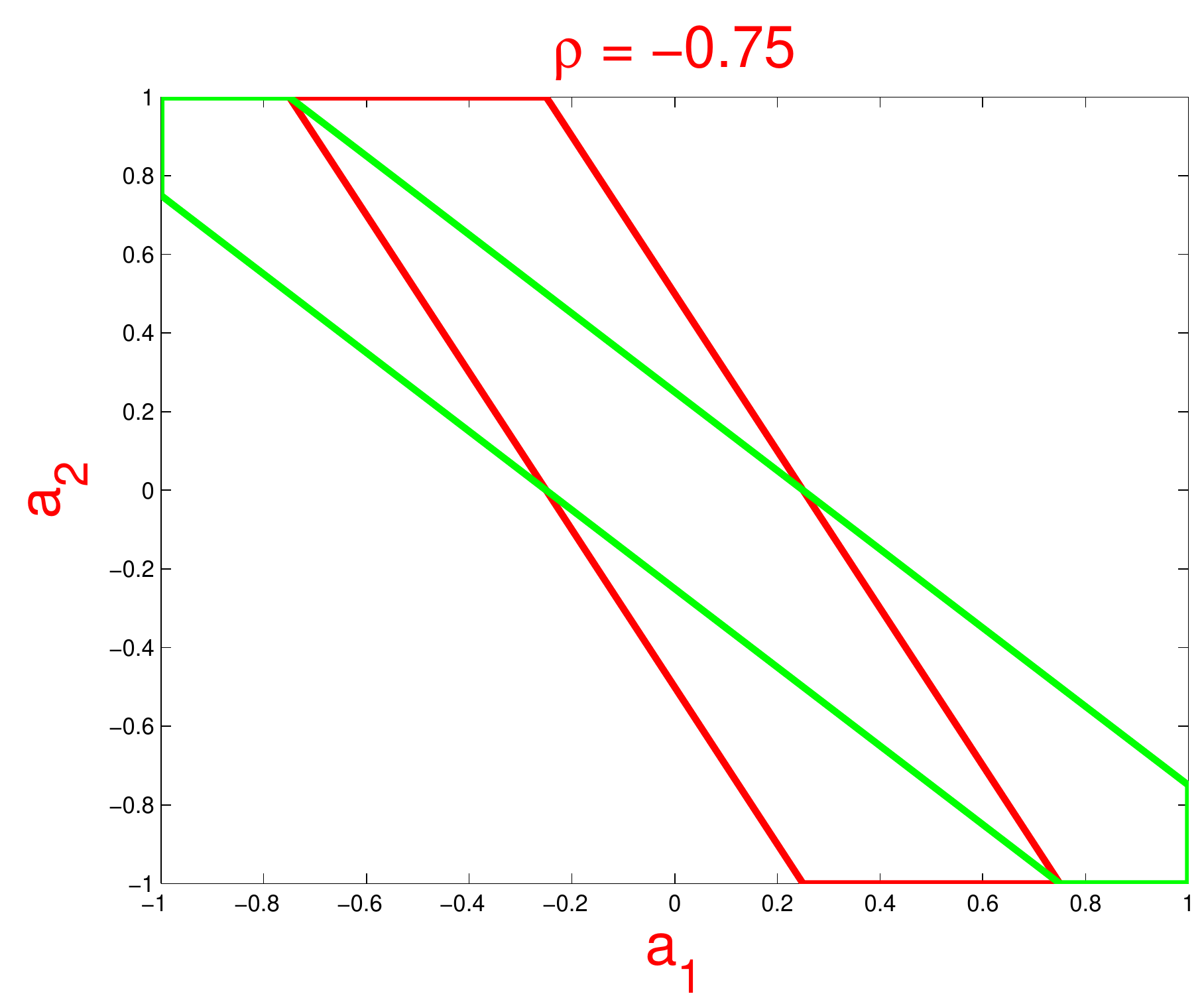}  \\
\caption{Let  $C_{SS}$ and  $\beta_S$ be as in Section \ref{subsec:example1}, where $\rho = -0.75$.   For a row of $C_{NS}$, say,  $(a_1, a_2)$.   The interior of the red box and that of the green box are the regions of   $(a_1, a_2)$ satisfying   the condition for the  lasso (i.e. (\ref{lassodisplay1}))  and for    marginal regression (i.e. (\ref{mrdisplay1})), respectively. }
\label{fig:Region1}
\end{centering}
\end{figure}

Consider an example where  
\[
C_{SS}  =   \left(  
\begin{array}{ll}
1  & \rho \\
\rho & 1 
\end{array} 
\right), \qquad  \beta_S  =  \left( \begin{array}{r} 
2 \\
1
\end{array} 
\right). 
\]
We investigate when $\rho$ ranges, which of two conditions  is weaker than the other.  Note that  $s = 2$ so  the matrix $C_{NS}$  only has two columns.  Fix a row of $C_{NS}$,  say  $a = (a_1, a_2)$.  Condition I requires  
\begin{equation} \label{lassodisplay1} 
|a_1 + a_2| \leq 1 + \rho, 
\end{equation} 
and Condition $F$ requires  
\begin{equation} \label{mrdisplay1}
|2 a_1 + a_2|  \leq \min\{(2 + \rho), (1 + 2 \rho)\}. 
\end{equation} 
Seemingly,  for many choices  of $\rho$,    two conditions (\ref{lassodisplay1}) and (\ref{mrdisplay1}) overlap  with each other.   Take $\rho = -0.75$ for illustration. In    Figure \ref{fig:Region1}, we display  the regions where  $(a_1, a_2)$ satisfy  (\ref{lassodisplay1})  and  (\ref{mrdisplay1}), respectively.     The figure shows that two regions are  overlapping.   As a result,  the condition for the lasso overlaps  with that for   marginal regression.   For different choices of $\rho$ and $\beta_S$,  those  regions in Figure \ref{fig:Region1}  may vary,  but to a large extent, two conditions continue to overlap with each other.

Examples for larger $s$ can be constructed  by letting  $C_{SS}$ be a  block diagonal matrix,  where the size of  each main diagonal block is small.     For each row of $C_{NS}$, the conditions 
for  the  lasso and  marginal regression are  similar to  those in  (\ref{lassodisplay1}) and 
(\ref{mrdisplay1}), respectively,   but maybe more complicated.  To save space, we omit further discussion along this line.

\subsubsection{In the special case of $\beta_S \propto 1_S$,  the condition for the lasso and the condition for marginal regression are  closely related}    \label{subsec:example1} 
Consider the special case  
$\beta_S  \propto 1_S$.  In this case,  
the main condition for the lasso (Condition I)  is 
\begin{equation} \label{C1}
|C_{NS}   \cdot  C_{SS}^{-1}  \cdot  1_S |   \;\;      \leq  \;\;   1_S, 
\end{equation} 
and the condition for marginal regression (Condition F)  is 
\begin{equation} \label{C2}
|C_{NS}  \cdot 1_S |   \;\;    \leq   \;\;   | C_{SS}  \cdot 1_S|,  
\end{equation} 
where both inequalities should be interpreted as hold component-wisely.  
Two  conditions are surprisingly similar:    removing  
 $C_{SS}^{-1}$ on the left side of (\ref{C1}) and adding  $C_{SS}$ to the right hand side of it gives 
   (\ref{C2}).  
  
Note that if in addition $1_S$ is an eigen-vector of $C_{SS}$, then  two conditions are  equivalent to each other.   This includes but is not limited to the case of  $s = 2$.

\subsubsection{In the special case of $C_{SS} = I$, the condition for the lasso is weaker}
Fix $n$ and consider the special case in which $C_{SS} = I$.
For the lasso, Condition E' (and thus E) is satisfied, Condition J' reduces to $\lambda < \rho$,
and Condition I becomes $\norm{C_{NS}} \le 1$.
Under these conditions, the lasso gives exact reconstruction, but Condition F can fail.
To see how, let $\tilde\beta\in\{-1,1\}^s$ be the vector such that
$\max |C_{NS} \tilde \beta| = \norm{C_{NS}}_\infty$
and let $\ell$ be the index of the row at which the maximum is attained,
choosing the row with the biggest absolute element if the maximum is not unique.
Let $u$ be the maximum absolute element of row $\ell$ of $C_{NS}$
with index $j$.
Define a vector $\delta$ to be zero except in component $j$, which
has the value $\rho\tilde\beta_j/(u \norm{C_{NS}}_\infty)$.
Let $\beta = \rho \tilde\beta + \rho \delta$.
Then,
\begin{align}
|(C_{NS}\beta)_\ell|
  &= \rho |(C_{NS}\tilde \beta)_\ell| + \rho \frac{1}{\norm{C_{NS}}_\infty} \\
  &= \rho \left(\norm{C_{NS}}_\infty + \frac{1}{\norm{C_{NS}}_\infty}\right) \\
  &> \rho.
\end{align}
It follows that $\max |C_{NS}\beta| > \rho = \min |\beta|$, so Condition F fails.

On the other hand, suppose Condition F holds for all $\beta_S\in\{-1,1\}^s$.
(It cannot hold for all $\cM_\rho$ by Theorem \ref{thm::marginal-Mrho}). 
Then, for all $\beta_S\in\{-1,1\}^s$, $\max |C_{NS} \beta_S| \le 1$,
which implies that $\norm{C_{NS}}_\infty \le 1$. Choosing $\lambda < \rho$,
we have Conditions E', I, and J' satisfied, showing by Theorem \ref{thm::lasso-Mrho}
that the lasso gives exact reconstruction.

It follows that the conditions for the lasso are weaker in this case.

\subsubsection{Small eigenvalues of $(X_S' X_S)$ may have an adverse  
effect on the performance of  the lasso, but not always on that for marginal regression}   \label{subsec:compareexample4}
For simplicity,  assume  
\[
\beta_S \propto 1_S.
\] 
(We remark that the phenomenon to be described below 
is not limited to the case of $\beta_S \propto 1_S$). 
For $1 \leq i \leq s$, let  $\lambda_i$ and $\xi_i$ be the $i$-th  eigenvalue and eigenvector   of $C_{SS}$.  Without loss of generality, we assume that $\xi_i$ have unit  $\ell^2$ norm.  By elementary algebra,  there are constants $a_1, \ldots, a_s$ such that 
$1_S        =   c_1 \xi_1 + c_2 \xi_2 + \ldots  + c_s   \xi_s$.  
It follows 
\[
C_{SS}^{-1}  \cdot  1_S    =  \sum_{i =1}^s  \frac{c_i}{\lambda_i} \xi_i   \qquad \mbox{and} \qquad   
C_{SS} \cdot  1_S    =  \sum_{i =1}^s (c_i \lambda_i) \xi_i.  
\]
Fix a row of $C_{NS}$, say, $a = (a_1, \ldots, a_s)$.   
Respectively, the conditions for the lasso  and marginal regression  require 
\begin{equation} \label{Editadd1}
|(a,  \sum_{i =1}^s  \frac{c_i}{\lambda_i} \xi_i)| \leq 1  \qquad   \mbox{and} \qquad |(a,1_S)|  \leq  |\sum_{i =1}^s (c_i \lambda_i) \xi_i|. 
\end{equation} 

Without loss of generality, we assume that  $\lambda_1$ is the smallest eigenvalue of $C_{SS}$.    Consider the case where $\lambda_1$ is  small, while all other eigenvalues have a magnitude comparable to $1$.    In this case, the  smallness of  $\lambda_1$ has a negligible effect on $ \sum_{i =1}^s (c_i \lambda_i) \xi_i$,  
and so has a  negligible effect on the condition for marginal regression.  However,  
the smallness of  $\lambda_1$ may have an adverse  effect on the performance of  the lasso.  To see the point, we note that 
$\sum_{i = 1}^s  \frac{c_i}{\lambda_i} \xi_i     \approx  \frac{c_1}{\lambda_1} \xi_1$. 
Compare this with the first term in (\ref{Editadd1}).  The condition for the lasso is roughly 
\[
|(a, \xi_1)| \leq \lambda_1, 
\]
which is rather  restrictive since $\lambda_1$ is small.   

We now further illustrate the point with an example.  
Let  \[ 
C_{SS}     =   \left(
\begin{array}{ccc} 
1       & -1/2    &c      \\
-1/2   &1      &0  \\ 
c    & 0  &1 
\end{array}  
 \right). 
\]
The smallest eigenvalue  of the matrix is  $\lambda_1  = \lambda_1(c)  
 = 1 -  \sqrt{c^2 + 1/4}$,   which is positive if and only if $c \leq \sqrt{3}/2$.  
Suppose $0 < c < \sqrt{3}/2$.   Fix a row of $C_{NS}$, say, $a = (a_1, a_2, a_3)$.   By direct calculations, the condition for marginal regression and the lasso  require  
\begin{equation} \label{compareexample4} 
|a_1 + a_2 + a_3| \leq 1/2,  \; \mbox{and} \;    |(6 - 4c) a_1 + (6 - 3c - 4 c^2) a_2 + (3 - 6 c) a_3| \leq (3 - 4 c^2) 
\end{equation} 
respectively.  As  $c$ approaches $\sqrt{3}/2$,  both  $\lambda_1(c)$ and the right hand side of the first inequality  in (\ref{compareexample4})  approach $0$.  As a result, the first inequality in (\ref{compareexample4})  becomes increasingly more restrictive, but the 
the second inequality remains the same for all $c$.  Therefore, a small $\lambda_1(c)$ has a negative effect on the broadness of the condition for the lasso, but not on that of  marginal regression. 

\begin{figure}
\begin{centering}
\includegraphics[height = 2.5 in,width = 3.5 in]{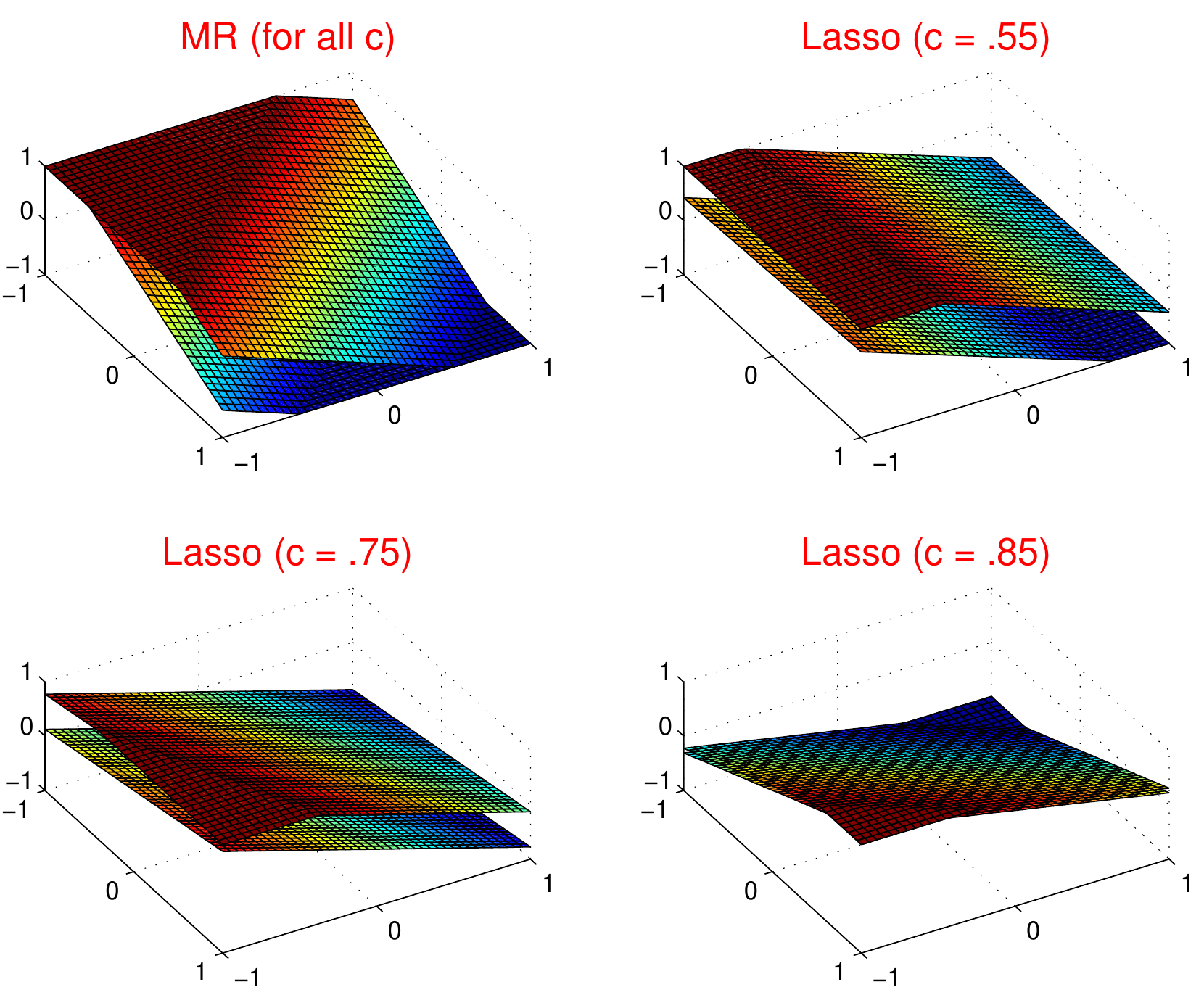}  \\
\caption{The regions sandwiched by two hyper-planes are  the regions of  $a = (a_1, a_2, a_3)$ satisfying the condition  for marginal regression  (first inequality in (\ref{compareexample4}); Panel 1)  and for the lasso (second inequality of (\ref{compareexample4}; Panel 2--4).  Details of the example are  in Section \ref{subsec:compareexample4}. Here,   $c = 0.55, 0.75, 0.85$ and the smallest eigenvalues of $C_{SS}$ are  $\lambda_1(c) = 0.29, 0.14, 0.014$.   As $c$ varies,  
the regions for marginal regression remain  the same so is  displayed only   in Panel 1 (MR stands for marginal regression).    
In contrast, the regions for   the lasso get substantially smaller as  $\lambda_1(c) $  decreases.    }
\label{fig:Region2}
\end{centering}
\end{figure}

 Figure \ref{fig:Region2}   displays  the regions where the vector $a$ satisfies the first and the second  inequality in (\ref{compareexample4}), respectively.  In this example,   $c = 0.55,
0.75, 0.85$, so   $\lambda_1(c)  = 0.29, 0.14,
0.014$ correspondingly.   To better visualize these regions,  we display their  2-D section   in Figure \ref{fig:Region2} (in the 2-D section, we set the first coordinate of $a$ to $0$). The figures suggest that when $\lambda_1(c)$ get increasingly smaller, the region corresponding to the lasso shrinks substantially, while that corresponding to marginal regression remains the same.

\begin{figure}
\begin{centering}
\includegraphics[height = 2.5 in,width = 3.5 in]{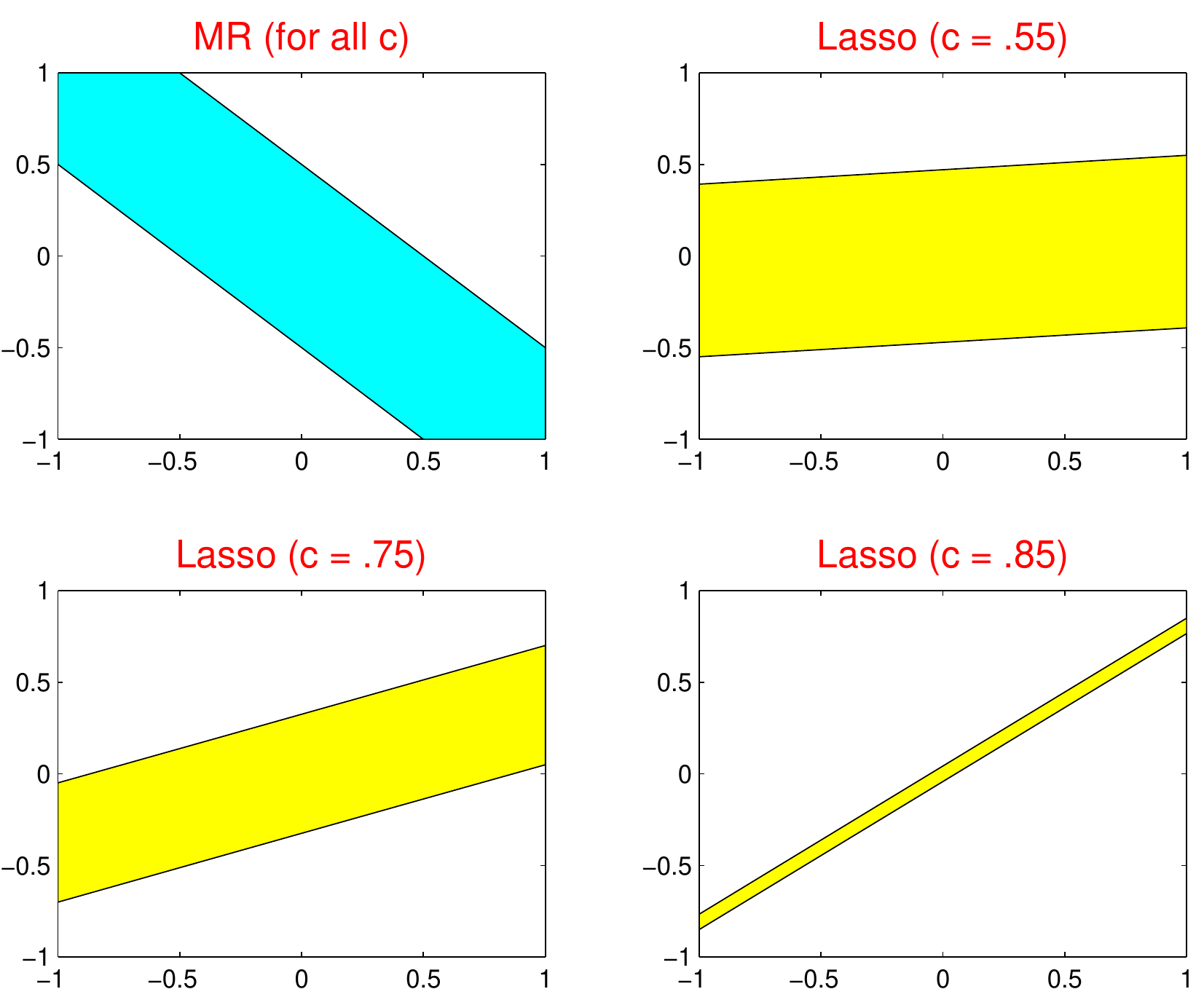}  \\
\caption{Displayed are  the 2-D sections of the regions in Figure \ref{fig:Region2},    where we set the first coordinate of $a$   to $0$.     As $c$ varies,  
the regions for marginal regression remain  the same, but those for the lasso get substantially smaller as  $\lambda_1(c)$ decrease.  $x$-axis: $a_2$. $y$-axis: $a_3$. }
\label{fig:Region3}
\end{centering}
\end{figure}

In conclusion, in the noise-free case,   the condition for the lasso and for marginal regression   are generally overlapping and closely related. On one hand, given a fixed matrix $X$, the conditions for the success of the lasso
hold over a larger class of $\beta$'s than that of marginal regression.  On the other hand,
marginal regression has a larger tolerance for collinearity than does the lasso.  Bear in mind  that for very large $p$,  the lasso is   computationally  much more demanding.

\subsection{Exact  reconstruction conditions for  marginal regression in the noisy case}   \label{subsec:noise} 
We now come back to Model (5)  and consider the noisy case:  
\begin{equation} \label{model2}
Y^{(n)} =  X^{(n)} \beta^{(n)}  + z^{(n)},  \qquad \mbox{where} \qquad  z^{(n)} \sim N(0, \sigma_n^2 \cdot I_n). 
\end{equation} 
To focus on variable selection, we suppose the parameter $\sigma_n^2$ is known.  The exact reconstruction condition  for  the lasso in the noisy 
case has been studied extensively in the literature (see for example \cite{Tibshirani:1996}).  So in this paper,  we focus on that for  marginal regression.  We address  two topics. 
First, we extend Condition F in the noise-free case  to  the noisy case, say Condition F'.    When Condition F' holds, we   show  that  with an appropriately chosen threshold $t$ (see (\ref{Definemarginal})),   marginal regression  fully recovers the support  with high probability.    Second,  we discuss how to determine  the threshold $t$  empirically.   

Recall that in the noise-free case,  Condition F is
\[
\max |C_{NS} \beta_S| \leq  \min |C_{SS} \beta_S|. 
\]
A natural extension of Condition F is the following. 
\begin{itemize}
\item[]
{\it
{\bf Condition F'}.~(Faithfulness)
\begin{equation} \label{eq::corrlearning1} 
\max |C_{NS} \beta_S|  +  2 \sigma_n   \sqrt{2 \log p}   \,<\,  \min |C_{SS} \beta_S|.  
\end{equation} 
}
\end{itemize}
When Condition F' holds,    it is possible to separate relevant variables from irrelevant variables with high probability. In detail, 
write 
\[
X = [x_1, x_2,\ldots, x_p], 
\]
where $x_i$ denotes  the $i$-th column of $X$.    Sort $|(Y, x_i)|$ in the descending order, and let $r_i = r_i(Y, X)$ be the ranks of $|(Y, x_i)|$ (assume no ties for simplicity). 
Introduce 
\[
\hat{S}_n(k)  =  \hat{S}_n(k;  X, Y,  p) =  \{i: \;  r_i(X, Y)  \leq  k \}, \qquad k = 1,2, \ldots, p. 
\]
Recall that $S(\beta)$ denotes the support of $\beta$ and $s = |S|$.     
The following lemma says that, if $s$ is known and Condition F' holds, then marginal regression is able to fully recover the support $S$  with high probability.  
\begin{lemma}  \label{lemma:F'}
Consider a sequence of regression models as    in (\ref{model2}). If for sufficiently large $n$,   Condition F' holds  and $p^{(n)} \geq n$, then 
\[
\lim_{n \goto \infty}  P \biggl(\hat{S}_n(s^{(n)}; X^{(n)}, Y^{(n)},  p^{(n)} )   \neq  S(\beta^{(n)})  \biggr)   =  0.   
\]  
\end{lemma}   
Lemma \ref{lemma:F'} is proved in the appendix. We remark that  if both $s$ and $(p-s)$ tend  to $\infty$ as $n$ tends to $\infty$,  then Lemma \ref{lemma:F'}  continues to hold if we replace  
$2 \sigma_n \sqrt{2 \log p}$ in (\ref{eq::corrlearning1}) by $\sigma_n(\sqrt{\log(p-s)} + \sqrt{\log s})$.  See the proof of the lemma for details.  

The key assumption of  Lemma \ref{lemma:F'} is that   $s$ is known so that   we know how to set the threshold $t$.  
Unfortunately,   $s$ is generally  unknown.  We  propose the following procedure  to estimate $s$.  Fix $1 \leq k \leq p$, let $i_k$ be the unique index satisfying  
\[
r_{i_k}(X, Y) =  k. 
\] 
Let $\hat{V}_n(k)  = \hat{V}_n(k; X, Y, p)$ be the linear space spanned by 
\[
x_{i_1}, x_{i_2}, \ldots, x_{i_k}, 
\]
and let $\hat{H}_n(k) = \hat{H}_n(k; X, Y,  p)$ be the projection matrix from $R^n$ to $\hat{V}_n(k)$ (here and below,  the $\hat{\;}$ sign  emphasizes the dependence of     indices $i_k$  on the data).  Define
\[
\hat{\delta}_n(k) = \hat{\delta}_n(k;  X, Y,  p) = \|(\hat{H}_n(k + 1) - \hat{H}_n(k))Y\|, \qquad 1 \leq k \leq p -1. 
\]
The term $\hat{\delta}_n^2(k)$ is closely related to  the F-test  for testing whether $\beta_{i_{k+1}} \neq 0$.   
We estimate $s$ by     
\[
\hat{s}_n =  \hat{s}_n(X, Y, p)   =  \max \biggl\{1 \leq  k \leq p: \;  \hat{\delta}_n(k) \geq \sigma_n \sqrt{2 \log n}  \biggr\}  + 1
\]
(in the case where $\hat{\delta}(k) < \sigma_n \sqrt{2 \log n}$ for all $k$, we define $\hat{s}_n = 1$).

Once $\hat{s}_n$ is determined, 
we estimate the support  $S$    by 
\[
\hat{S}(\hat{s}_n, X, Y, p) =  \{i_k:  \, k  = 1, 2, \ldots,   \hat{s}_n \}. 
\]
It turns out that under mild conditions,  $\hat{s}_n = s$ with high probability.   In detail, 
suppose that the support $S(\beta)$ consists of indices $j_1, j_2, \ldots, j_s$.  
Fix $1 \leq k \leq s$.  Let $\tilde{V}_S$ be the linear space spanned  by $x_{j_1}, \ldots, x_{j_s}$, and let 
$\tilde{V}_{S, (-k)}$ be the linear  space  spanned by $x_{j_1}, \ldots, x_{j_{k-1}}, x_{j_{k+1}}, \ldots,  x_{j_s}$.   Project $\beta_{j_k} x_{j_k}$ to the linear  space $\tilde{V}_S \cap \tilde{V}_{S, (-k)}^{\perp}$.    Let  
$\Delta_n(k,  \beta,    X,     p)$ 
be the  $\ell^2$  norm  of the resulting vector, and let  
\[
\Delta_n^*(\beta, X, p)  =  \min_{1 \leq k \leq s}  \Delta_n(k, \beta, X, p).  
\]
The following theorem says that if $\Delta_n^*(\beta, X, p)$ is slightly larger than $\sigma_n \sqrt{2 \log n}$, then $\hat{s}_n = s$ and $\hat{S}_n = S$ with high probability.  In other words,  marginal regression  fully recovers the support with high probability. Theorem \ref{thm:noise} is proved in the appendix. 
\begin{theorem} \label{thm:noise}
Consider a sequence of regression models as    in (\ref{model2}).  Suppose that   for sufficiently large $n$,     Condition F' holds,  $p^{(n)}  \geq n$, and 
\[
\lim_{n\goto \infty}   \biggl( \frac{\Delta_n^*(\beta^{(n)}, X^{(n)}, p^{(n)})}{\sigma_n} - \sqrt{2 \log n}  \biggr)   = \infty.  
\]
Then 
\[
\lim_{n \goto \infty}  P \biggl( \hat{s}_n(X^{(n)}, Y^{(n)}, p^{(n)})  \neq  s^{(n)} \biggr)  \goto 0, 
\]
and 
\[
\lim_{n\goto \infty} \biggl(\hat{S}_n(\hat{s}_n(X^{(n)},Y^{(n)}, p^{(n)}); X^{(n)}, Y^{(n)}, n, p^{(n)} )   \neq  S(\beta^{(n)})  \biggr)  \goto 0.   
\]  
\end{theorem} 
Theorem \ref{thm:noise}  says  that  the tuning parameter for marginal regression (i.e.  the threshold $t$) can be set 
successfully in a data driven fashion. In comparison, how to set the tuning parameter $\lambda$ for the lasso  
has been a withstanding open problem in the literature.

\section{The Deterministic Design, Random Coefficient Regime}   \label{sec:exten} 
Recall that the Faithfulness Condition is 
\[
\max | C_{NS} \beta_S| \leq  \min |C_{SS} \beta_S|. 
\]
In this section, we study how broad the Faithfulness Condition holds.  We approach this by modeling $\beta$ as random (the matrix $X$ is kept deterministic), and find out conditions under which  the Faithfulness Condition holds with high probability. 

The discussion in this section  is closely related to the 
work by Donoho and Elad \cite{Donoho-Elad:2003} on the Incoherence Condition.  Compared to the Faithfulness Condition, the advantage of the Incoherence Condition is that it does not involve the unknown support of $\beta$, 
so it is checkable in practice.  The downside of the Incoherence Condition is that it aims to control the worst case so it is conservative.   In this section, we derive a condition---Condition F''---
which can be viewed as a middle  ground between the Faithfulness Condition and the Incoherence Condition: 
it is not tied to the unknown support so it is more tractable than the Faithfulness Condition,  and it is also much less stringent than the  Incoherence Condition.

In detail, we model $\beta$ as follows.   Fix $\eps \in (0,1)$, $a > 0$, and a distribution $\pi$, where   
\begin{equation} \label{pi1}
\mbox{the support of $\pi$}  \qquad   \subset   \qquad  (-\infty, -a] \cup [a, \infty). 
\end{equation} 
For each $1 \leq i \leq p$, we draw a sample $B_i$ from $\mathrm{Bernoulli}(\eps)$.  When $B_i  = 0$, we set $\beta_i = 0$. When $B_i = 1$, we draw $\beta_i \sim  \pi$. Marginally,   
\begin{equation} \label{Definebeta}
\beta_i \stackrel{iid}{\sim} (1 - \eps) \nu_0 + \eps \pi,   
\end{equation} 
where $\nu_0$ denotes  the point mass at $0$.    We study for which quadruplets $(X, \eps, \pi, a)$      the Faithfulness Condition holds with high probability.  

Recall that the design matrix $X = [x_1, \ldots, x_p]$, where  $x_i$ denotes  the $i$-th column. 
Fix $t \geq 0$ and $\delta > 0$.   Introduce  
\[
g_{ij}(t)  =    E_{\pi} [e^{t  u (x_i, x_j)}] - 1,     
 \qquad  \bar{g}_i(t) =   \sum_{j \neq i} g_{ij}(t), 
\]
where the random variable $u \sim \pi$. As before,  we have suppressed the superscript $\,^{(n)}$ for $g_{ij}(t)$ and $\bar{g}_i(t)$.  Define   
\[
A_n(\delta, \eps, \bar{g}) =  A_n(\delta, \eps, \bar{g}; X, \pi) =  \min_{t > 0}  \biggl( e^{- \delta t}   \sum_{i =1}^p [e^{ \eps \bar{g}_i(t)} +  e^{\eps \bar{g}_i(-t)}] \biggr), 
\]
where $\bar{g}$ denotes the vector $(\bar{g}_1, \ldots, \bar{g}_p)^T$.  
The following lemma is proved in the appendix. 
\begin{lemma}  \label{lemma:coherent}  
Fix $n$, $X$,  $\delta > 0$,   $\eps \in (0,1)$,      and distribution $\pi$.   Then 
\begin{equation} \label{lemma1} 
P( \max |C_{NS} \beta_S|  \geq \delta )  \leq (1 - \eps)  A_n(\delta, \eps, \bar{g}; X, \pi), 
\end{equation} 
and 
\begin{equation} \label{lemma2} 
\qquad P(\max |(C_{SS} - I_S) \beta_S| \geq \delta)  \leq \eps  A_n(\delta, \eps,\bar{g};   X, \pi). 
\end{equation} 
\end{lemma}

Now, suppose  the distribution $\pi$ satisfies (\ref{pi1}) for some $a > 0$. Take $\delta = a/2$ on the right hand side of  (\ref{lemma1})-(\ref{lemma2}).  Except for a probability of 
$A_n(a/2, \eps, \bar{g})$, 
\[
\max |C_{NS} \beta_S| \leq a/2,  \qquad \min|C_{SS}\beta_S|  \geq \min |\beta_S|  - \max|(C_{SS} - I) \beta_S|   \geq a/2, 
\]
so $\max|C_{NS} \beta_S| \leq \min |C_{SS}\beta_S|$ and the Faithfulness Condition holds.  This motivates the following condition,   where $(a, \eps, \pi)$ may depend on $n$. 
\begin{itemize}
\item[]
{\it
{\bf Condition F''}.~(Faithfulness)
\begin{equation} \label{eq::corrlearning2} 
\lim_{n \goto \infty}  A_n(a_n/2, \eps_n, \bar{g}^{(n)}; X^{(n)}, \pi_n)   = 0. 
\end{equation} 
}
\end{itemize}
The following theorem  says that if   Condition F'' holds, then Condition F holds with high probability.  
\begin{theorem} \label{thm:faithful2}
Consider a sequence of noise-free regression models as    in (\ref{model2}),  where the noise component $z^{(n)} = 0$ and $\beta^{(n)}$ is generated as in (\ref{Definebeta}).   Suppose  Condition F'' holds. Then as $n$ tends to $\infty$,   except for a probability that tends to $0$,   
\[
\max |C_{NS} \beta_S |  \leq  \min |C_{SS} \beta_S|. 
\]
\end{theorem} 
Theorem \ref{thm:faithful2} is the direct result of Lemma \ref{lemma:coherent} so we omit the proof.

\subsection{Comparison of Condition F''  with the Incoherence Condition} 
Introduced in Donoho and Elad \cite{Donoho-Elad:2003} (see also \cite{Donoho-Huo:2001}), 
the {\it Incoherence} of a matrix $X$ is defined as \cite{Donoho-Elad:2003}
\[
   \max_{i \neq j}   |C_{ij}|, 
\] 
where $C = X^TX$ is the Gram matrix as before. The notion is motivated by the study in recovering a sparse 
signal from an over-complete dictionary.  
In the special case where $X$ is the concatenation 
two orthonormal bases  (e.g.  a  Fourier basis and a wavelet basis), 
$ \max_{i \neq j}   |C_{ij}|$ measures how coherent two bases are and so the term of incoherence;
see \cite{Donoho-Elad:2003, Donoho-Huo:2001} for details.  Consider Model      (\ref{eq::basic-regression})  in the case where both $X$ and $\beta$ are deterministic, and the noise component $z = 0$.  The following
results are proved in \cite{Chen-Donoho:1998, Donoho-Elad:2003,
Donoho-Huo:2001}.
\begin{itemize} 
\item Lasso   yields exact variable selection if   $s  <  \frac{1 + \max_{i \neq j}|C_{ij}|}{2 \max_{i\neq j} |C_{ij}|}$.    
\item Marginal regression yields exact variable selection if  $s < \frac{c}{2 \max_{i \neq j} |C_{ij}|}$ for some constant $c \in (0,1)$, and that the nonzero coordinates of $\beta$ have comparable  magnitudes (i.e. 
the ratio between  the largest and the smallest  nonzero coordinate of $\beta$  is bounded away from $\infty$). 
\end{itemize} 

In comparison,   the Incoherence Condition only depends   on $X$ so it is checkable.    Condition F depends on  the unknown support of $\beta$. Checking such a condition is almost as hard as estimating the support $S$. 
Condition F'' provides a middle ground. It depends on $\beta$ only through $(\eps, \pi)$. In cases where we either have a good knowledge of $(\eps, \pi)$ or we can estimate them,    Condition F'' is checkable.

At the same time,  the Incoherence Condition is  conservative, especially when $s$ is large.  In fact, in order for either the lasso or 
marginal regression   to  have  an exact  variable selection, 
it is required that 
\begin{equation}  \label{Incoherence} 
\max_{i \neq j}  |C_{ij}|  \leq O\biggl(\frac{1}{s}\biggr),  
\end{equation}  
In other words,  all coordinates of the Gram matrix $C$ need to be no greater than $O(1/s)$.    This is 
much more conservative than Condition F. 

However, we must note that the Incoherence Condition aims to control the worst case:  it sets out to guarantee  {\it uniform} success of a
procedure across all $\beta$ under minimum constraints.    In comparison, Condition F aims to control  a single case, and Condition F''  aims to control almost all the cases in a specified class.   As such,      Condition F''   provides  a middle ground between  Condition F and the Incoherence Condition, 
applying more broadly than the former, while being less conservative 
than the later.

Below, we use two examples to illustrate that Condition F'' is much less conservative  than the Incoherence Condition. 
   In the first example, we consider a weakly dependent case where  $\max_{i \neq j }  |C_{ij}| \leq O(1/\log(p))$.     In the second example, we suppose the matrix $C$  is  sparse,  but the nonzero coordinates of $C$  may be large.

\subsubsection{The weakly dependent case} 
Suppose that for sufficiently large $n$,  there are two sequence of  positive numbers $a_n \leq  b_n$ such that
\[
\mbox{the support of $\pi_n$} \qquad \subset  \qquad  [-b_n, -a_n] \cup [a_n, b_n],  
\] 
and that 
\[
\frac{b_n}{a_n} \cdot \max_{i\neq j} |C_{ij}|  \leq c_1 /\log(p),  \qquad \mbox{$c_1 > 0$ is a constant}.   
\]
For $k \geq 1$, denote the $k$-th  moment of $\pi_n$ by 
\begin{equation} \label{pi2}
\mu_n^{(k)}  = \mu_n^{(k)}(\pi_n). 
\end{equation} 
Introduce $m_n = m_n(X)$ and $v_n^2 = v_n^2(X)$ by 
\[
m_n(X)  =   p \eps_n   \cdot   \max_{1 \leq i \leq p}  \biggl\{ \biggl| \frac{1}{p}  \sum_{j \neq i} C_{ij}  \biggr|  \biggr\}, \qquad 
 v_n^2(X)     =  p \eps_n \cdot  \max_{1 \leq i \leq p}  \biggl\{ \frac{1}{p}   \sum_{j \neq i} C^2_{ij}  \biggr\}. 
\]
\begin{cor} \label{lemma:weakdependent} 
Consider a sequence of  regression models as  in (\ref{model2}),  where the noise component $z^{(n)} = 0$ and $\beta^{(n)}$ is  generated as in  (\ref{Definebeta}). 
If there are constants   $c_1 > 0$ and $c_2 \in (0,1/2)$ such that 
\[
\frac{b_n}{a_n} \cdot \max_{i \neq j } \{  |C_{ij}|  \} \leq c_1 /\log(p^{(n)}), 
\]
and 
\begin{equation} \label{weakcondition} 
\limsup_{n \goto \infty}  \biggl( \frac{\mu_n^{(1)}(\pi_n)}{a_n}     m_n(X^{(n)}) \biggr)  \leq  c_2,  \qquad \limsup_{n \goto \infty}  \biggl(\frac{\mu_n^{(2)}(\pi_n)}{a_n^2}      v_n^2(X^{(n)})  \log(p^{(n)}) \biggr)  = 0,  
\end{equation} 
then 
\[
\lim_{n \goto \infty}  A_n(a_n/2, \eps_n, \bar{g}^{(n)}; X^{(n)}, \pi_n)   = 0, 
\]
and Condition F'' holds. 
\end{cor} 
Corollary  \ref{lemma:weakdependent} is proved in the appendix. 
For interpretation, we consider   the special case where there is a generic constant $c  > 0$  such that $b_n \leq c  a_n$. As a result,    $\mu_n^{(1)}/a_n \leq c $,  $\mu_n^{(2)}/a_n^2 \leq c^2$. 
The conditions   reduce  to that,    for sufficiently large $n$ and all $1 \leq i \leq p$,  
\[
| \frac{1}{p}   \sum_{j \neq i}^p  C_{ij}|   \leq  O(\frac{1}{p \eps_n}), \qquad   \frac{1}{p} \sum_{j \neq i}^p C^2_{ij}  = o(1/p \eps_n).
\]
Note that by (\ref{Definebeta}),  $s  = s^{(n)} \sim \mathrm{Binomial}(p, \eps_n)$,    so  $s\approx p \eps_n$.  
Recall that  the  Incoherence Condition is   
\[
\max_{i \neq j}  |C_{ij}|   \leq O(1/s). 
\]
In comparison, the   Incoherence Condition requires that each   coordinate of $(C-I)$ is no greater than 
$O(1/s)$, while  Condition F''  only requires that 
the average of each row of $(C-I)$   is no greater than  $O(1/s)$.
The latter  is much less  conservative. 

\subsubsection{The sparse case} 
Let $N_n^*(C)$ be the maximum number of nonzero off-diagonal coordinates of $C$: 
\[ 
 N_n^*(C) = \max_{1 \leq i \leq p} \{N_n(i)\}, \qquad N_n(i) = N_n(i; C) = \#\{j: \; j  \neq i,   C_{ij}  \neq  0\}. 
\]
Suppose there is a constant $c_3 > 0$ such that 
\begin{equation} \label{sparse2}
\liminf_{n \goto \infty}  \biggl( \frac{-\log(\eps_n N_n^*(C))}{\log (p^{(n)})}  \biggr) \geq  c_3. 
\end{equation} 
Also, suppose there is a constant $c_4 > 0$ such  that for sufficiently large $n$,  
\begin{equation} \label{sparse1} 
\mbox{the support of $\pi_n$ is contained in $[-c_4 a_n,  a_n] \cup [a_n, c_4 b_n]$}. 
\end{equation} 
The following corollary is proved in the appendix. 

\begin{cor} \label{lemma:sparsespikes} 
Consider a sequence of noise-free regression models as    in (\ref{model2}),  where the noise component $z^{(n)} = 0$ and $\beta^{(n)}$ is randomly generated as in (\ref{Definebeta}).  Suppose  (\ref{sparse2})-(\ref{sparse1}) hold. 
If there is a constant $\delta > 0$ such that 
\begin{equation} \label{sparsecondition} 
\max_{i  \neq j} |C_{ij}|  \leq   \delta, \qquad \mbox{and} \qquad \delta < \frac{c_3}{2c_4}, 
\end{equation} 
then 
\[
\lim_{n \goto \infty}  A_n(a_n/2, \eps_n, \bar{g}^{(n)}; X^{(n)}, \pi_n)   = 0, 
\]
and Condition F'' holds. 
\end{cor}

For interpretation, consider a special case where  
\[
  \eps_n = p^{-\vartheta}.  
\]
In this case, 
the condition   reduces to     
\[ 
 N_n^*(C)  \ll p^{\vartheta -   2 c_4 \delta}.  
\]
As a result,   Condition F'' is satisfied if each row of $(C-I)$ contains no more than $p^{\vartheta - 2 c_4  \delta}$  nonzero coordinates each of which   $\leq  \delta$. 
Compared to the  Incoherence Condition  $\max_{i \neq j} |C_{ij}|  \leq O(1/s) = O(p^{-\vartheta})$,  our  condition   is much weaker.  

In conclusion, if we alter our attention from the worst-case scenario to the average scenario, and alter our aim  from exact variable selection   to exact variable selection with probability $\approx 1$,  then the   
condition required for success---Condition F''---is 
much more relaxed  than the Incoherence Condition.


\section{Hamming Distance when $X$ is Gaussian;  Partition  of the Phase Diagram}  \label{sec:hamm} 
So far, we have focused on  exact variable selection.  In many applications, exact variable selection 
is not possible. Therefore,   it is of interest to study  the  Type I and Type II errors of variable selection (a Type I error is a misclassified $0$ coordinate of $\beta$,  and a Type II error is a misclassified nonzero coordinate).  

In this section, we use the Hamming distance to measure the  variable selection errors. 
Back to Model (\ref{eq::basic-regression}),  
\begin{equation} \label{Gaussian0}
Y = X \beta + z, \qquad z   \sim N(0,  I_n),   
\end{equation} 
where without loss of generality, we assume 
$\sigma_n = 1$. 
As in the preceding section (i.e. (\ref{Definebeta})),  we suppose 
\begin{equation} \label{Definebeta1}
\beta_i \stackrel{iid}{\sim} (1 -\eps) \nu_0  + \eps \pi.    
\end{equation} 
For any variable selection procedure 
$\hat{\beta}  =  \hat{\beta}(Y; X)$,    the  Hamming distance between $\hat{\beta}$ and the true $\beta$ is  
\[
d(\hat{\beta}  |    X)  =  d(\hat{\beta}; \eps, \pi |X) =  \sum_{j = 1}^p E_{\eps, \pi} (E_z[1(\sgn(\hat{\beta}_j)  \neq  \sgn(\beta_j)) |  X]). 
\]
Note that by Chebyshev's inequality,  
\[
P( \mbox{non-exact variable selection by $\hat{\beta}(Y; X)$} ) \leq d(\hat{\beta} | X).  
\]
So a small Hamming distance guarantees exact variable selection with high probability.

How to characterize precisely  the Hamming distance is a challenging problem.  We approach this by  modeling $X$ as random.  Assume  that  the coordinates of $X$ are iid samples from $N(0,1/n)$: 
\begin{equation} \label{Gaussian2}
X_{ij}  \stackrel{iid}{\sim}  N(0,1/n).  
\end{equation} 
The choice of the variance ensures that most diagonal coordinates of the Gram matrix $C = X^TX$ are approximately $1$. Let $P_X(x)$ denote the joint density of  the coordinates of $X$.    The expected Hamming distance is then 
\[
d^*(\hat{\beta})  =  d^*(\hat{\beta}; \eps, \pi) =  \int   d(\hat{\beta}; \eps, \pi  |    X  = x) P_X(x) dx. 
\]

We adopt an asymptotic framework where we calibrate $p$ and $\eps$ with   
\begin{equation} \label{Gaussian3}
p = n^{1/\theta}, \qquad   \eps_n = n^{(1 - \vartheta)/
\theta} \equiv p^{1 - \vartheta}, \qquad 0 < \theta,   \vartheta < 1. 
\end{equation} 
This models a situation where $p \gg n$ and the vector $\beta$ gets increasingly sparse as  $n$ grows. 
Also,  we assume $\pi_n$ in (\ref{Definebeta})  is a point mass  
\begin{equation} \label{Gaussian4}
\pi_n = \nu_{\tau_n}.  
\end{equation} 
Despite its seemingly idealistic, the model was found to be  subtle and rich in theory (e.g.  \cite{CJL, AoS, RoySoc, PNAS1, PNAS2, Rice}). 
In addition, compare two experiments, in one of them   $\pi_n = \nu_{\tau_n}$, and in the other the support of $\pi_n$ is contained in $[\tau_n, \infty)$.  Since the second model is easier for inference than the first one, the optimal Hamming distance for the first one gives an upper bound for that for the second one.    

With $\eps_n$ calibrated as above, the most interesting range for $\tau_n$ is $O(\sqrt{2 \log p})$ \cite{AoS}:  when $\tau_n \gg \sqrt{2 \log p}$,   exact variable selection can be easily achieved by either the lasso or marginal regression. When $\tau_n \ll \sqrt{2 \log p}$,  no variable selection procedure can achieve exact variable selection.  In light of this, we calibrate 
\begin{equation} \label{Gaussian5}  
\tau_n =    \sqrt{2 (r/\theta)  \log n} \equiv \sqrt{2 r  \log p},   \qquad r > 0. 
\end{equation} 
With these calibrations, we can rewrite  
\[
d_n^*(\hat{\beta}; \eps, \pi) = d_n^*(\hat{\beta}; \eps_n, \tau_n). 
\]

\bed
Denote $L(n)$ by a  multi-log term which satisfies that   $ \lim_{n \goto \infty} (L(n) \cdot   n^{\delta})  =  \infty$ and that $\lim_{n \goto \infty}  (L(n) \cdot n^{-\delta})  = 0$ for any  $\delta  > 0$. 
\eed
We are now ready to spell out the main results.  Define 
\[
\rho(\vartheta) =  (1 + \sqrt{1 - \vartheta})^2, \qquad 0 < \vartheta < 1.  
\]
The following theorem is proved in the appendix, which gives the lower bound for the Hamming distance. 
\begin{theorem} \label{thm:LB}
Fix $\vartheta \in (0,1)$, $\theta > 0$, and $r > 0$ such that  $\theta > 2 ( 1  - \vartheta)$.   Consider a sequence of regression models as in (\ref{Gaussian0})-(\ref{Gaussian5}).    As $n \goto \infty$, for any variable selection procedure $\hat{\beta}^{(n)}$, 
\[
d_n^*(\hat{\beta}^{(n)}; \eps_n, \tau_n)   \geq   \left\{ \begin{array}{ll}
L(n) p^{1-\frac{(\vartheta + r)^2}{4 r}},  &\qquad r \geq  \vartheta,  \\
(1 + o(1)) \cdot p^{1 -\vartheta},  &\qquad  0 < r < \vartheta.  
\end{array}  
\right. 
\]
\end{theorem}

\begin{figure}
\includegraphics[width=3in]{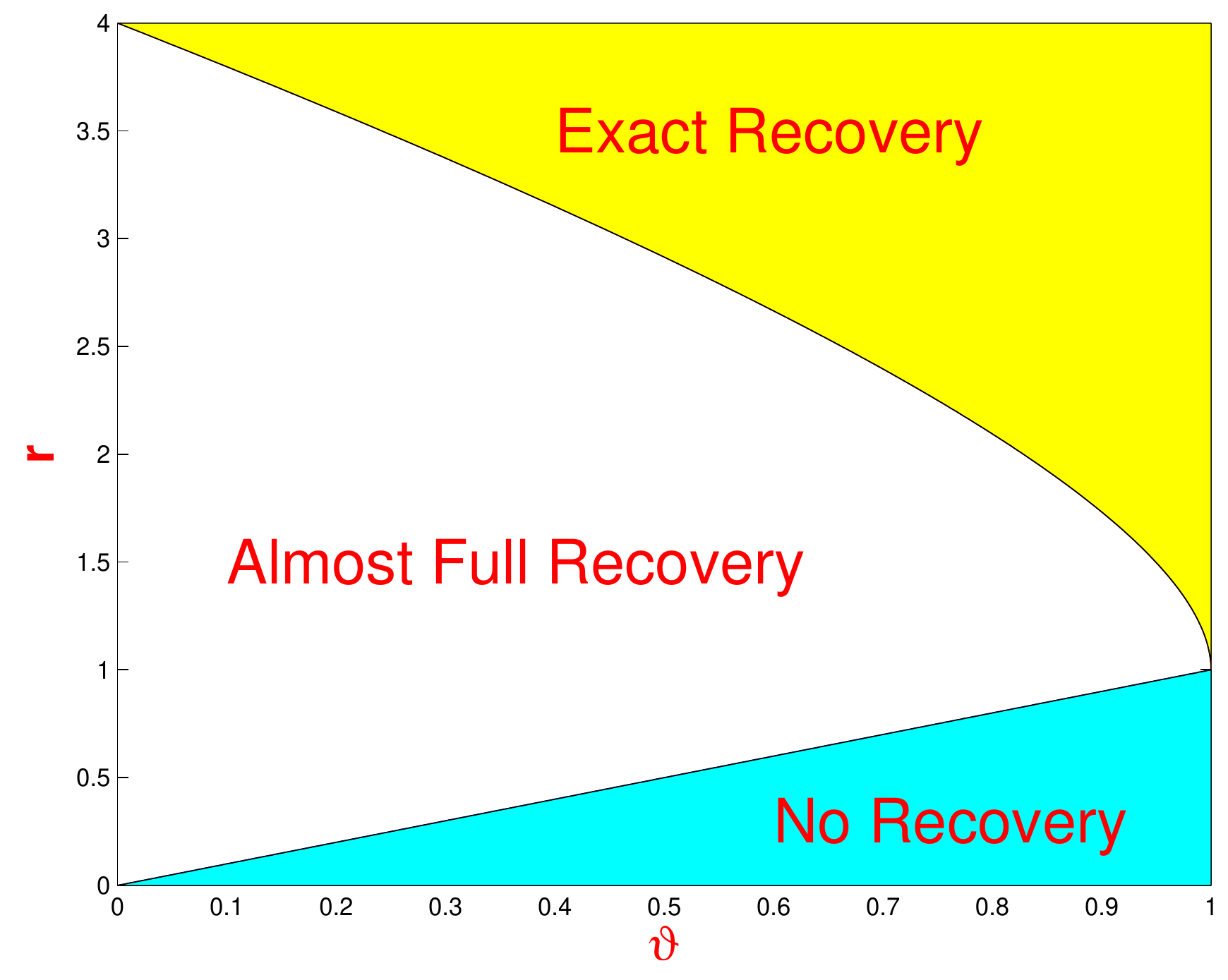}
\caption{There regions as in Section \ref{sec:hamm}. In  Region of Exact Recovery, both the lasso and marginal regression yield exact recovery with high probability. 
In Region of Almost  Full Recovery,   it is impossible to have large probability for exact variable selection, but the Hamming distance of both the lasso and marginal regression $\ll p \eps_n$. In Region of No Recovery,  optimal Hamming distance $\sim p \eps_n$ and all variable selection procedures fail completely.  Displayed is the part of the plane corresponding to  $0 < r < 4$ only.    }
\label{fig::Phase}
\end{figure}

At the same time,  let $\hat{\beta}_{mr}$ be the estimate  of using marginal regression  with   threshold 
\[
 t_n    =  (\frac{\vartheta + r}{2 \sqrt{r}} \wedge \sqrt{r}) \cdot   \sqrt{2 \log p}. 
\]
We have the following theorem. 
\begin{theorem}  \label{thm:mr}
Fix $\vartheta \in (0,1)$, $r > 0$, and $\theta >  (1 - \vartheta)$.  Consider a sequence of regression models as in (\ref{Gaussian0})-(\ref{Gaussian5}).  
As $p \goto \infty$, the Hamming distance of  marginal regression with the threshold $t_n =  (\frac{\vartheta + r}{2 \sqrt{r}} \wedge \sqrt{r}) \cdot   \sqrt{2 \log p^{(n)}}$ satisfies 
\[
d_n^*(\hat{\beta}_{mr}^{(n)}; \eps_n, \tau_n)  \leq   \left\{ \begin{array}{ll}
L(n) p^{1-\frac{(\vartheta + r)^2}{4 r}},  &\qquad r \geq  \vartheta,  \\
(1 + o(1)) \cdot  p^{1 -\vartheta},  &\qquad  0 < r < \vartheta.  
\end{array}  
\right. 
\]
\end{theorem} 

Similarly, choosing the tuning parameter 
$\lambda_n  = 2  (\frac{\vartheta + r}{2 \sqrt{r}} \wedge \sqrt{r})  \sqrt{2 \log p}$ in the lasso, 
we have the following theorem. 
\begin{theorem}  \label{thm:lasso}
Fix $\vartheta \in (0,1)$, $r > 0$, and $\theta >  (1 - \vartheta)$. Consider a sequence of regression models as in (\ref{Gaussian0})-(\ref{Gaussian5}).  
As $p \goto \infty$, the Hamming distance of the lasso with the tuning parameter $\lambda_n = 2(\frac{\vartheta + r}{2 \sqrt{r}} \wedge \sqrt{r}) \cdot   \sqrt{2 \log p^{(n)}}$ satisfies 
\[
d_n^*(\hat{\beta}_{lasso}^{(n)}; \eps_n, \tau_n)  \leq   \left\{ \begin{array}{ll}
L(n) p^{1-\frac{(\vartheta + r)^2}{4 r}},  &\qquad r \geq  \vartheta,  \\
(1 + o(1)) \cdot p^{1 -\vartheta},  &\qquad  0 < r < \vartheta.  
\end{array}  
\right. 
\]
\end{theorem} 
The proofs of  Theorems \ref{thm:mr}-\ref{thm:lasso} are routine   and we omit them.

Theorems \ref{thm:LB}-\ref{thm:lasso}  say  that in the $\vartheta$-$r$ plane, we have three  different regions, as displayed in Figure \ref{fig::Phase}. 
\begin{itemize} 
\item Region I ({\it   Exact Reovery}):    $0 < \vartheta < 1$ and $r > \rho(\vartheta)$. 
\item Region II ({\it  Almost Full Recovery}):     $0 < \vartheta < 1$ and $\vartheta < r < \rho(\vartheta)$. 
\item Region III ({\it No Recovery}):    $0 < \vartheta < 1$ and $0 < r < \vartheta$.  
\end{itemize}
In the Region of Exact Recovery, the Hamming distance for both marginal regression and the lasso    are  algebraically small. Therefore, except for a probability that is algebraically small, both marginal regression and the lasso 
give exact recovery. 
 
In the Region of Almost Full Recovery, both the Hamming distance of marginal regression and the lasso are  much smaller than the number of relevant variables (which $\approx p \eps_n$).  Therefore, almost all relevant variables have been recovered. Note also that  the number of misclassified  irrelevant variables  is comparably much smaller than $p \eps_n$.  In this region, the optimal Hamming is algebraically large, so  for any variable selection procedure, the probability of  exact recovery is algebraically small.  

In the Region of No Recovery,   the Hamming distance   $\sim p \eps_n$.  In this region,  asymptotically,  it is impossible to distinguish relevant variables from irrelevant variables, and any variable selection procedure  fails completely.

The results improve on those  by  Wainwright  \cite{Wainwright:2006}. It was shown in \cite{Wainwright:2006} that   
 there are constants $c_2 > c_1 > 0$ such that in the region of $\{0 < \vartheta < 1,  r > c_2\}$,  
 the lasso yields exact variable selection with overwhelming probability, and that in the region of $\{0 < \vartheta < 1, r < c_2\}$,  no procedure could yield exact variable selection.   
 Our results not only provide the exact rate of the Hamming distance, but also  tighten the constants $c_1$ and $c_2$ so that 
 $c_1 = c_2 = (1 + \sqrt{1 - \vartheta})^2$.

\section{Simulations and Examples}  
\label{sec:simul} 

\begin{figure}
\includegraphics[height=2in, width=2.5in]{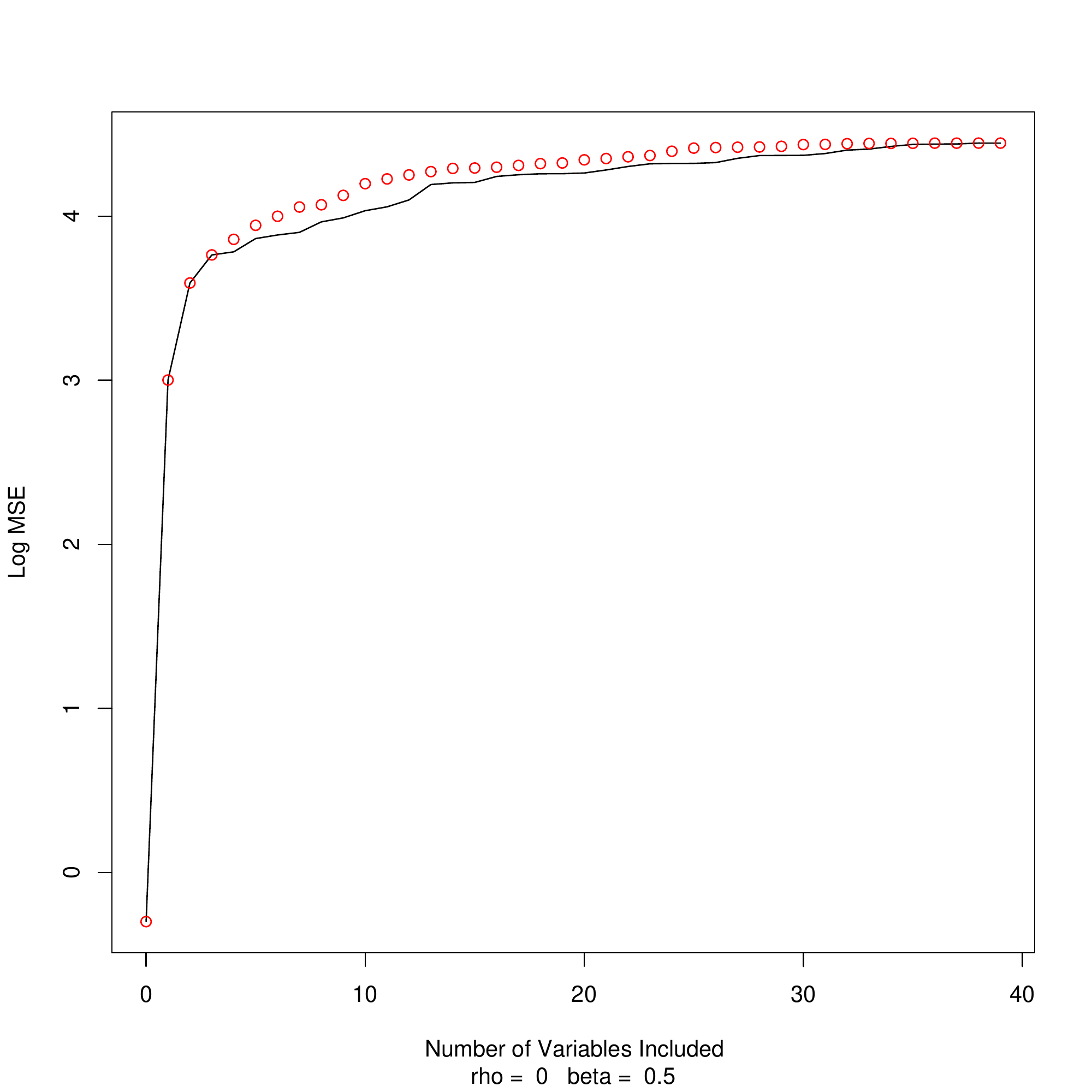}\includegraphics[height=2in, width=2.5in]{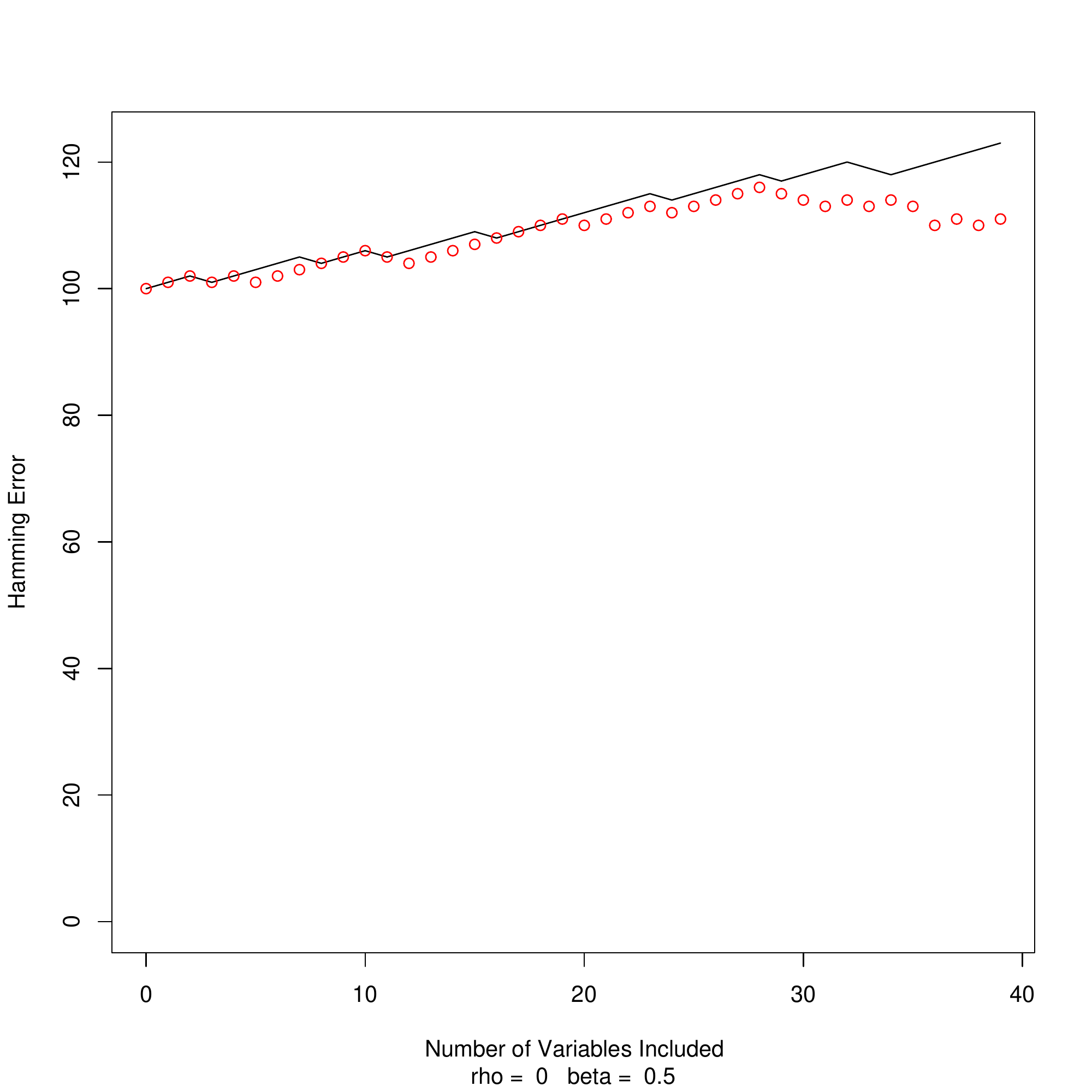}\\
\includegraphics[height=2in, width=2.5in]{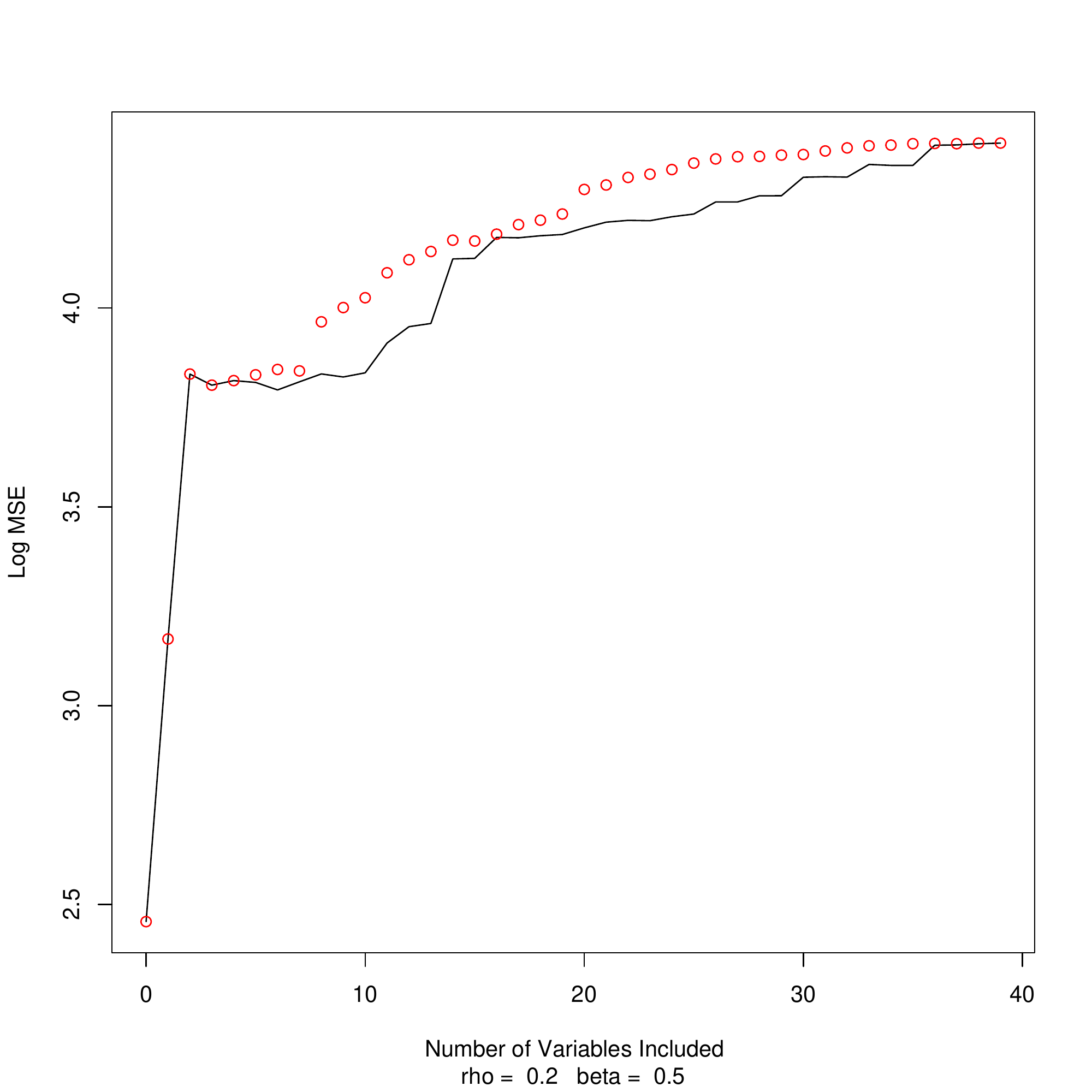}\includegraphics[height=2in, width=2.5in]{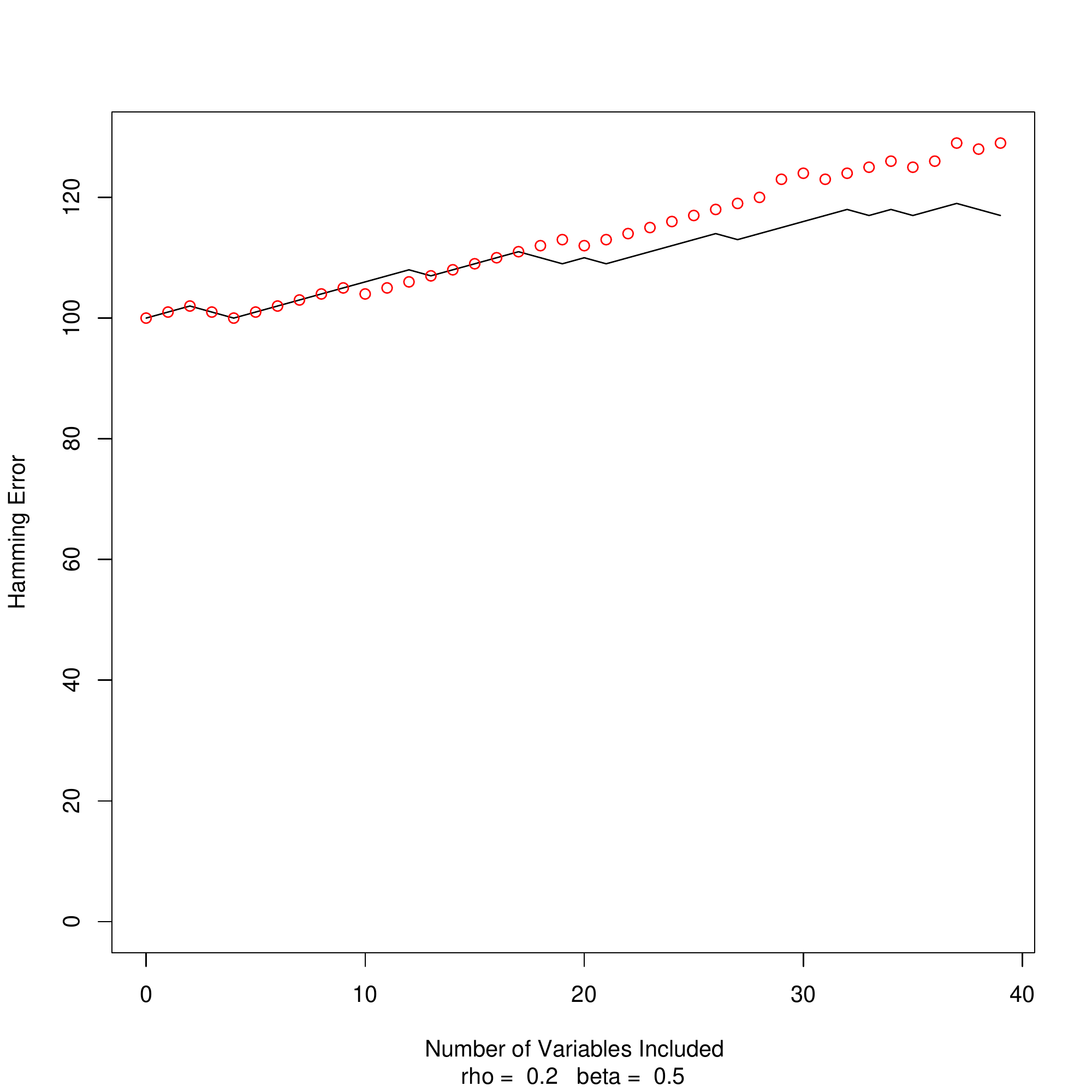}\\
\includegraphics[height=2in, width=2.5in]{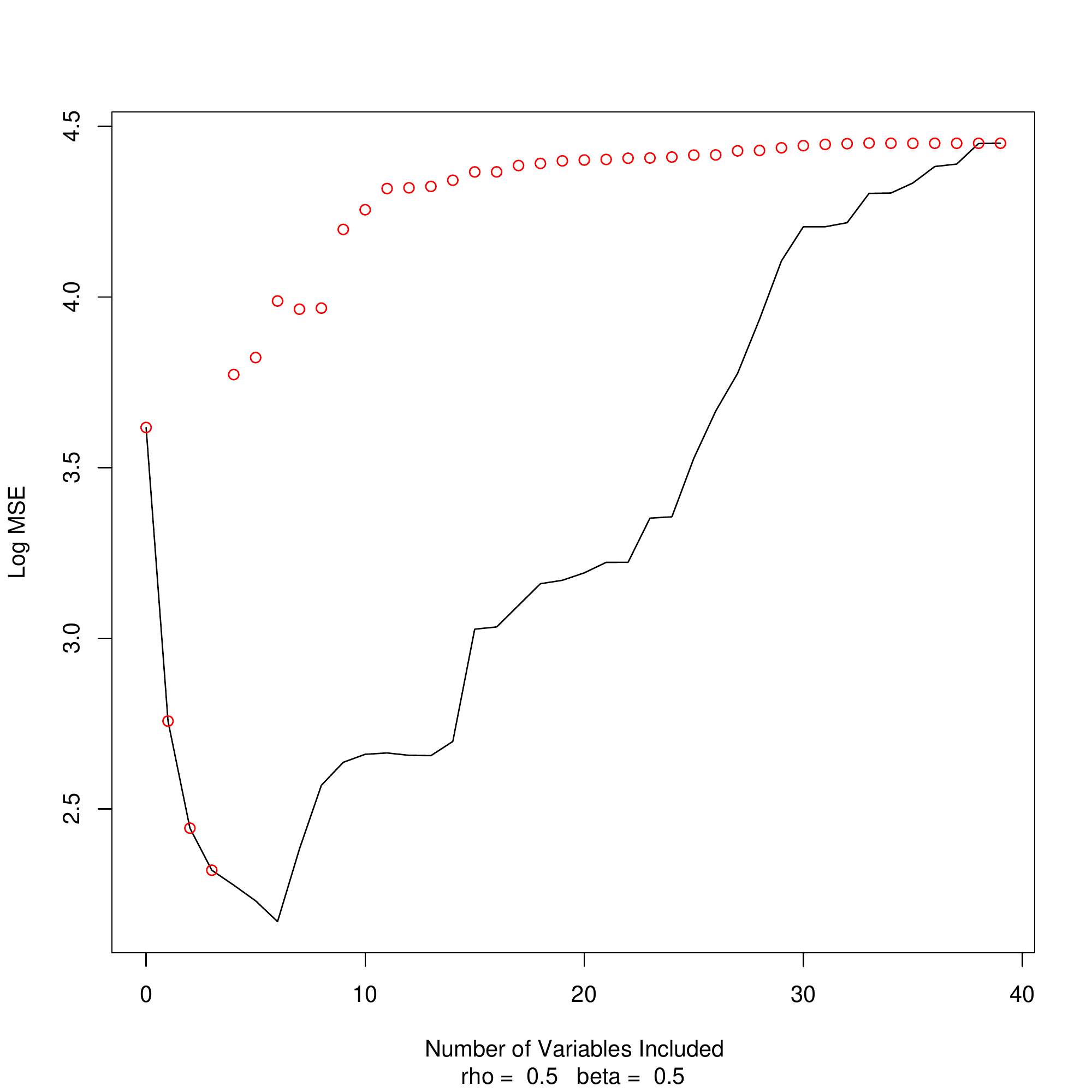}\includegraphics[height=2in, width=2.5in]{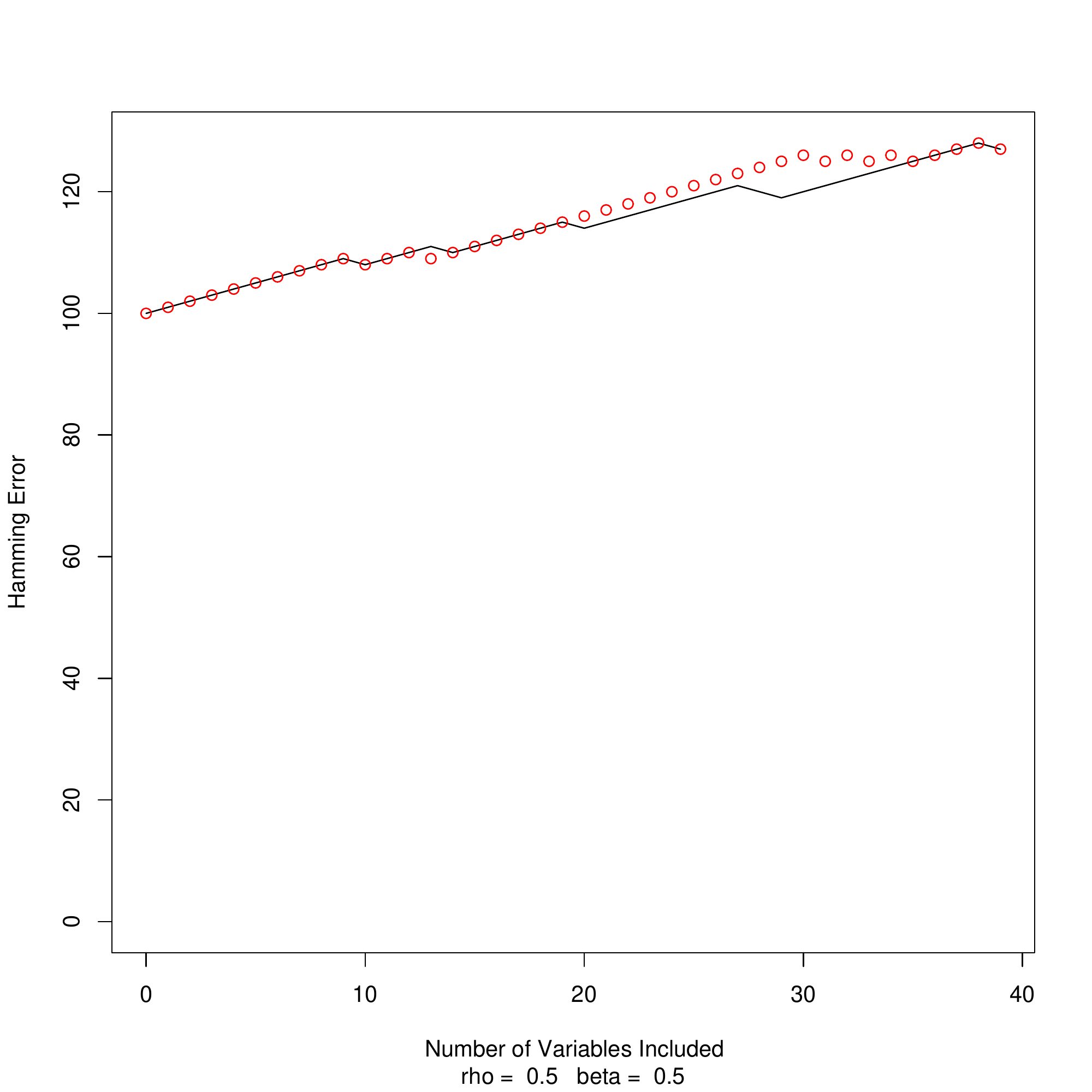}\\
\includegraphics[height=2in, width=2.5in]{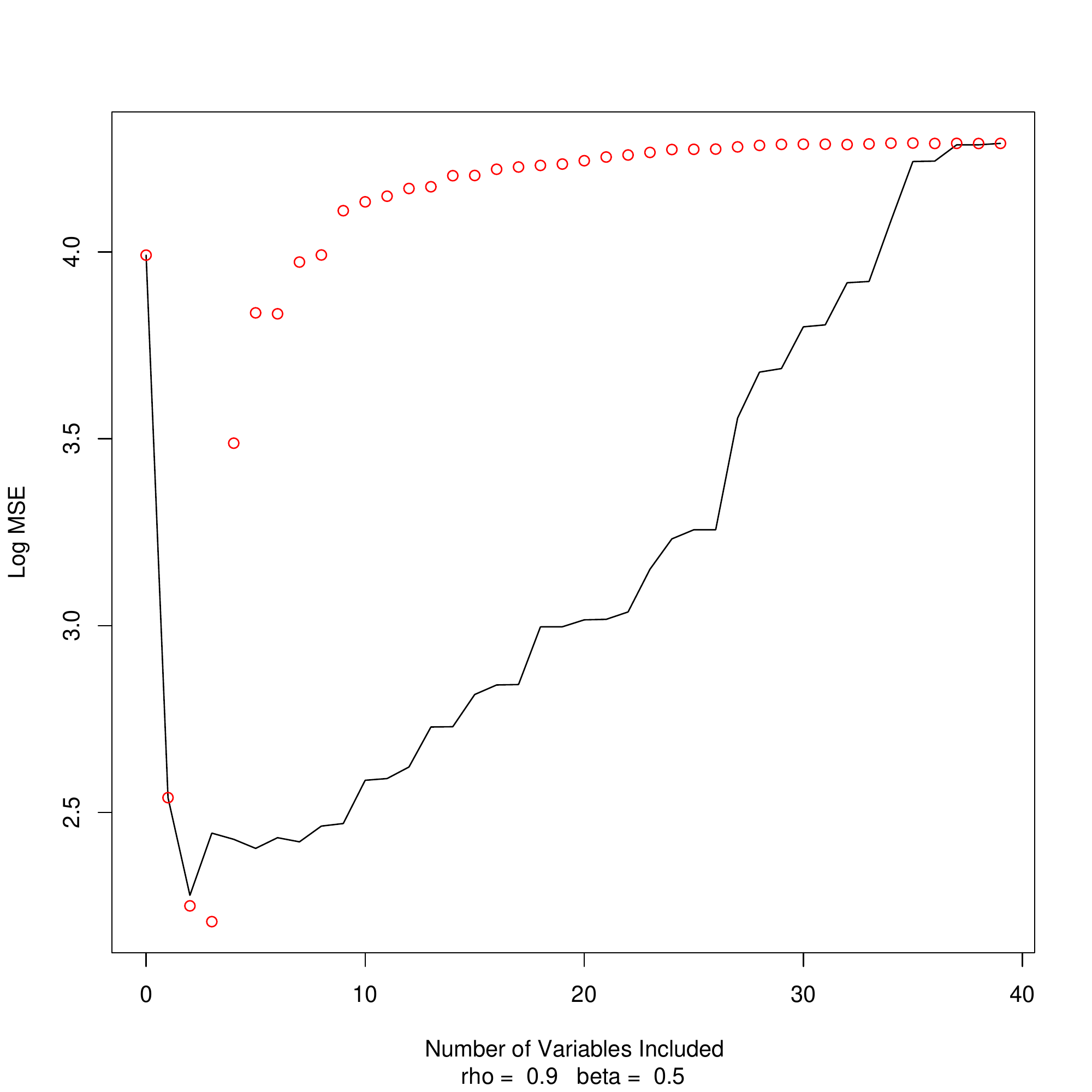}\includegraphics[height=2in, width=2.5in]{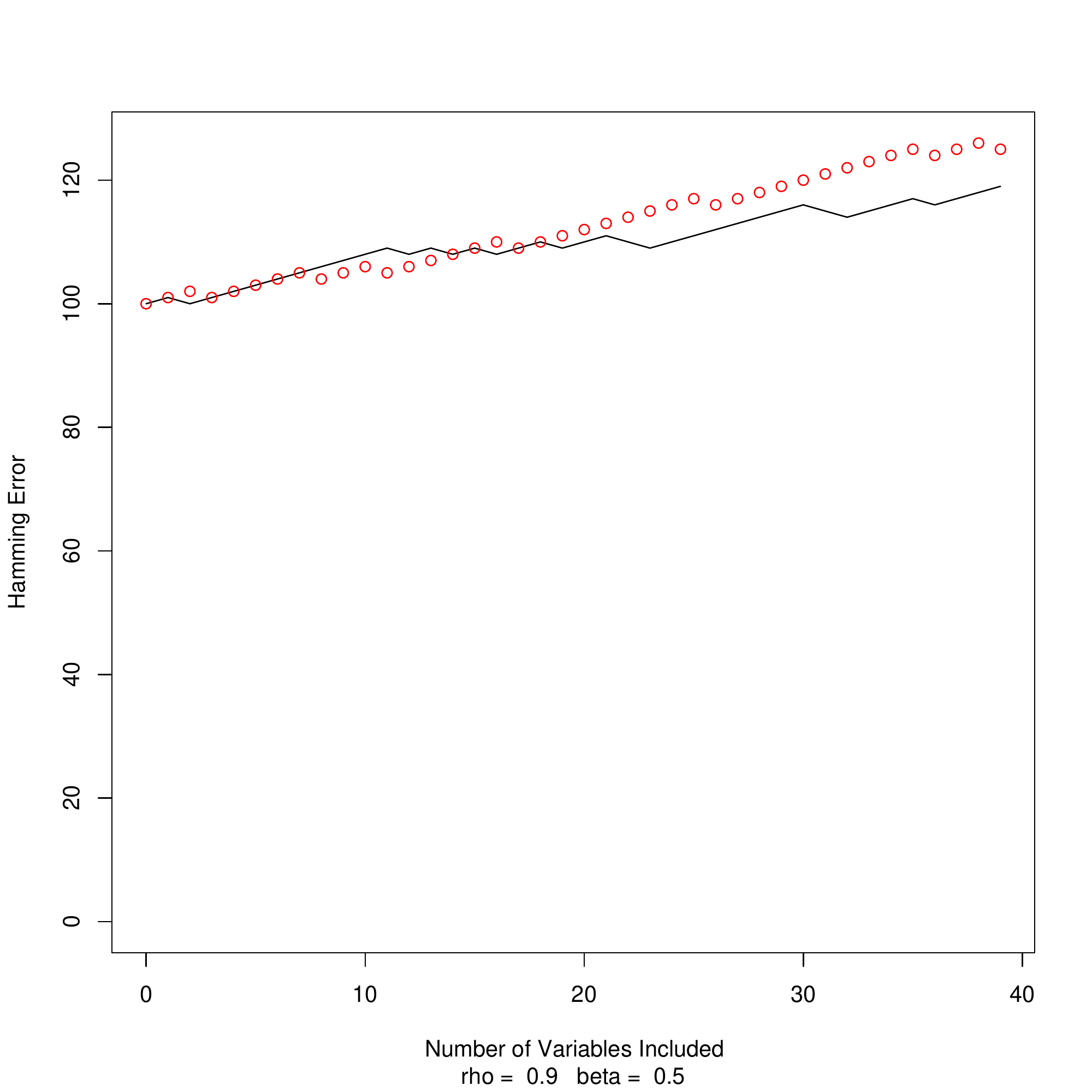}\\
\caption{Prediction error (left column) and Hamming error (right column) for the lasso (red circle) and marginal regression (solid line). The  x-axis displays the number of variables included.  Row 1-4: $\rho = 0, 0.2, 0.5, 0.9$. Nonzero $\beta_i$: $0.5$. Results are based on one replication.    }
\label{fig::ex1}
\end{figure}

\begin{figure}
\includegraphics[height=2in, width= 2.5in]{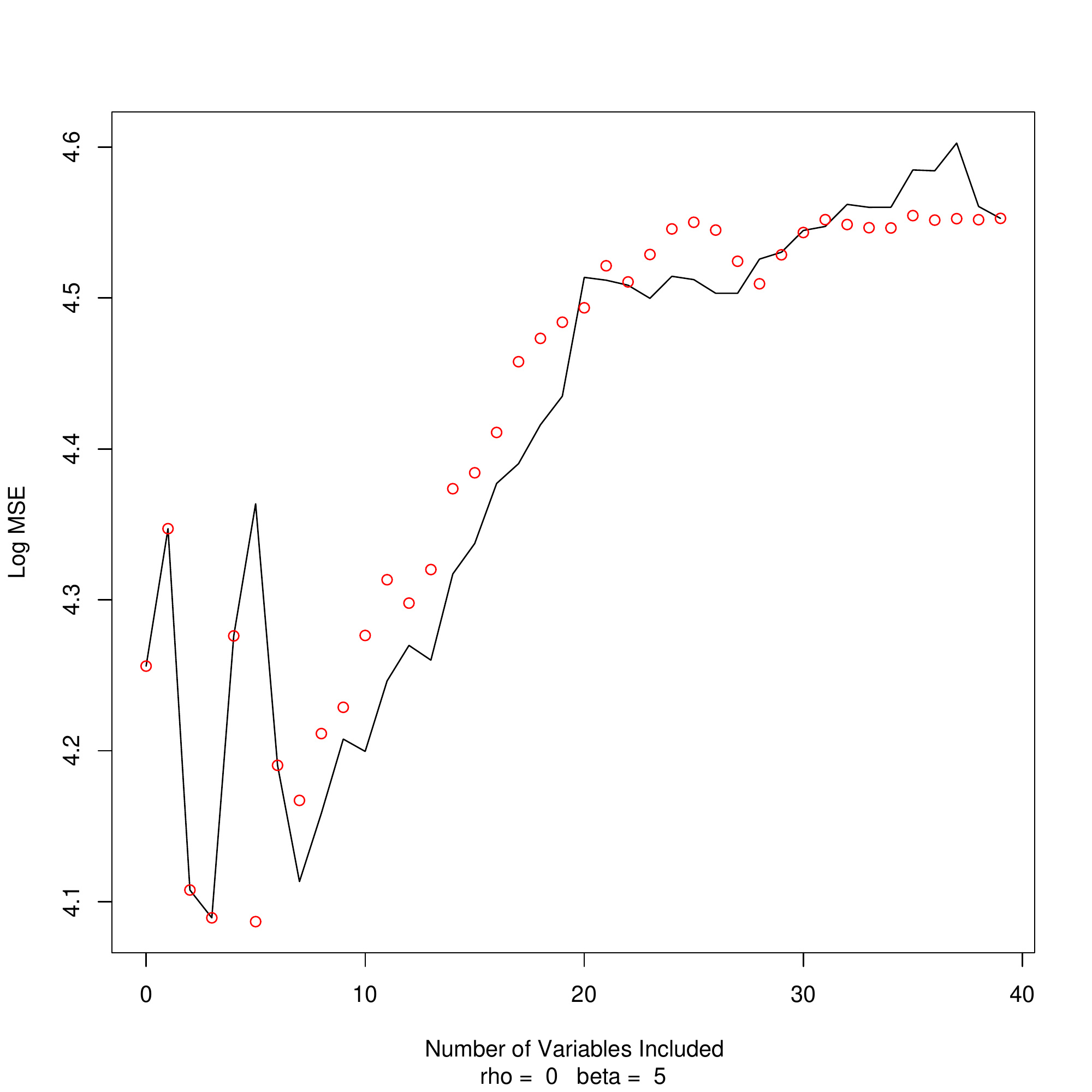}\includegraphics[height=2in, width= 2.5in]{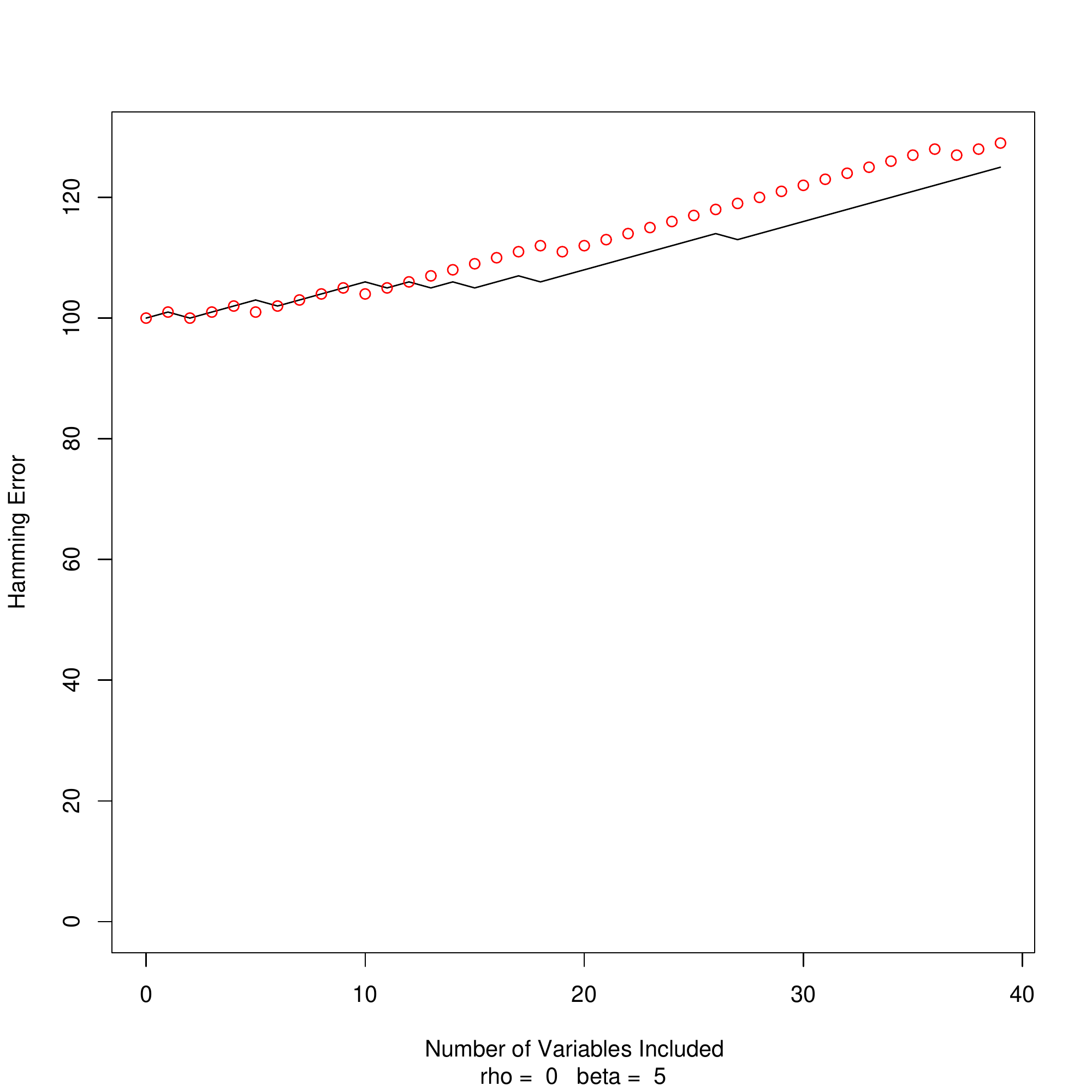}\\
\includegraphics[height=2in, width= 2.5in]{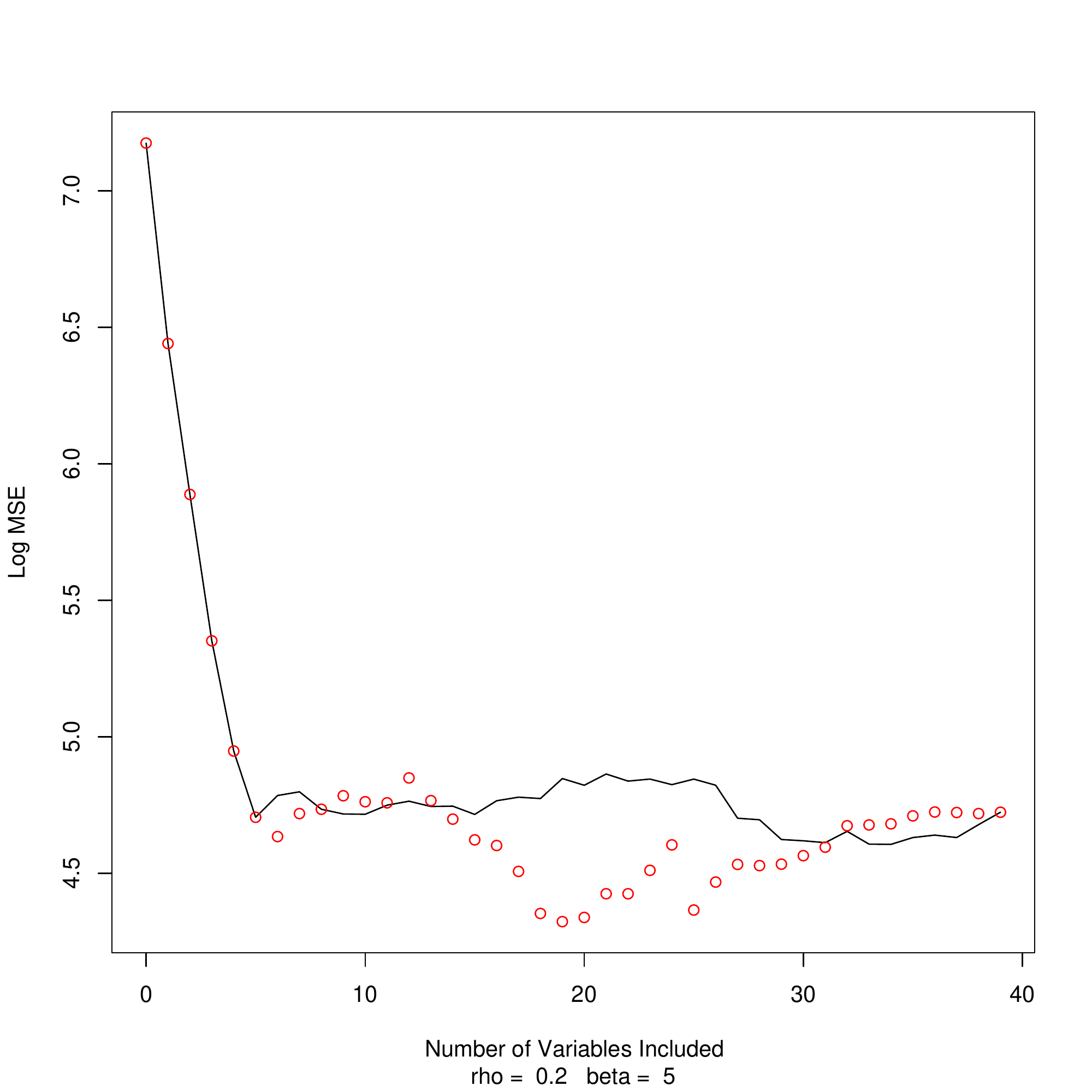}\includegraphics[height=2in, width= 2.5in]{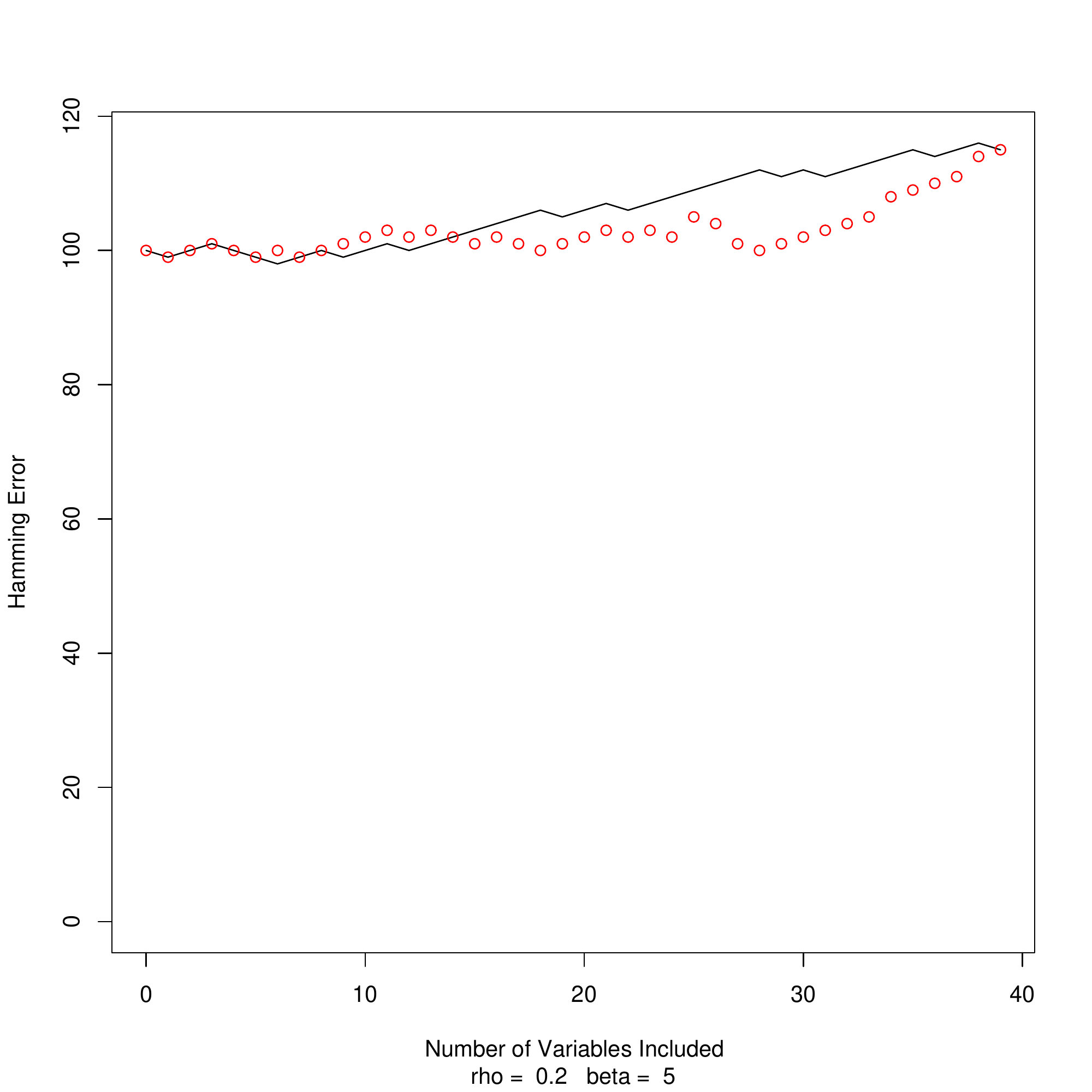}\\
\includegraphics[height=2in, width= 2.5in]{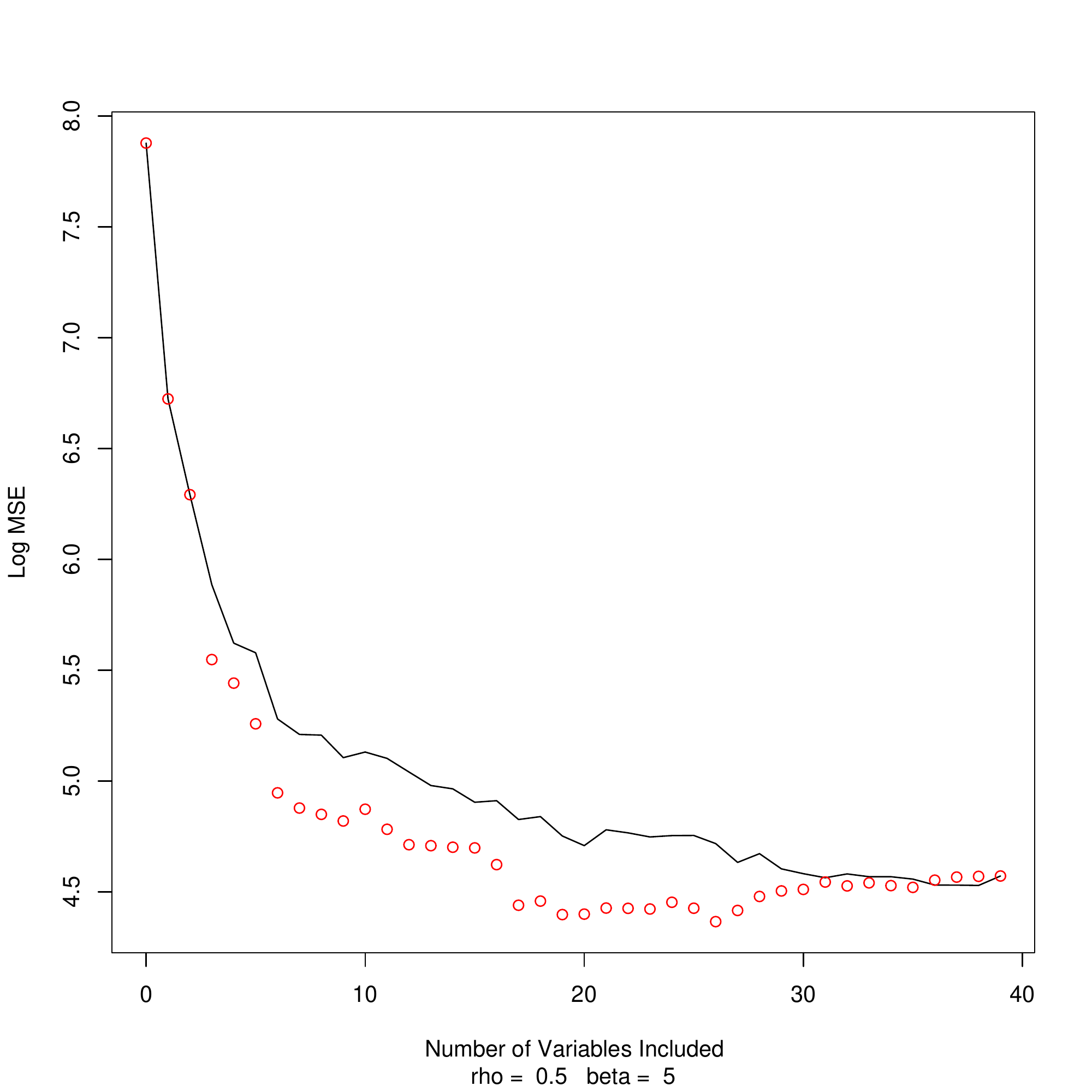}\includegraphics[height=2in, width= 2.5in]{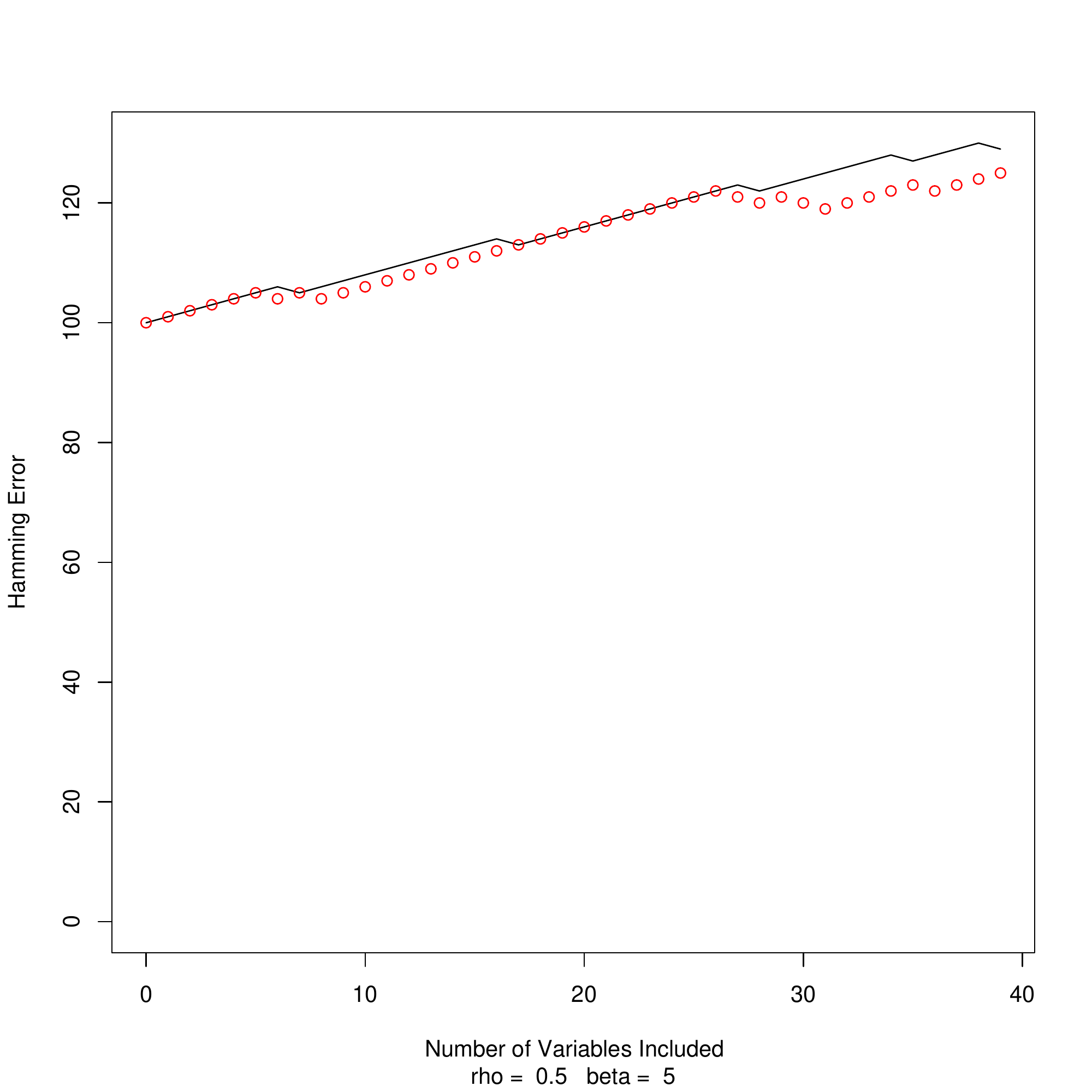}\\
\includegraphics[height=2in, width= 2.5in]{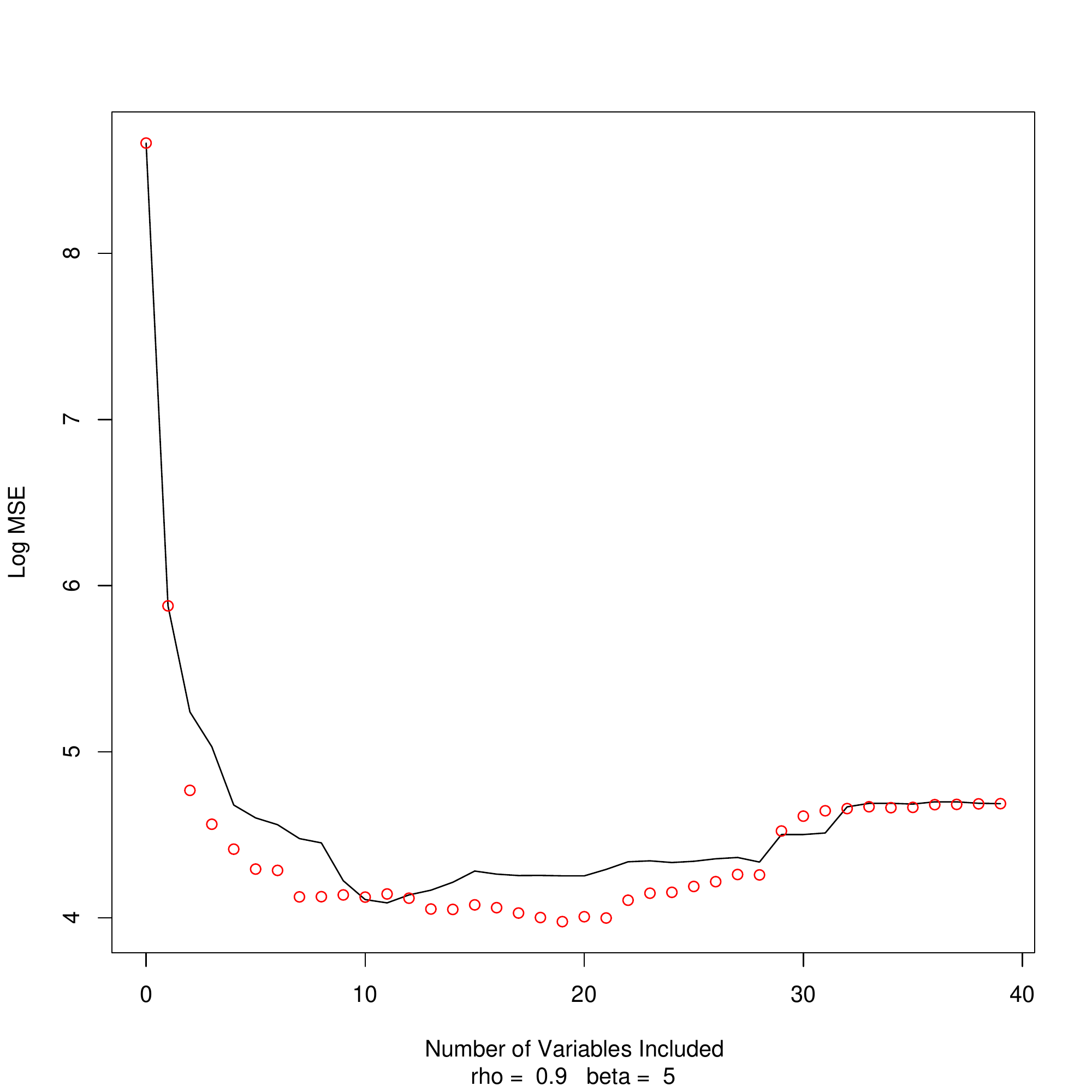}\includegraphics[height=2in, width= 2.5in]{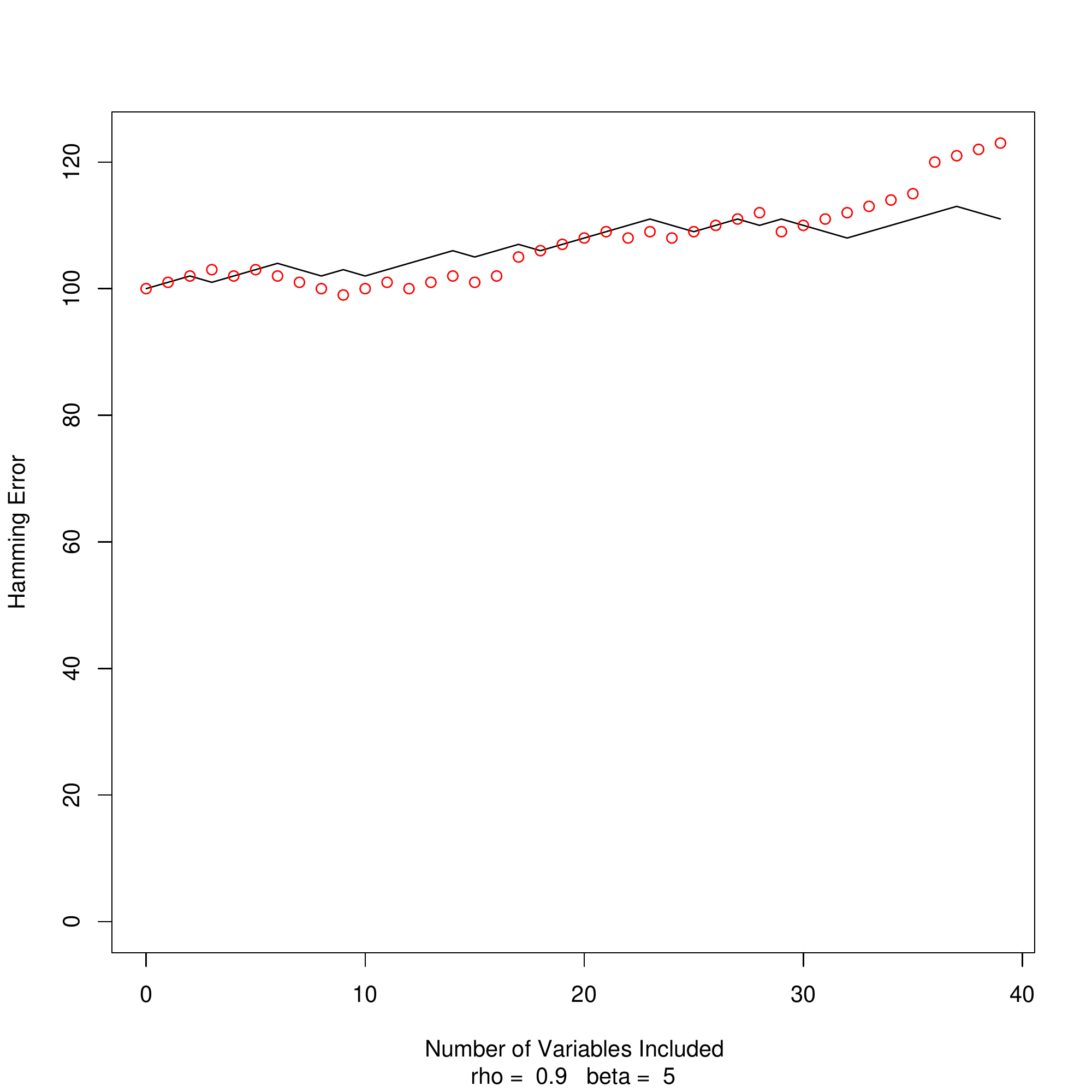}\\
\caption{Same as in Figure \ref{fig::ex1}, but all nonzero $\beta_i$ equal $5$.  }
\label{fig::ex2}
\end{figure}

\begin{figure}
\includegraphics[height=2in, width= 2.5in]{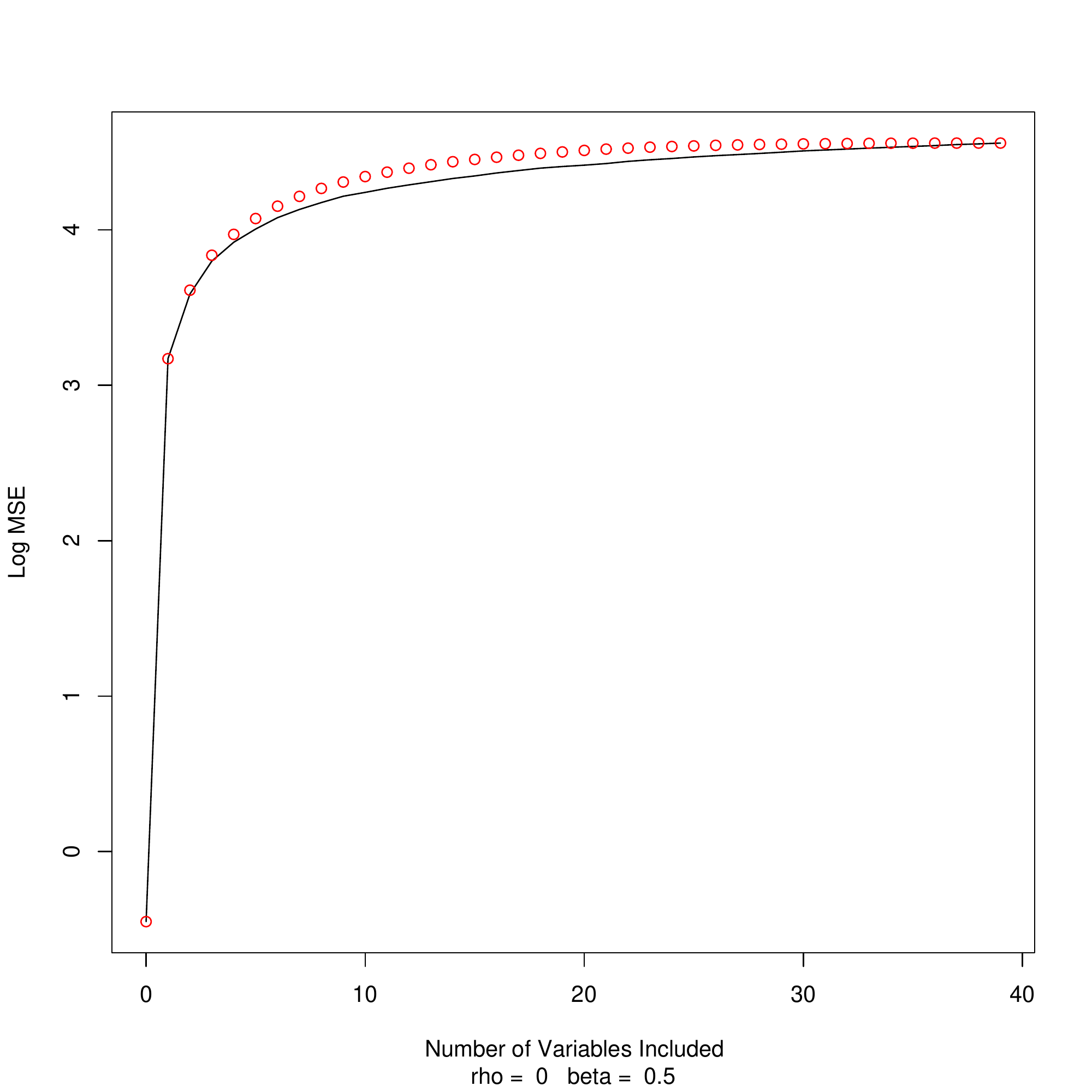}\includegraphics[height=2in, width=2.5in]{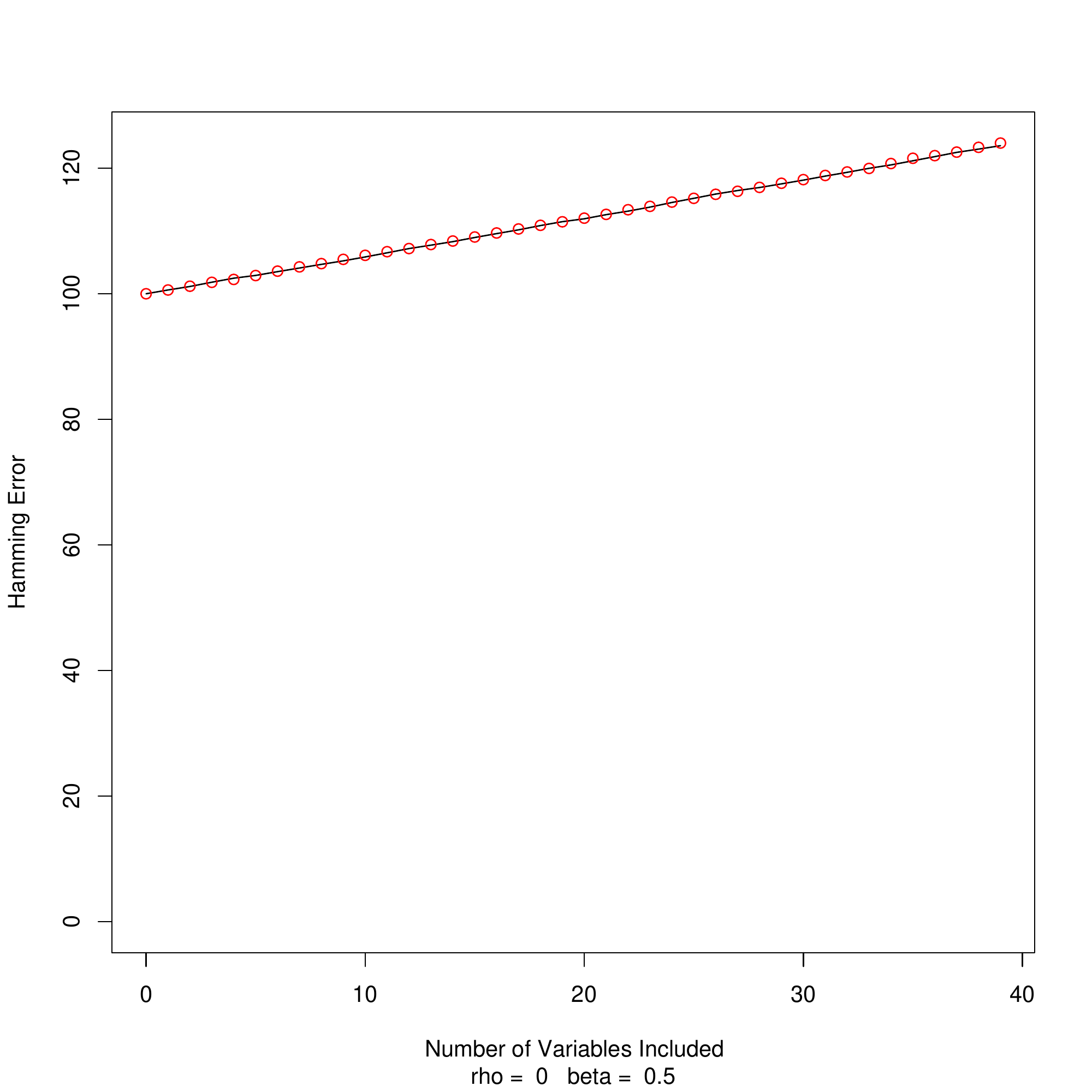}\\
\includegraphics[height=2in, width= 2.5in]{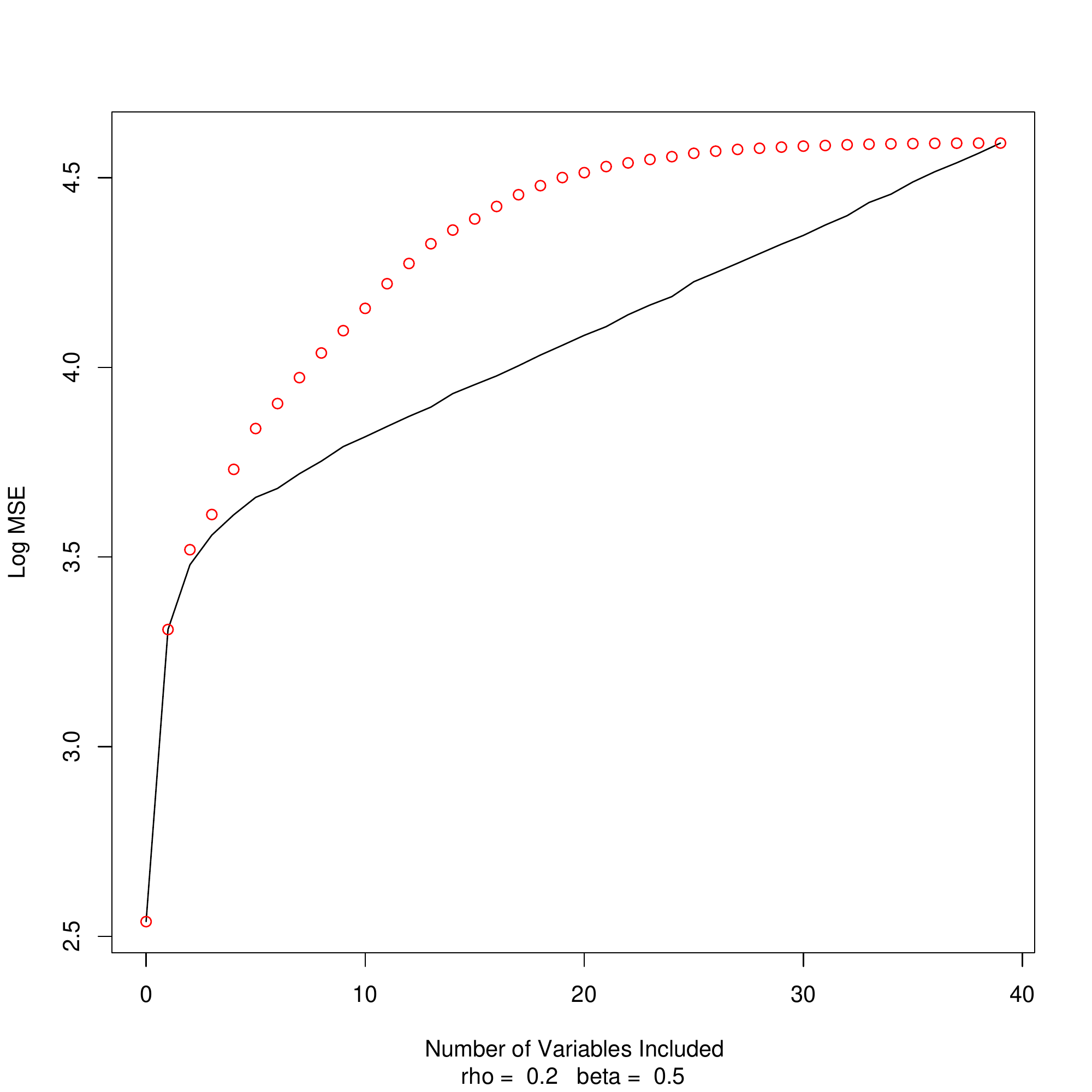}\includegraphics[height=2 in, width=2.5in]{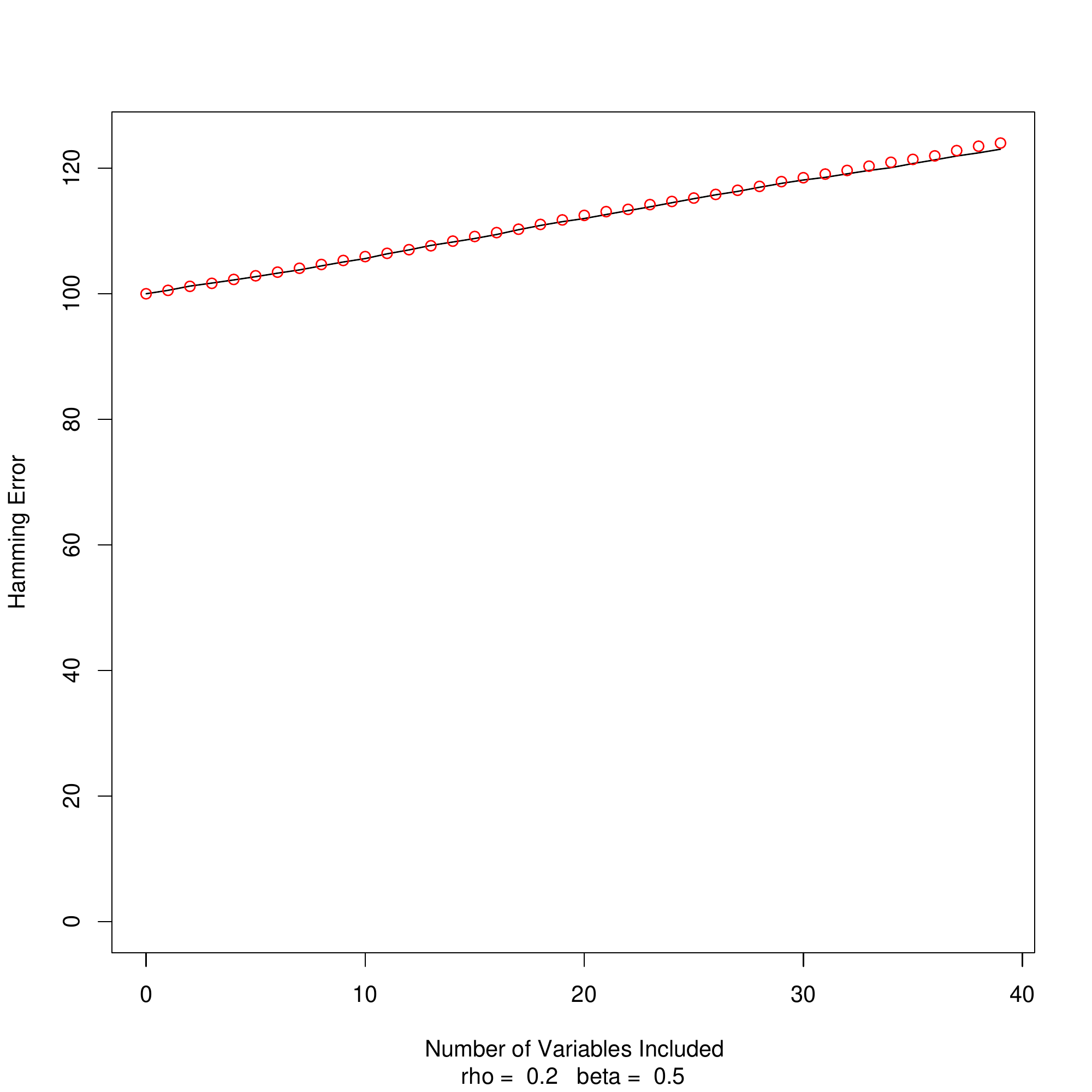}\\
\includegraphics[height=2in, width= 2.5in]{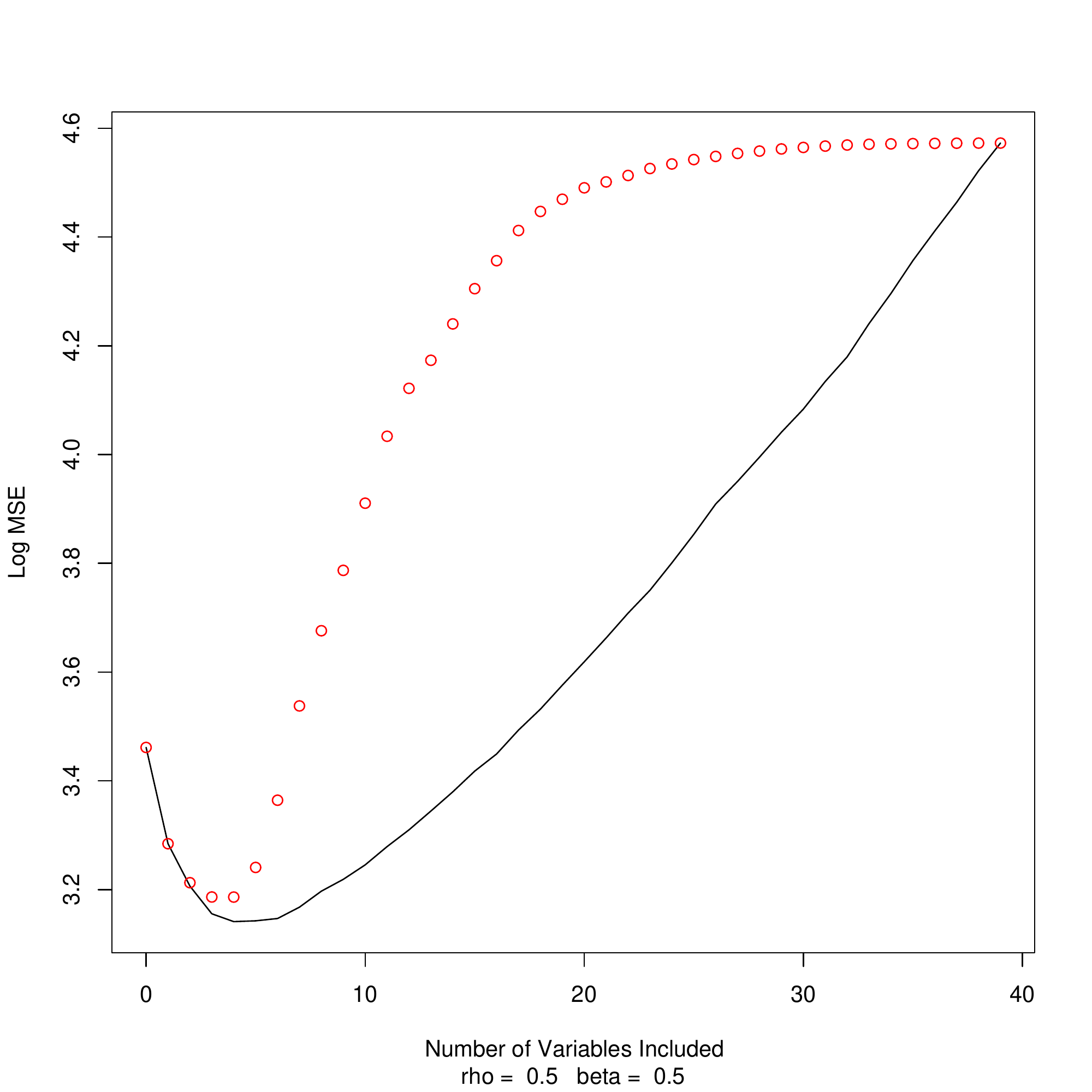}\includegraphics[height=2in, width=2.5in]{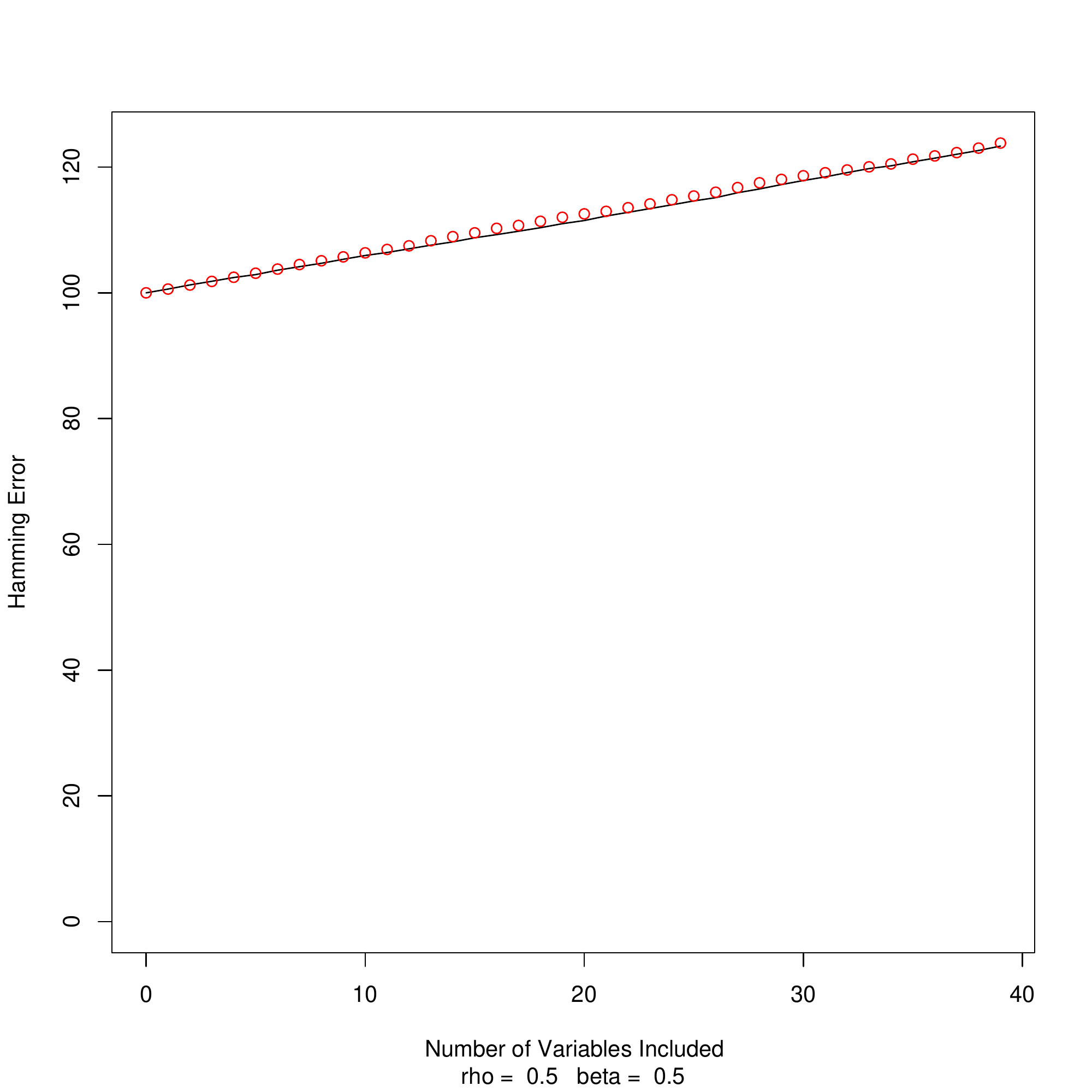}\\
\includegraphics[height=2in, width= 2.5in]{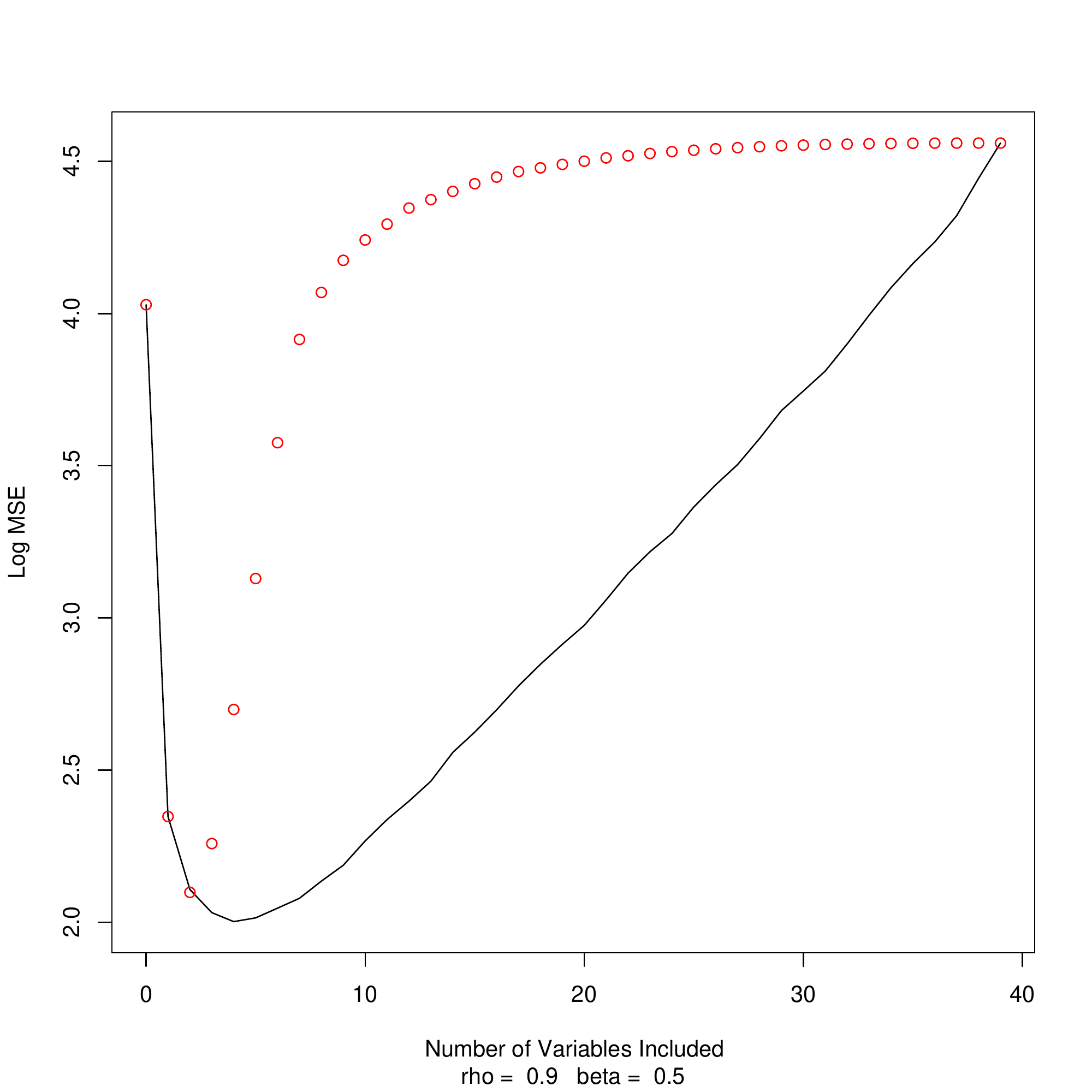}\includegraphics[height=2in, width=2.5in]{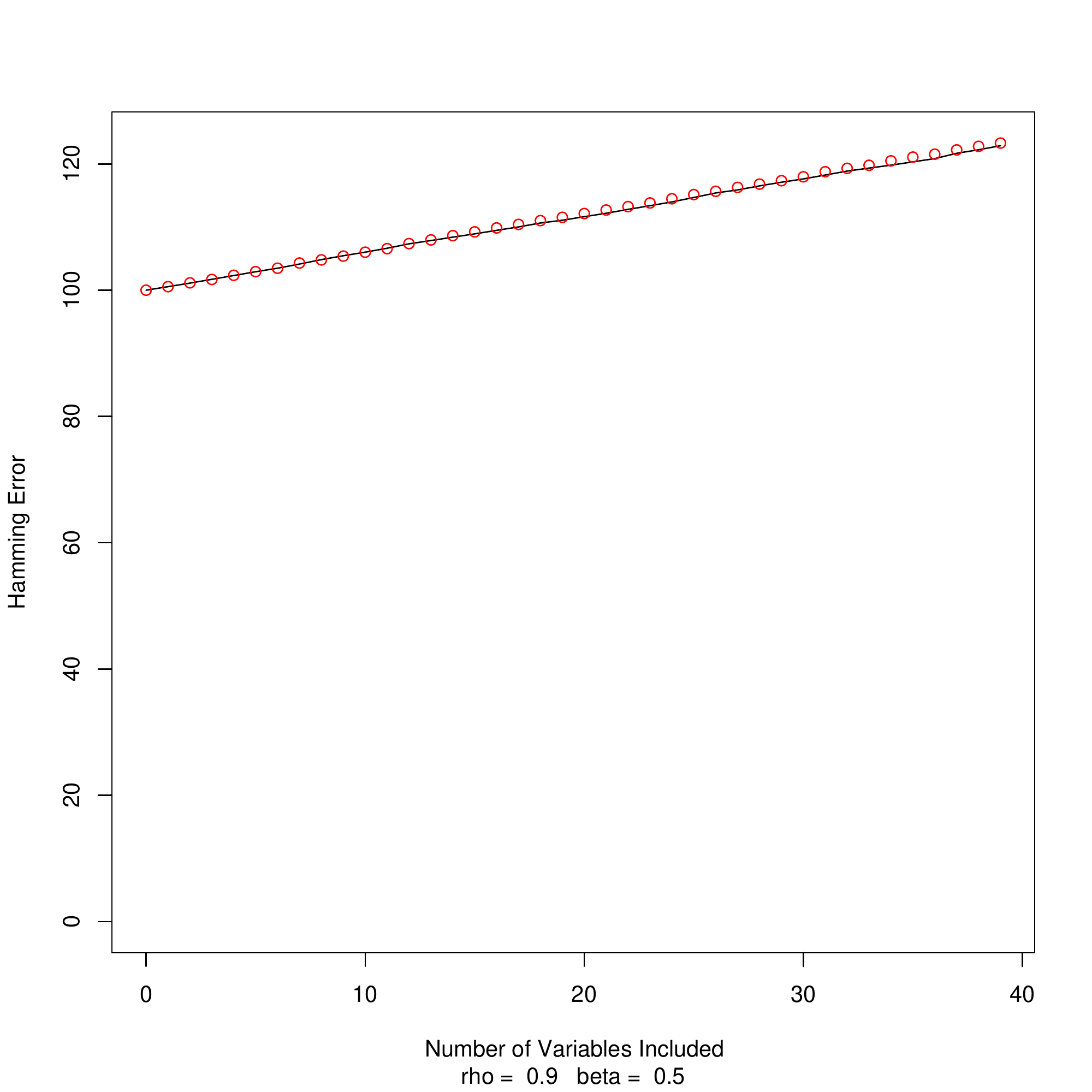}\\
\caption{Same as in Figure \ref{fig::ex1}, but displayed are  the average errors across $100$ replications.}
\label{fig::sim1}
\end{figure}

\begin{figure}
\includegraphics[height=2in, width= 2.5in]{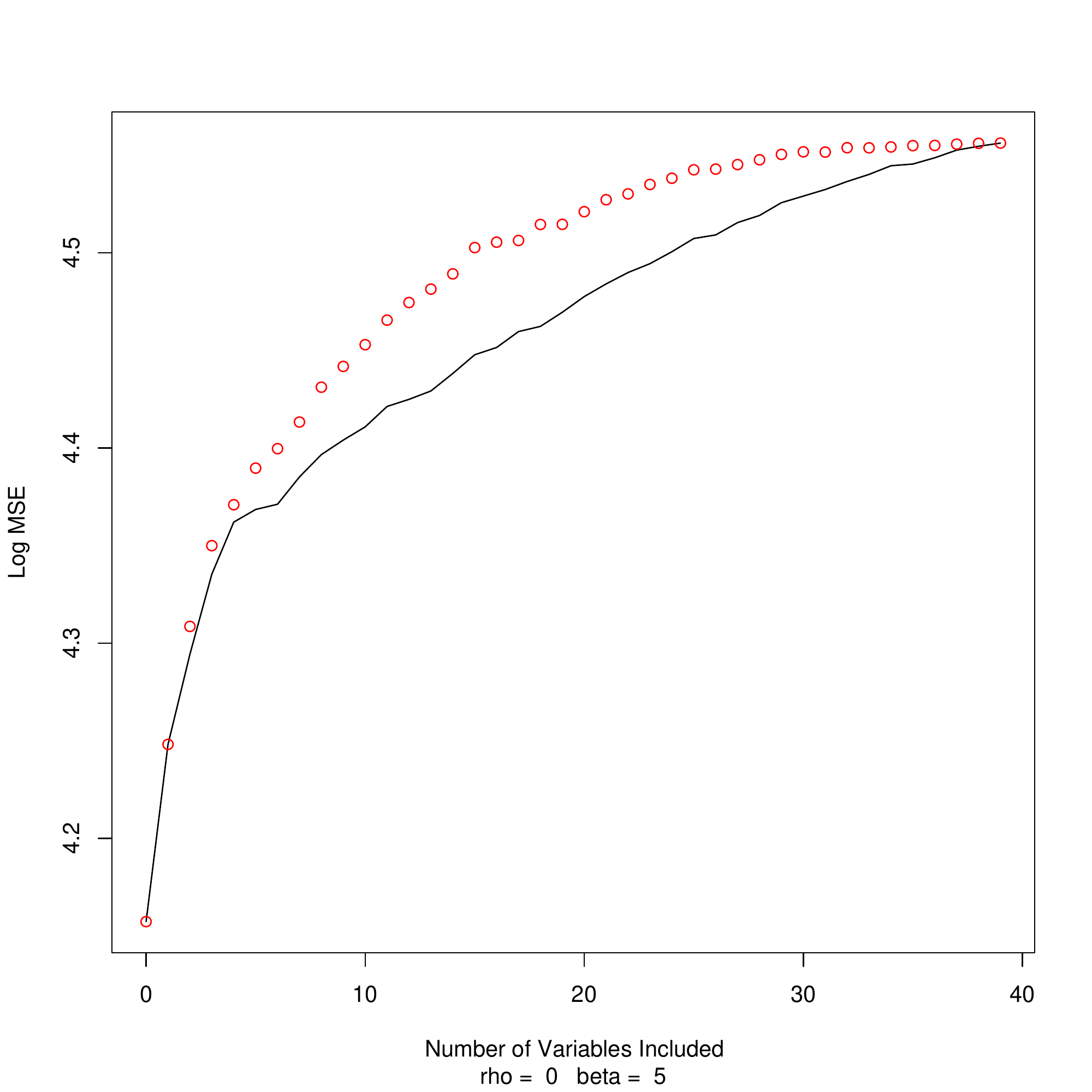}\includegraphics[height=2in, width= 2.5in]{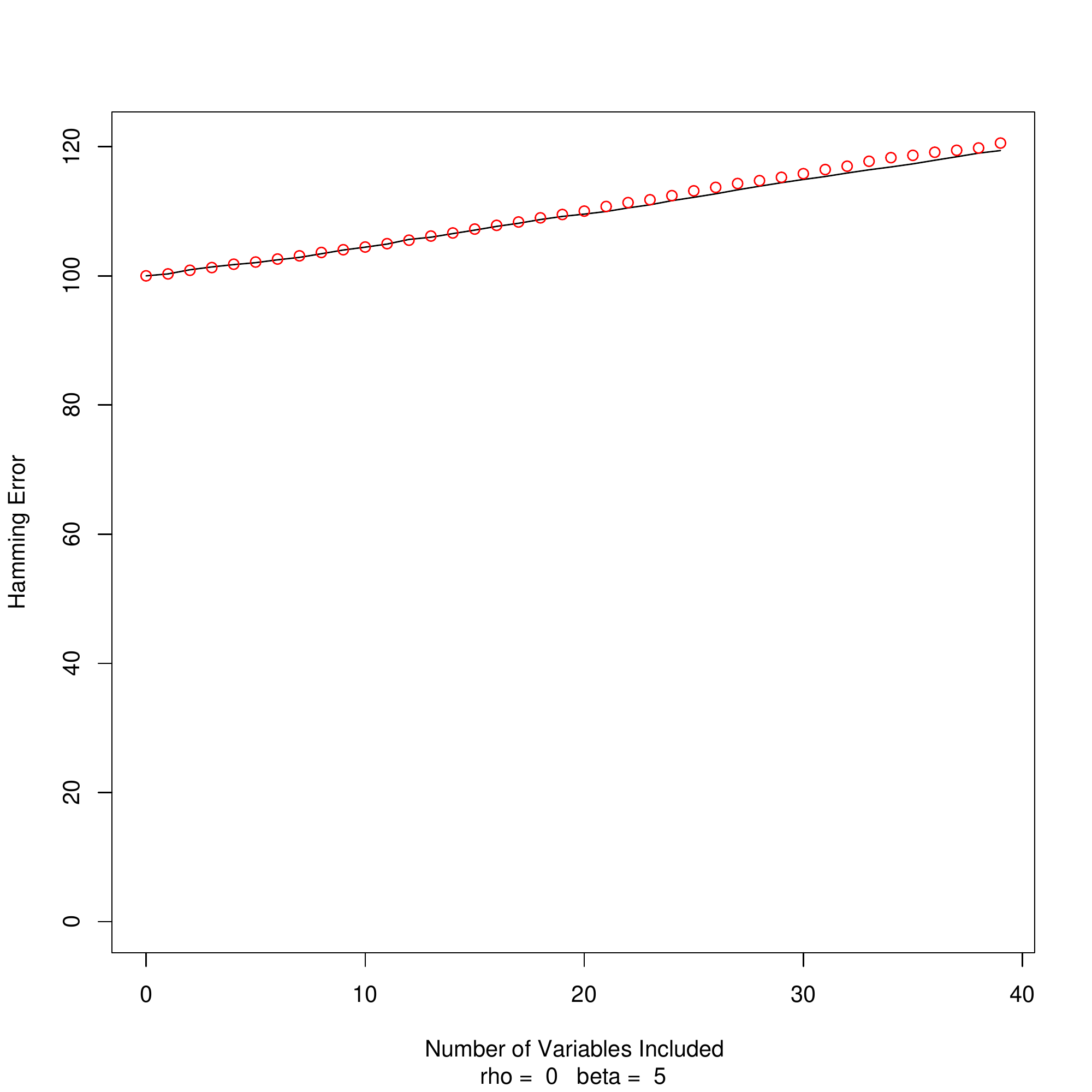}\\
\includegraphics[height=2in, width= 2.5in]{plot6_Part1.pdf}\includegraphics[height=2 in, width=2.5in]{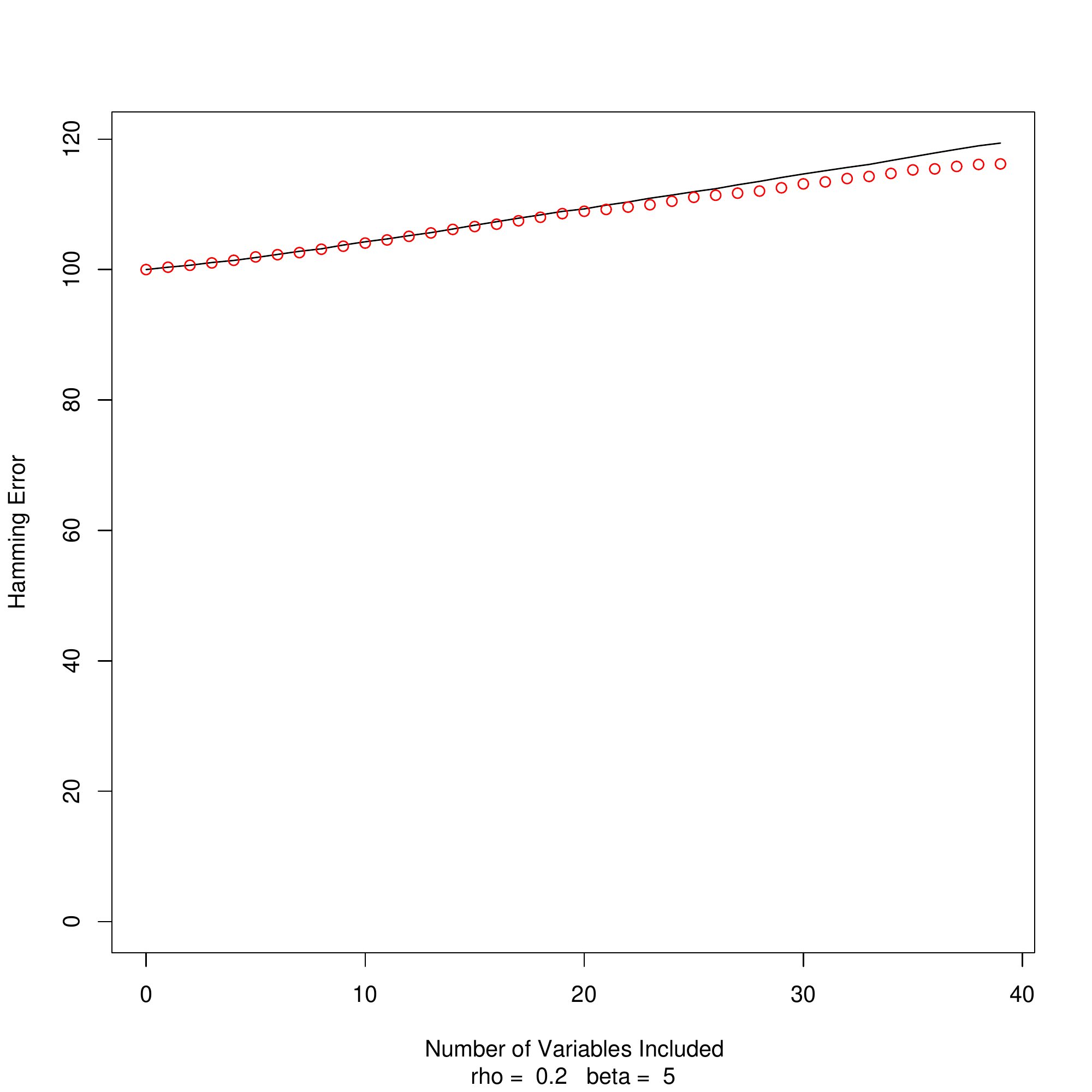}\\
\includegraphics[height=2in, width= 2.5in]{plot7_Part1.pdf}\includegraphics[height=2in, width= 2.5in]{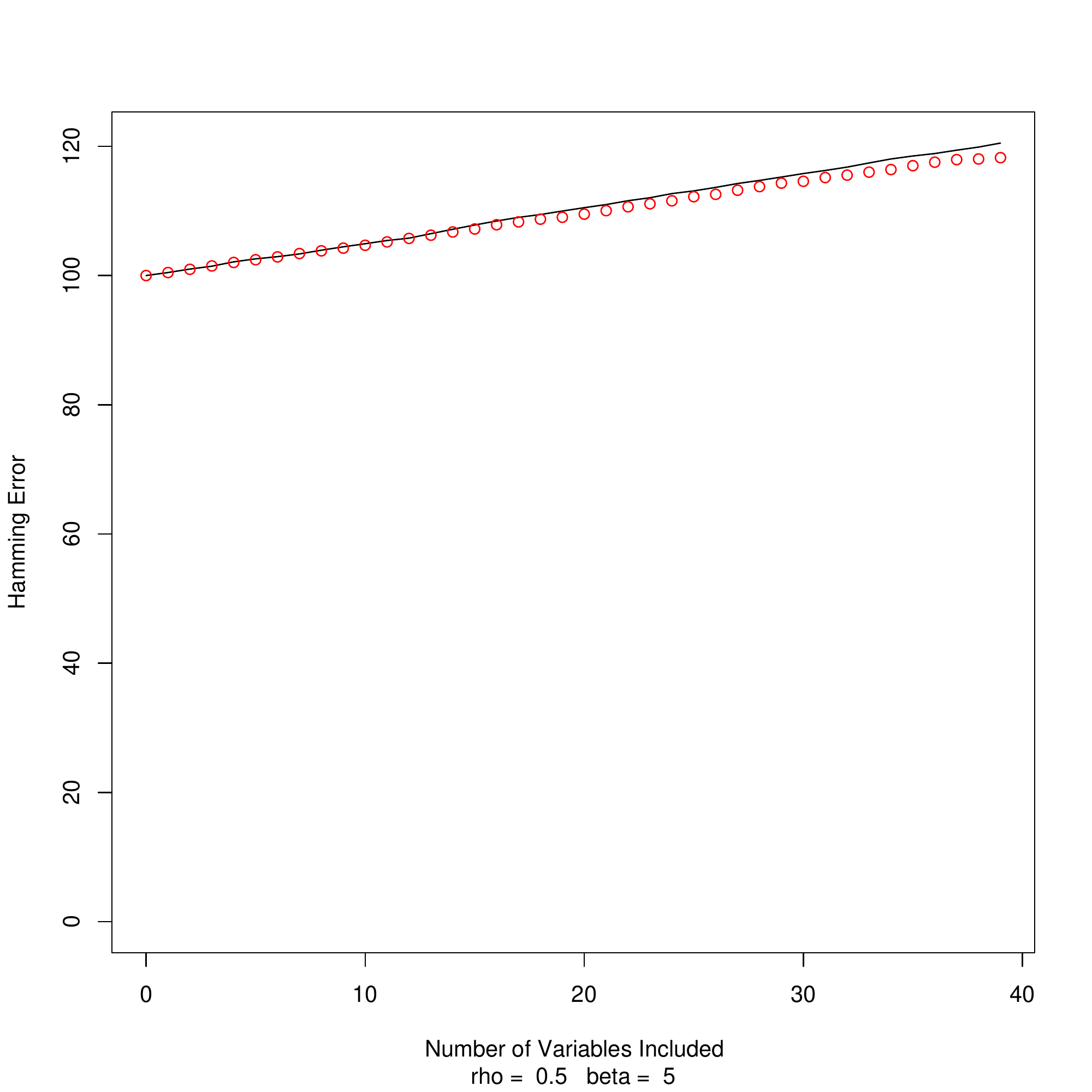}\\
\includegraphics[height=2in, width= 2.5in]{plot8_Part1.pdf}\includegraphics[height=2in, width= 2.5 in]{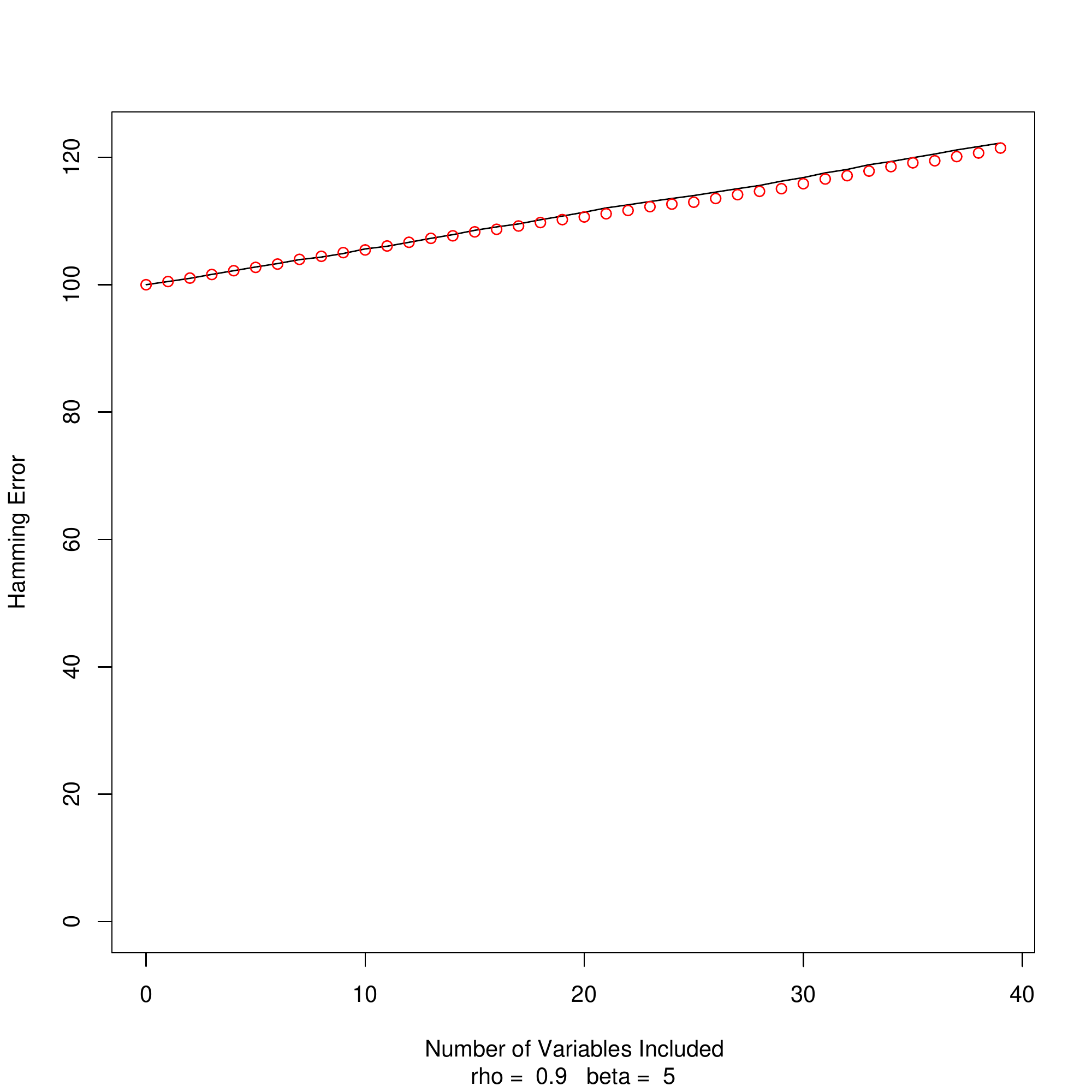}\\
\caption{Same as in Figure \ref{fig::ex2}, but displayed are the average errors across $100$ replications.}
\label{fig::sim2}
\end{figure}

In this section we consider some numerical examples.
Figures \ref{fig::ex1} and \ref{fig::ex2}
show the prediction error and the
Hamming error for the lasso and marginal regression,
as a function of the number of variables selected.
In all cases, $n=40$, $p=500$,
$\sigma = 10$ and $s =100$. All nonzero $\beta_j$'s
are set equal to either $.5$ or 5. Each row of the matrix $X$ is generated 
independently from $N(0, \Sigma(\rho))$, where $\Sigma(\rho)$ is a $p$ by $p$ matrix 
with $1$ on the diagonal and $\rho$ elsewhere.   We take  $\rho=0, 0.2, 0.5,0 .9$.
Figures \ref{fig::sim1} and \ref{fig::sim2}
are the same but are averaged over 100 replications.

We see that in virtually all cases,
marginal regression is competitive with the lasso
and in some cases is much better.


\section{Proofs}  
\label{sec:proof} 
\subsection{Proof of Theorem \ref{thm::marginal-Mrho}} 
First, let $k_i$ denote the number of non-zero diagonal entries in row $i$ of $C$.
Because $C$ is symmetric but not diagonal, at least two rows must have non-zero $k_i$.
Assume without loss of generality that the rows and columns of $C$ are arranged 
so that the rows with non-zero $k_i$ form the initial minor.
It follows that the initial minor is itself a positive definite symmetric matrix.
And because any such matrix $A$ satisfies $|A_{ij}| < \max_k C_{kk}$ for $j \ne i$,
there exists a row $i$ of $C$ with $k_i > 0$ and $|C_{ij}| < C_{ii}$ for any $j \ne i$.

Define $\beta$ as follows:
\begin{equation}
\beta_j = \begin{cases}
\frac{\rho C_{ii}}{C_{ij}} & \text{if $j \ne i$ and $C_{ij} \ne 0$} \\
\rho                      & \text{if $j \ne i$ and $C_{ij} = 0$} \\
-k_i \rho                 & \text{if $j = i$.}.
\end{cases}
\end{equation}
Because $|C_{ij}| \le C_{ii}$, this satisfies $|\beta_j| \ge \rho$, so $\beta \in \cM^s_\rho$.
Moreover,
\begin{equation}
(C\beta)_i = \sum_j C_{ij} \beta_j = -k_i C_{ii} \rho + \sum_{{j\ne i\atop C_{ij} \ne 0}}\frac{\rho C_{ii}}{C_{ij}} C_{ij} = 0.
\end{equation} 
This proves the theorem.

\subsection{Proof of Lemma \ref{lemma:F'}} 
By the definition of $\hat{S}_n(s)$, it is sufficient to show that 
except for a probability that tends to $0$,  
\[
\max |X_N^T Y| <  \min |X_S^TY|. 
\]
Since $Y = X \beta + z = X_S \beta_S + z$, we have $X_N^T  Y = X_N^T (X _S \beta_S  + z) = C_{NS} \beta_S  + X_N^T z$.  Note that $x_i^T z \sim N(0, \sigma_n^2)$.  By Boolean algebra and elementary statistics, 
\[
P ( \max | X_N^T  z | > \sigma_n \sqrt{2 \log p })    \leq   \sum_{ i \in N}  P(|(x_i, z)| \geq \sigma_n \sqrt{2 \log p})  \leq \frac{C}{\sqrt{\log p}} \frac{p-s}{p}. 
\]
It follows that except for a probability of $o(1)$,  
\[
\max| X_N^T Y| \leq \max|C_{NS} \beta_S|   + \max | X_N^T  z|  \leq  \max | C_{NS} \beta_S|  + \sigma_n \sqrt{2 \log p}. 
\]
Similarly, except for a probability of $o(1)$,  
\[
\min|X_S^T Y|  \geq \min|C_{SS} \beta_S|  - \max |X_S^T z| \geq  \min |C_{SS} \beta_S|  - \sigma_n  \sqrt{2 \log p}. 
\]
 Combining  these gives the claim.  
 \qed

\subsection{Proof of Theorem \ref{thm:noise}}   
Once  the first claim is proved, the second claim follows from Lemma \ref{lemma:F'}.  So we only show the first claim.  Write for short $\hat{S}_n(s) = \hat{S}_n(s^{(n)}; X^{(n)}, Y^{(n)}, p^{(n)})$, $s = s^{(n)}$, and $S = S(\beta^{(n)})$.  All we need to show is 
\[
\lim_{n \goto \infty} P( \hat{s}_n \neq  s)  = 0. 
\]
Introduce the event 
\[
D_n =  \{\hat{S}_n(s)   =  S \}.  
\]
It follows from Lemma \ref{lemma:F'} that 
\[
P(D_n^c)  \goto 0.  
\]
Write 
\[
P( \hat{s}_n \neq s) \leq  P(D_n) P(\hat{s}_n \neq s | D_n) + P(D_n^c).  
\]
It is sufficient to show $\lim_{n \goto \infty}  P(\hat{s}_n \neq s | D_n)   = 0$, or equivalently, 
\begin{equation} \label{noisepf0}
\lim_{n \goto \infty} P( \hat{s}_n  >  s | D_n)  = 0 \qquad \mbox{and}  \qquad \lim_{n \goto \infty} P( \hat{s}_n <   s | D_n)  = 0.   
\end{equation} 

Consider the first claim of (\ref{noisepf0}).  Write for short $t_n =  \sigma_n \sqrt{2 \log n}$. 
Note that the event $\{\hat{s}_n > s | D_n\}$ is contained in the event of $\cup_{k = s}^{p-1} \{ \hat{\delta}_n(k)  \geq t_n | D_n\}$.   Recalling  $P(D_n^c) = o(1)$, 
\begin{equation} \label{noisepf1}
P( \hat{s}_n > s) \leq \sum_{k = s}^{p-1} (\hat{\delta}(k) \geq  t_n |D_n)  \lesssim \sum_{k=s}^{p-1}  P(\hat{\delta}_n(k) \geq  t_n), 
\end{equation} 
where we say two positive sequences $a_n \lesssim b_n$ if $\limsup_{n \goto \infty} (a_n/b_n)  \leq 1$.

Fix $s \leq k \leq p -1$. By definitions,   $\hat{H}(k+1) - \hat{H}(k)$ is the projection matrix from $R^n$ to $\hat{V}_n(k+1) \cap \hat{V}_n(k)^{\perp}$.  So conditional on the event $\{ \hat{V}_n(k+1)  = \hat{V}_n(k)\}$, 
$\delta_n(k) = 0$, and conditional on the event $\{ \hat{V}_n(k+1)  \subsetneq  \hat{V}_n(k)\}$,  $\delta_n^2(k)  \sim \sigma_n^2 \chi^2(1)$.    Note  that $P(\chi^2(1) \geq \ 2 \log n) = o(1/n)$. It follows that 
\begin{align}   
\sum_{k = s}^{p-1}   P( \hat{\delta}_n(k) \geq t_n)  &=   \sum_{k = s}^{p-1}  P(\hat{\delta}_n(k) \geq t_n |  \hat{V}_n(k) \subsetneq \hat{V}_n(k+1)) P(\hat{V}_n(k) \subsetneq \hat{V}_n(k+1))   \nonumber
\\ 
&= o(\frac{1}{n})  \sum_{k = s}^{p-1}   P(\hat{V}_n(k) \subsetneq \hat{V}_n(k+1)).  \label{noisepf2} 
\end{align} 
Moreover, 
\begin{align*} 
\sum_{k =s}^{p-1} P(\hat{V}_n(k) \subsetneq  \hat{V}_n(k+1))  &=  \sum_{k =s}^{p-1} E[1( \mathrm{dim}(\hat{V}_n(k+1))  >  \mathrm{dim}(\hat{V}_n(k)))]  \\
  &= E[ \sum_{k = s}^{p-1} 1( \mathrm{dim}(\hat{V}_n(k+1))  >  \mathrm{dim}(\hat{V}_n(k)))]. 
\end{align*} 
Note that for any realization of the sequences $\hat{V}_n(1), \ldots, \hat{V}_n(p)$,   
$\sum_{k = s}^{p-1}   1(\mathrm{dim}(\hat{V}_n(k+1)) >  \mathrm{dim}(\hat{V}_n(k)))  \leq n$. 
It follows that 
\begin{equation} \label{noisepf3}
\sum_{k = s}^{p-1} P(\hat{V}_n(k) \subsetneq  \hat{V}_n(k+1)) \leq n. 
\end{equation} 
Combining  (\ref{noisepf1})-(\ref{noisepf3}) gives  the claim. 

Consider the second claim of (\ref{noisepf0}). 
By the definition of $\hat{s}_n$,     the event 
$\{ \hat{s}_n < s | D_n)\}$   
is contained in the event 
$\{ \hat{\delta}_n(s-1)    <  t_n   | D_n\}$. By definitions,   $\hat{\delta}_n(s-1) = \|  (\hat{H}(s)  - \hat{H}(s-1)) Y  \|$, where $\|\cdot\| = \|\cdot\|_2$ denotes the $\ell^2$ norm.  So all we need to show is     
\begin{equation} \label{noisepf7}
\lim_{n \goto \infty} P(  \|  (\hat{H}(s) - \hat{H}(s-1)) Y  \|  <   t_n | D_n)  = 0. 
\end{equation}

Fix $1 \leq k \leq p$. Recall that $i_k$ denotes the index at  which the rank of $|(Y, x_{i_k})|$ among all $|(Y, x_j)|$ is $k$.  Denote  $\tilde{X}(k)$ by the $n$ by $k$ matrix $[x_{i_1}, x_{i_2}, \ldots, x_{i_k}]$, and denote $\tilde{\beta}(k)$ by the $k$-vector $(\beta_{i_1}, \beta_{i_2}, \ldots, \beta_{i_k})^T$.    Conditional on the event $D_n$,    $\hat{S}_n(s) = S$, and $\beta_{i_1}, \beta_{i_2},  \ldots, \beta_{i_s}$ are all the nonzero coordinates of $\beta$.  So according to our notations, 
\begin{equation} \label{noisepf6.9}
X \beta =  \tilde{X}(s)  \tilde{\beta}(s) =  \tilde{X}(s-1)  \tilde{\beta}(s-1)+  \beta_{i_s} x_{i_s} 
\end{equation} 
Now, first, note that $\hat{H}(s)  \tilde{X}(s)     =  \tilde{X}(s)$ and $\hat{H}(s-1)  \tilde{X}(s-1) =   \tilde{X}(s-1)$. Combine this with (\ref{noisepf6.9}). It follows from direct calculations that  
\begin{equation} \label{noisepf7.0}
(\hat{H}(s)  - \hat{H}(s-1)) X \beta  =   (I - \hat{H}(s-1)) x_{i_s}.
\end{equation} 
Second,  since 
$x_{i_s} \in \hat{V}_n(s)$,   $(I - \hat{H}(s)) x_{i_s} = 0$. So 
\begin{equation} \label{noisepf7.1} 
(I - \hat{H}_{s-1}) x_{i_s} = (I  - \hat{H}(s)  ) x_{i_s} + (\hat{H}(s)   - \hat{H}(s-1)) x_{i_s}  =  (\hat{H}_s - \hat{H}_{s-1}) x_{i_s}. 
\end{equation} 
Last, split $x_{i_s}$ into two terms, 
$x_{i_s} = x_{i_s}^{(1)} + x_{i_s}^{(2)}$ 
such that $x_{i_s}^{(1)} \in \hat{V}_n(s-1)$ and $x_{i_s}^{(2)}  \in \hat{V}_n(s)  \cap (\hat{V}_n(s-1))^{\perp}$.   It follows that  $(\hat{H}(s)  - \hat{H}(s-1)) x_{i_s}^{(1)} = 0$, and  so  
\begin{equation} \label{noisepf7.2}
(\hat{H}(s)  - \hat{H}(s-1)) x_{i_s}  =   (\hat{H}(s)   - \hat{H}(s-1)) x_{i_s}^{(2)}. 
\end{equation} 
Combining (\ref{noisepf7.0})-(\ref{noisepf7.2}) gives 
\begin{equation}  \label{noisepf7.3}
(\hat{H}(s)  - \hat{H}(s-1)) X \beta  = (\hat{H}(s)  - \hat{H}(s-1)) x_{i_s}^{(2)}.  
\end{equation} 
Recall that $Y  = X \beta + z$, it follows that 
\begin{equation} \label{noisepf7.4} 
(\hat{H}_s - \hat{H}_{s-1}) Y    = (\hat{H}(s)  - \hat{H}(s-1)) ( \beta_{i_s} x_{i_s}^{(2)}    + z). 
\end{equation} 

Now,  take an orthonormal basis of $R^n$, say $\hat{q}_1, \hat{q}_2, \ldots, \hat{q}_n$, such that $\hat{q}_1 \in  \hat{V}_n(s) \cap \hat{V}_n(s-1)^{\perp}$,  $\hat{q}_2, \ldots, \hat{q}_s \in \hat{V}_n(s-1)$, and $\hat{q}_{s+1}, \ldots, \hat{q}_n \in \hat{V}_n(s)^{\perp}$.   Recall  that $x_{i_s}^{(2)}$ is contained in the one dimensional linear space $\hat{V}_n(s) \cap \hat{V}_n(s-1)^{\perp}$,  so without loss of generality, assume 
$(x_{i_s}^{(2)}, \hat{q}_1)  = \| x_{i_s}^{(2)}\|$.   Denote the square matrix $[\hat{q}_1, \ldots, \hat{q}_n]$ by $\hat{Q}$.  Let  $\tilde{z} = \hat{Q}z$ and let  $\tilde{z}_1$ be the first coordinate of $\tilde{z}$.      Note that marginally $\tilde{z}_1 \sim N(0, \sigma_n^2)$.
Over the event $D_n$,  it follows from the construction of $\hat{Q}$ and basic algebra that 
\begin{equation} \label{noisepf7.5}
\|(\hat{H}(s) - \hat{H}(s-1)) (\beta_{i_s} x_{i_s}^{(2)} + z) \|^2 =   ( \|\beta_{i_s} x_{i_s}^{(2)}\| + \tilde{z}_1)^2. 
\end{equation} 
Combine (\ref{noisepf7.4}) and (\ref{noisepf7.5}),   
\[
\|(\hat{H}(s)  - \hat{H}(s-1)) Y\|^2    =  ( \|\beta_{i_s} x_{i_s}^{(2)}\| + \tilde{z}_1)^2, \qquad \mbox{over the event $D_n$}. 
\]
As a result,  
\begin{equation} \label{noisepf7.6}
P(\|(\hat{H}(s)  - \hat{H}(s-1)) Y\| < t_n |D_n)    =   P(  (\|\beta_{i_s} x_{i_s}^{(2)}\| + \tilde{z}_1)^2  < t_n  |D_n). 
\end{equation}

Recall that conditional on the event $D_n$,  $\hat{S}_n(s) = S$.  So by the definition of $\Delta_n^* = \Delta_n(\beta, X, p)$,  
\[
\| \beta_{i_s} x_{i_s}^{(2)}\|  \geq \Delta_n^*, 
\]
and 
\begin{equation} \label{noisepf7.7}
P(  (\|\beta_{i_s} x_{i_s}^{(2)}\| + \tilde{z}_1)^2  < t_n  |D_n)  \leq  P(  \|\beta_{i_s} x_{i_s}^{(2)}\| + \tilde{z}_1  < t_n  |D_n) \leq  P ( \Delta_n^* + \tilde{z}_1  < t_n  |D_n).  
\end{equation} 
Recalling that $\tilde{z}_1 \sim N(0, \sigma_n^2)$ and that $P(D_n^c) = o(1)$, 
\begin{equation}  \label{noisepf7.8} 
P( \Delta_n^* + \tilde{z}_1 < t_n | D_n) \leq P(\Delta_n^* + \tilde{z}_1 < t_n) + o(1).  
\end{equation} 
Note that by  the assumption of $(\frac{\Delta_n^*}{\sigma_n} - t_n) \goto \infty$,  $   P(\Delta_n^* + \tilde{z}_1 < t_n) = o(1)$. 
Combining this with  (\ref{noisepf7.7})-(\ref{noisepf7.8}) gives
\begin{equation} \label{noisepf7.9} 
P(  (\|\beta_{i_s} x_{i_s}^{(2)}\| + \tilde{z}_1)^2  < t_n^2   |D_n)  = o(1). 
\end{equation} 
Inserting (\ref{noisepf7.9}) into (\ref{noisepf7.6}) gives (\ref{noisepf7}).  
\qed

\subsection{Proof of Lemma \ref{lemma:coherent}} 
For $1 \leq i \leq p$,  introduce the random variable 
\[
Z_i =  \sum_{j \neq  i}^p  \beta_j  (x_i, x_j). 
\]
When $B_i = 0$,  $\beta_i = 0$, and so  $Z_i = \sum_{j =1}^p \beta_j (x_i,x_j)$. 
By the definition of $C_{NS}$,  
\[
\max |C_{NS} \beta_S|   = \max_{1 \leq i \leq p}  \{ (1 - B_i) \cdot  |\sum_{j = 1}^p \beta_j (x_i, x_j)|\}  = \max_{1 \leq i \leq p}  \{(1 - B_i) |Z_i|\}.     
\]
Also, recalling  that the columns of matrix $X$  are  normalized such that $(x_i,x_i)  =1$,     
the diagonal coordinates of $(C_{SS} - I)$ are $0$. Therefore, 
\[
\max |(C_{SS} - I) \beta_S| = \max_{1 \leq i \leq p}  \{  B_i \cdot  |\sum_{j \neq i }  \beta_j (x_i, x_j)|\}    =  \max_{1 \leq i \leq p}  \{ B_i \cdot  |Z_i| \}. 
\]
Note that $Z_i$ and $B_i$ are independent and  that $P(B_i = 0) = (1 - \eps)$. It follows that   
\[
P( \max|C_{NS} \beta_S| \geq \delta)  \leq \sum_{i =1}^p P(B_i = 0) P(|Z_i| \geq \delta | B_i = 0) = (1 - \eps)  \sum_{i = 1}^p  P(|Z_i| \geq \delta), 
\]
and 
\[
P(\max|(C_{SS} - I) \beta_S| \geq \delta) \leq   \sum_{i =1}^p P(B_i = 1) P(|Z_i| \geq \delta | B_i = 1) =  \eps  \sum_{i = 1}^p  P(|Z_i| \geq \delta). 
\]
Compare these with the lemma. It is sufficient to show 
\begin{equation}  \label{cohpf0.0} 
P(|Z_i| \geq \delta)  \leq  e^{-\delta t} [e^{\eps\bar{g}_i(t)} + e^{\eps \bar{g}_i(-t)}]. 
\end{equation}

Now, by the definition of $g_{ij}(t)$,   the moment generating function of $Z_i$ satisfies  that 
\begin{equation} \label{cohpf0}
E[e^{t Z_i}]   = E[e^{t \sum_{j \neq  i} \beta_j (x_i, x_j)}] =   \Pi_{j \neq i}  [1 + \eps g_{ij}(t)].
\end{equation} 
Since $1 + x \leq e^x$ for all $x$,   $1  + \eps g_{ij}(t)    \leq e^{\eps  g_{ij}(t)}$, so   by the definition of $\bar{g}_i(t)$,  
\begin{equation} \label{cohpf3.1}
E[e^{t Z_i}]  \leq   \Pi_{j \neq i} e^{\eps  g_{ij}(t)}  = e^{ \eps \bar{g}_i(t) }.   
\end{equation} 
It follows from 
 Chebyshev's inequality that  
\begin{equation} \label{cohpf4.1}
P (Z_i  \geq \delta) \leq     e^{-\delta t} E[e^{t Z_i}] \leq e^{-\delta t} e^{\eps \bar{g}_i(t)}.  
\end{equation}
Similarly, 
\begin{equation} \label{cohpf4.2}
P(Z_i < -\delta) \leq e^{-\delta t} e^{\eps \bar{g}_i( - t)}
\end{equation} 
Inserting (\ref{cohpf4.1})-(\ref{cohpf4.2}) into (\ref{cohpf0.0}) gives the claim.   \qed 

\subsection{Proof of Corollary \ref{lemma:weakdependent}} 
Choose a constant $q$ such that $q/2 - c_2q > 1$ and let $t_n =  q \log (p) / a_n$.   By the definition of $A_n(a_n/2, \eps_n, \bar{g})$,   it is sufficient to show that for all $1 \leq i \leq p$,   
\[
e^{- a_n t_n/2} e^{\eps_n \bar{g}_i(t_n)}  = o(1/p), \qquad e^{- a_n t_n/2} e^{\eps_n \bar{g}_i(- t_n)} = o(1/p).
\]
The proofs are similar, so we only show the first one.  Let $u$ be a random variable such that $u  \sim \pi_n$. Recall that  the support of $|u|$ is contained in $[a_n, b_n]$.  By the assumptions and the choice of $t_n$, for all fixed $i$ and $j \neq i$,   
$|t_n u (x_i, x_j)| \leq q \log (p) (b_n/a_n) |(x_i, x_j)| \leq  c_1 q$.    Since $e^{x} - 1 \leq x + e^{x} x^2/2$,   it follows from Taylor expansion that 
\[
 \eps_n \bar{g}_i(t_n)   = \eps_n [e^{t_n u (x_i, x_j)} -1]  \leq  \eps_n \sum_{j \neq i}  E_{\pi_n} [t_n u (x_i, x_j) + \frac{e^{c_1q}}{2}  t_n^2 u^2 (x_i, x_j)^2]. 
\]
By definitions of $m_n(X)$ and $v_n^2(X)$,   $\eps_n \sum_{j \neq i} E_{\pi} [t_n u (x_i, x_j)] = t_n \mu_n^{(1)}  m_n(X)$, and $\eps_n \sum_{j \neq i} 
E_{\pi_n} [t_n^2 u^2 (x_i, x_j)^2]  = t_n^2 \mu_n^{(2)} v_n^2(X)$.   It follows from (\ref{weakcondition}) that  
\[
\eps_n \bar{g}_i(t_n) \leq  q \log(p)     \cdot [\frac{\mu_n^{(1)}}{a_n} m_n(X)   + \frac{e^{c_1q}}{2} \frac{\mu_n^{(2)}}{a_n^2}  v_n^2(X)  q  \log(p)]  \lesssim  q c_2 \log(p). 
\]
Therefore,  
\[
e^{-a_n t_n/2} e^{\eps_n \bar{g}_i(t_n)}  \leq e^{- [ q/2   - c_2q  + o(1)] \log(p)},  
\]
and claim follows by the choice of $q$.  \qed 

\subsection{Proof of Corollary \ref{lemma:sparsespikes}}   
Choose a constant $q$ such that $2 < q  < \frac{c_3}{c_4 \delta}$.   Let $t_n  =  a_n q \log(p)$, and $u$ be a random variable such that $u \sim \Pi_n$.  Similar to the proof of Lemma \ref{lemma:weakdependent}, we only show that 
\[
e^{-a_n t_n/2} e^{\eps_n \bar{g}_i(t_n)} = o(1/p),  \qquad \mbox{for all $1 \leq i \leq p$}. 
\]
Fix $i \neq j$. When  $(x_i, x_j)  = 0$,  $e^{t u (x_i, x_j)} - 1 = 0$.   When $(x_i,x_j) \neq 0$,  
$e^{t_n u (x_i, x_j) } -1  \leq e^{t_n   (b_n/a_n)\delta} \leq e^{c_4 q \delta \log p}$. Also,  $\eps_n N_n^* \leq e^{- [c_3 + o(1)] \log(p)}$.    Therefore,  
\[
 \eps_n \bar{g}_i(t) \leq \eps_n N_n^* e^{c_4 q \delta \log(p)} \leq  e^{ -[c_3 - c_4 q \delta + o(1)] \log p}.   
\]
By the choice of $q$,   $c_3 - c_4 q \delta > 0$, so   $\eps_n \bar{g}_i(t) = o(1)$.  
It  follows that 
\[
e^{-a_n t_n/2} e^{\eps_n \bar{g}_i(t_n)} \leq o( e^{-a_n t_n/2})  =  o(e^{-q \log(p)/2}), 
\]
which gives the claim by $q > 2$.   \qed

\subsection{Proof of Theorem \ref{thm:LB}} 
Write 
\[
X = [x_1, \tilde{X}],  \qquad \beta = (\beta_1, \tilde{\beta})^T. 
\]
Fix a constant $c_0 > 3$.   Introduce the event  
\begin{equation} \label{DefineA}
D_n(c_0)   = \{ 1_S^T   \tilde{X}_S^T  \tilde{X}_S  1_S  \leq |S|  [1 +  \sqrt{\frac{|S|}{n}}  (1 +    \sqrt{2 c_0 \log p})]^2, \mbox{for all $S$}  \}. 
\end{equation} 
The following lemma is proved in Section \ref{subsec:A}. 
\begin{lemma}  \label{lemma:A}
Fix  $c_0 > 3$. As  $p \goto \infty$, 
\[
P(D_n^c(c_0))   = o(1/p^2). 
\]
\end{lemma} 
Since $d_n(\hat{\beta} |X) \leq p$ for any variable selection procedure $\hat{\beta}$,   Lemma \ref{lemma:A} implies that  the overall contribution of $D_n^c$ to the Hamming distance $d_n^*(\hat{\beta})$  is $o(1/p)$.    In addition, write 
\[
d_n(\hat{\beta}  | X)   =   \sum_{j =1}^p  E[ 1(\hat{\beta}_j  \neq \beta_j)].
\]
 By symmetry,  it is sufficient to show that for any realization of $(X, \beta) \in D_n(c_0)$,   
\begin{equation} \label{LBtoshow0}
E[1(\hat{\beta}_j  \neq \beta_j)]  \geq   \left\{ \begin{array}{ll}
L(n) p^{ -\frac{(\vartheta + r)^2}{4 r}},  &\qquad r \geq  \vartheta,  \\
p^{  -\vartheta},  &\qquad  0 < r < \vartheta, 
\end{array}  
\right. 
\end{equation} 
where $L(n)$  is a multi-log term that does not depend on $(X, \beta)$.

We now show (\ref{LBtoshow0}).  Toward this end, we relate the estimation problem  to the problem of  testing the null hypothesis of $\beta_1 = 0$   versus the alternative hypothesis of $\beta_1 \neq 0$. 
Denote $\phi$ by the density of $N(0,1)$. Recall that 
$X = [x_1, \tilde{X}]$ and $\beta = (\beta_1, \tilde{\beta})^T$.  
The joint density associated with  the  null hypothesis   is 
\[
f_0(y) =  f_0(y; \eps_n, \tau_n, n | X)  \phi( y  - \tilde{X} \tilde{\beta} )  d \tilde{\beta} = \phi(y) \int e^{y^T \tilde{X} \tilde{\beta} - | \tilde{X} \tilde{\beta}|^2/2} d \tilde{\beta},   
\]
and the joint density associated with  the alternative hypothesis is 
\begin{align} 
f_1(y) = f_1(y; \eps_n, \tau_n, n |X)  &= \int \phi( y - \tau_n   x_1 - \tilde{X} \tilde{\beta} )  d \tilde{\beta}  \nonumber  \\
&= \phi(y - \tau_n  x_1) \int e^{y^T \tilde{X} \tilde{\beta} - | \tilde{X} \tilde{\beta}|^2/2}  e^{ - \tau_n  x_1^T \tilde{X} \tilde{\beta}} d \tilde{\beta}.      \label{Defineff1}
\end{align} 
Since the prior probability that the null hypothesis is true is $(1 - \eps_n)$,   the optimal test is  
the Neyman-Pearson test that rejects the null  if and only if 
\[
\frac{f_1(y)}{f_0(y)} \geq \frac{(1 - \eps_n)}{\eps_n}. 
\]
The optimal testing error is equal to 
\[
1 - \|(1 - \eps_n) f_0 -  \eps_n  f_1\|_1. 
\]
Compared  to (\ref{Definelasso}),  $\| \cdot \|_1$ stands for  the $L^1$-distance between two functions, not the $\ell^1$ norm of a vector.   

We need to modify $f_1$ into a more tractable form, but with negligible difference in $L^1$-distance.  
Toward this end, let $N_n(\tilde{\beta})$ be the number of nonzeros coordinates of   
$\tilde{\beta}$.  Introduce  the event  
\[
B_n  = \{ |N_n(\tilde{\beta}) - p \eps_n| \leq \frac{1}{2} p \eps_n \}.  
\]
Let  
\begin{equation} \label{Defineany}
a_n(y)  = a_n(y; \eps_n, \tau_n | X) =  \frac{ \int  (e^{y^T  \tilde{X} \tilde{\beta}  - | \tilde{X} \tilde{\beta}|^2/2})   (e^{-\tau_n x_1^T  \tilde{X} \tilde{\beta}})      \cdot 1_{\{B\}}  d \tilde{\beta}}{ \int  (e^{-y^T  \tilde{X} \tilde{\beta}  - | \tilde{X} \tilde{\beta}|^2/2})     \cdot 1_{\{B\}}  d \tilde{\beta}}.  
\end{equation}
Note that the only difference between the numerator and the denominator is the term $ e^{-\tau_n x_1^T  \tilde{X} \tilde{\beta}}$ which $\approx 1$ with high probability. 
Introduce
\begin{equation} \label{Definefff1} 
\tilde{f}_1(y) =  a_n(y) \phi(y - \tau_n x_1)  \int  e^{y^T \tilde{X}  \tilde{\beta} - | \tilde{X}  \tilde{\beta}|^2/2}     d \tilde{\beta}. 
\end{equation} 
The following lemma is proved in  Section \ref{subsec:f1}. 
\begin{lemma}  \label{lemma:f1} 
As $p \goto \infty$, there is a generic constant $c > 0$ that does not depend on $y$ such that 
$|a_n(y) -1| \leq  c \log(p) p^{(1 - \vartheta) - \theta/2}$ and $  
\|f_1 - \tilde{f}_1\|_1 = o(1/p)$. 
\end{lemma} 

We now ready to show the claim. 
Define $\Omega_n = \{y:  a_n(y) \phi(y - \tau_n x_1) \geq \phi(y)\}$. 
Note that by the definitions of $f_0(y)$ and $\tilde{f}_1(y)$,  $y \in \Omega_n$ if and only if 
\[
\frac{\eps_n \tilde{f}_1(y)}{(1 - \eps_n) f_0(y)}   \geq 1.  
\]
By  Lemma \ref{lemma:f1}, 
\[
|\int \tilde{f}_1(y) dy - 1|     \leq  \| \tilde{f}_1  - f_1\|_1 \leq  o(1/p). 
\]
It follows from elementary calculus that 
\[
1 - \|(1 - \eps_n) f_0 - \eps_n \tilde{f}_1\|_1  =   \int_{\Omega_n} (1  - \eps_n)  f_0(y) dy  + \int_{\Omega_n^c} \eps_n \tilde{f}_1(y) dy + o(1/p). 
\]
Using Lemma \ref{lemma:f1} again, we can replace $\tilde{f}_1$ by $f_1$ on the right hand side, so 
\[
1 - \|(1 - \eps_n) f_0 - \eps_n \tilde{f}_1\|_1  =   \int_{\Omega_n} (1 - \eps_n)  f_0(y) dy  + \int_{\Omega_n^c} \eps_n f_1(y) dy + o(1/p).  
\]

At the same time,  let $\delta_p  = c \log(p)  p^{(1 - \vartheta)  - \theta/2}$ be as in Lemma \ref{lemma:f1},  
and let 
\[
t_0 = t_0(\vartheta, r) =  \frac{\vartheta + r}{2 \sqrt{r}} \sqrt{2 \log p}. 
\]
be the unique solution of the equation 
$\phi(t) =  \eps_n \phi(t - \tau_n)$.  It follows from   Lemma \ref{lemma:f1} that, 
\[
\{\tau_n x^T  y  \geq   t_0  (1 + \delta_p)\}  \subset   \Omega_n \subset   \{ \tau_n x_1^T y  \geq t_0 (1  - \delta_p)\}. 
\]
As a result, 
\[
 \int_{\Omega_n} f_0(y) dy  \geq   \int_{\tau_n x_1^T  y \geq  t_0(1  +    \delta_p)} f_0(y) 
\equiv  P_0( \tau_n x_1^T  Y \geq t_0 (1 + \delta_p)), 
\]
and 
\[
 \int_{\Omega_n^c} f_1(y) dy  \geq   \int_{\tau_n x_1^T y \leq  t_0(1  -     \delta_p)} f_1(y)   \equiv  P_1( \tau_n x_1^T  Y \leq t_0 (1 - \delta_p)). 
\]
Note that under the null,  
$x_1^T  Y  = x_1^T  \tilde{X} \tilde{\beta} + x_1^T  z$.  
It is seen that given $x_1$,  $x_1^T z \sim N(0, |x_1|^2)$,   and $|x_1|^2 = 1 + O(1/\sqrt{n})$. Also, it is seen that  except for a probability of $o(1/p)$,  $x_1^T  \tilde{X} \tilde{\beta}$ is algebraically small. 
It follows that 
\[
P_0( \tau_n x_1^T Y \geq t_0 (1 + \delta_p))  \lesssim \bphi(t_0) =  L(n)  p^{ -\frac{(\vartheta + r)^2}{4r}}, 
\] 
where $\bar{\Phi} = 1 - \Phi$ is the survival function of $N(0,1)$. 
Similarly,  under the alternative, 
\[
x_1^T  y   =  \tau_n (x_1, x_1) + x_1^T  \tilde{X} \tilde{\beta} + x_1^T z,  
\]
where $(x_1, x_1) = 1 + O(1/\sqrt{n})$. So   
\[
\eps_n  P_1( \tau_n x_1^T  y \leq t_0(1 - \delta_p)) \lesssim  \Phi(t_0 - \tau_n) = \left\{ 
\begin{array}{ll}
L(n) p^{-\frac{(\vartheta + r)^2}{4 r}},  &\qquad r \geq  \vartheta,  \\
L(n) p^{-\vartheta},  &\qquad  0 < r < \vartheta, 
\end{array}  
\right. 
\] 
Combine these gives the theorem.  \qed 

\subsubsection{Proof of Lemma \ref{lemma:A}}           \label{subsec:A}
It is seen that 
\[
P(D_n^c(c_0)) \leq \sum_{k = 1}^p  P\biggl(1_S^T X^T X 1_S \geq  k  [1 + \sqrt{\frac{k}{n}}  (1 +    \sqrt{2 c_0 \log p})]^2,  \mbox{for all $S$ with $|S| = k$}\biggr). 
\]
Fix $k \geq 1$. There are $\binom{p}{k}$ different $S$ with $|S| = k$.    It follows from \cite[Lecture 9]{Vershynin} that except a probability of 
$2\exp( - c_0  \log(p) \cdot k)$
that the largest eigenvalue of $X_S^T  X_S$ is no greater than 
$[1 + \sqrt{\frac{k}{n}} (1  + \sqrt{2 c_0 \log p})]^2$.   So for any $S$ with $|S|  = k$,  it follows from basic algebra that 
\[
P(1_S^T X^T X 1_S \geq  k  [1 + \sqrt{\frac{k}{n}}  (1 +    \sqrt{2 c_0 \log p})]^2)  \leq  2\exp( - c_0  \log(p) \cdot k). 
\]
Combining these with $\binom{p}{k} \leq p^k$ gives 
\[
P(D_n^c(c_0)) \leq 2 \sum_{k =1}^p \binom{p}{k} \exp(- c_0 (\log p) k)   \leq 2 \sum_{k = 1}^p \exp( - (c_0 -1) \log(p) k). 
\]
The claim follows by $c_0 > 3$. \qed

\subsubsection{Proof of Lemma \ref{lemma:f1}}   \label{subsec:f1}
 First,  we claim that 
for any $X$ in event $D_n(c_0)$, 
\begin{equation} \label{lemmaapf}
|x_1^T  \tilde{X} \tilde{\beta}| \leq  c \log(p)  (N(\tilde{\beta})/\sqrt{n}), 
\end{equation} 
where $c > 0$ is a generic constant.  
Suppose $N_n(\tilde{\beta}) = k$ and the nonzero coordinates of $\tilde{\beta}$ are $i_1, i_2, \ldots, i_k$.  Denote the $(k+1) \times (k+1)$ submatrix of $X^T X$ containing the $1^{st}$, $(1 +  i_1)$-th, $\ldots$, and $(1+i_k)$-th  rows and columns  by $U_{k+1}$. Let $\xi_1$ be the $(k+1)$-vector with $1$ on the first coordinate and $0$ elsewhere, let $\xi_2$ be the $(k+1)$-vector with $0$ on the first coordinate and $1$ elsewhere. Then 
\[
x_1^T  \tilde{X} \tilde{\beta} =   \tau_n \xi_1^T   U_{k+1} 
\xi_2 \equiv  \tau_n  \xi_1^T  (U_{k + 1} - I_{k+1}) \xi_2. 
\]
Let  $(U_{k+1} - I_{k+1})= Q_{k +1 } \Lambda_{k + 1} Q_{k+1}^T$ be the orthogonal decomposition.  By the definition of $D_n(c_0)$,     all eigenvalues of $(U_{k+1} - I_{k+1})$   are no greater than 
$(1 +  \sqrt{c \log(p) k/n})^2 - 1 \leq \sqrt{c  \log p} \sqrt{k/n}$ in absolute value.  As a result,    all diagonal coordinates of $\Lambda_{k+1}$ are no greater than 
\[
\sqrt{c  \log p} \sqrt{k/n} 
\]
in absolute value, 
and 
\[
\|\xi_1^T (U_{k+1} - I_{k+1}) \xi_2\| \leq   \|\xi_1^T  Q_{k+1} \Lambda_{k+1} \| \cdot  \|Q_{k+1}  \xi_2\| \leq  \sqrt{c  \log p} \sqrt{k/n} \|\xi_1^T  Q_{k+1} \| \cdot  \|Q_{k+1}  \xi_2\|. 
\] 
The claim follows from $\|\xi_1^T  Q_{k+1} \| = 1$ and $\| Q_{k+1}  \xi_2\| = \sqrt{k}$.

We now show the lemma. Consider the first claim.  
Consider a realization of $X$ in the event  $D_n(c_0)$ and a realization of $\tilde{\beta}$ in the event  $B_n$. 
 By the definitions of $B_n$,    $N_n(\tilde{\beta})   \leq p \eps_n +  \frac{1}{2} p \eps_n$. Recall that $p \eps_n = p^{1 - \vartheta}$, $n  = p^{\theta}$.  It follows  that  $\log(p) N(\tilde{\beta}) /\sqrt{n}  \leq c \log(p) p \eps_n/\sqrt{n} = c \log(p) p^{1 - \vartheta - \theta/2}$.  Note that by the assumption of $(1 - \vartheta) < \theta/2$, the exponent is negative.  Combine this with 
 (\ref{lemmaapf}),  
\begin{equation} \label{lemmaaa}
| e^{- \tau_n x_1^T \tilde{X} \tilde{\beta}}  - 1|  \leq c \log(p) (N (\tilde{\beta}) /\sqrt{n}), 
\end{equation} 
Now, note that in the definition of $a_n(y)$  (i.e. (\ref{Defineany})),   the only difference between the integrand on the top  and that  on the bottom is the term $e^{- \tau_n x_1^T \tilde{X} \tilde{\beta}}$. Combine this with    (\ref{lemmaaa}) gives the claim.

Consider the second claim.    By the definitions of $\tilde{f}_1(y)$ and $a_n(y)$,  
\begin{align*} 
\tilde{f}_1(y) &= a_n(y)  \phi(y - \tau_n x_1)   \cdot    
  \biggl[  \int [e^{y^T \tilde{X}  \tilde{\beta} - | \tilde{X}  \tilde{\beta}|^2/2}  1_{B_n}]  d \tilde{\beta}  +   \int [e^{y^T  \tilde{X}  \tilde{\beta} - | \tilde{X}  \tilde{\beta}|^2/2}     1_{B_n^c}]  d \tilde{\beta}\biggr] \\
&=  \phi(y - \tau_n x_1)   \cdot    
  \biggl[    \int [e^{y^T \tilde{X}  \tilde{\beta} - | \tilde{X}  \tilde{\beta}|^2/2}  e^{-\tau_n x_1^T \tilde{X} \tilde{\beta}} 1_{B_n^c}]  d \tilde{\beta}  +  a_n(y)  \int [e^{y^T  \tilde{X}  \tilde{\beta} - | \tilde{X}  \tilde{\beta}|^2/2}     1_{B_n^c}]  d \tilde{\beta}\biggr]. 
\end{align*} 
By the definition of $f_1(y)$, 
\[
f_1(y)  =   \phi(y - \tau_n x_1)   \cdot    
  \biggl[  \int [e^{y^T \tilde{X}  \tilde{\beta} - | \tilde{X}  \tilde{\beta}|^2/2}  e^{-\tau_n x_1^T \tilde{X}\tilde{\beta}} 1_{B_n}]  d \tilde{\beta}  +   \int [e^{y^T  \tilde{X}  \tilde{\beta} - | \tilde{X}  \tilde{\beta}|^2/2}  e^{-\tau_n x_1^T \tilde{X}\tilde{\beta}}     1_{B_n^c}]  d \tilde{\beta}\biggr]. 
\]
Compare two equalities and recall that $a_n(y) \sim 1$ (Lemma \ref{lemma:A}), 
\begin{align} 
\|f_1 - \tilde{f}_1\|_1  &\lesssim  \int \phi(y - \tau_n x_1) [ \int(  e^{y^T  \tilde{X}  \tilde{\beta} - | \tilde{X}  \tilde{\beta}|^2/2}     +      e^{y^T \tilde{X}  \tilde{\beta} - | \tilde{X}  \tilde{\beta}|^2/2} e^{-\tau_n x_1^T \tilde{X}\tilde{\beta}})     1_{B_n^c}   d \tilde{\beta}] dy \nonumber \\
&=  \int \int  \phi(y - \tau_n x_1 - \tilde{X} \tilde{\beta})  [   e^{\tau_n x_1^T \tilde{X} \tilde{\beta}} + 1] 1_{B_n^c} d \tilde{\beta} dy.  \label{Wellner1}
\end{align} 
Integrating over $y$, the last term is equal to 
$\int  [1 +   e^{ \tau_n x_1^T  \tilde{X} \tilde{\beta}}] \cdot 1_{B_n^c} d \tilde{\beta}$.

At the same time,  by  (\ref{lemmaapf}) and the definition of $B_n^c$,   
\begin{equation} \label{Wellner2}
\int  [1 +   e^{ \tau_n x_1^T  \tilde{X} \tilde{\beta}}] \cdot 1_{B_n^c} d \tilde{\beta} \leq   \sum_{\{k:  |k - p \eps_n| \geq \frac{1}{2} p \eps_n \}} [1 +  e^{c\log(p)  k /\sqrt{n}}]   P( N(\tilde{\beta})   = k).  
\end{equation} 
Recall that $p \eps_n = p^{1 - \vartheta}$, $n  = p^{\theta}$, and $(1 - \vartheta) < \theta/2$.   Using Bennett's inequality for $P(N(\tilde{\beta}) = k)$ (e.g. \cite[Page 440]{Wellner}),  it follows from elementary calculus that 
\begin{equation} \label{Wellner3} 
 \sum_{\{k:  |k - p \eps_n| \geq \frac{1}{2} p \eps_n \}} [1 +  e^{c\log(p)  k /\sqrt{n}}]   P( N(\tilde{\beta})   = k)  = o(1/p). 
\end{equation} 
Combining (\ref{Wellner1})--(\ref{Wellner3})  gives the claim.   \qed

\medskip\noindent
{\bf Acknowledgement}:  We  would like to thank David Donoho and Robert Tibshirani for helpful discussion.  CG  was supported in part by 
 NSF grant DMS-0806009 and NIH grant R01NS047493,    JJ was supported in part by 
NSF CAREER award  DMS-0908613, and LW was supported in part by NSF grant DMS-0806009.      

\end{document}